\documentclass[envcountsame,envcountsect]{svmult}

\usepackage{mathptmx}       
\usepackage{helvet}         
\usepackage{courier}        

\usepackage{amsmath,amssymb}

\smartqed

\newcommand\reel{{\mathbb{R}}}
\newcommand\complex{{\mathbb{C}}}

\newcommand\relatif{{\mathbb{Z}}}

\newcommand\Abar{{\overline A}}
\newcommand\Wbar{{\overline W}}
\newcommand\Bbar{{\overline B}}
\newcommand\Bcheck{{\check B}}
\newcommand\Ybar{{\overline Y}}
\newcommand\esp{{\mathbb{E}}}
\newcommand\eps{\varepsilon}
\newcommand\Itilde{{\widetilde I}}
\newcommand\Pitilde{{\widetilde\Pi}}
\newcommand\tribuf{{\mathcal{F}}}
\newcommand\proba{{\mathbb{P}}}
\newcommand\hilbert{{\mathcal{H}}}
\newcommand\hilberthat{{\widehat\hilbert}}
\newcommand\hilbertbar{{\overline\hilbert}}
\newcommand\wiener{{\mathcal{W}}}

\newcommand\inform{{\mathcal I}}

\newcommand\Btilde{{\widetilde B}}
\newcommand\holder{{\mathbb{H}}}

\newcommand\Bhat{{\widehat B}}

\newcommand\Ihatp{{\widehat I}_+}
\newcommand\Ibarp{{\overline I}_+}

\author{Jean Picard}
\institute{J. Picard (1,2)\at(1) Clermont Universit\'e, Universit\'e Blaise Pascal, Laboratoire de Math\'ematiques, BP 10448,
\mbox{F-63000} CLERMONT-FERRAND \smallskip\\
(2) CNRS, UMR 6620, LM, F-63177 AUBI\`ERE\\
\email{\texttt{Jean.Picard@math.univ-bpclermont.fr}}}

\title*{Representation formulae for the fractional Brownian motion}

\begin{document}

\maketitle

\abstract{ We discuss the relationships between some classical representations of the fractional
Brownian motion, as a stochastic integral with respect to a standard Brownian motion, or as a
series of functions with independent Gaussian coefficients. The basic notions of fractional
calculus which are needed for the study are introduced. As an application, we also prove some
properties of the Cameron-Martin space of the fractional Brownian motion, and compare its law with
the law of some of its variants. Several of the results which are given here are not new; our aim
is to provide a unified treatment of some previous literature, and to give alternative proofs and
additional results; we also try to be as self-contained as possible.}

\section{Introduction}

Consider a fractional Brownian motion $(B_t^H;t\in\reel)$ with Hurst parameter $0<H<1$. These
processes appeared in 1940 in \cite{kolmogorov40}, and they generalise the case $H=1/2$ which is
the standard Brownian motion. A huge literature has been devoted to them since the late 60's. They
are often used to model systems involving Gaussian noise, but which are not correctly explained
with a standard Brownian motion. Our aim here is to give a few basic results about them, and in
particular to explain how all of them can be deduced from a standard Brownian motion.

The process $B^H$ is a centred Gaussian process which has stationary increments and is
$H$-self-similar; these two conditions can be written as
\begin{equation}\label{statiosim}
   B_{t+t_0}^H-B_{t_0}^H\simeq B_t^H,\qquad
   B_{\lambda t}^H\simeq\lambda^HB_t^H
\end{equation}
for $t_0\in\reel$ and $\lambda>0$, where the notation $Z_t^1\simeq Z_t^2$ means that the two
processes have the same finite dimensional distributions. We can deduce from \eqref{statiosim} that
$B_{-t}^H$ and $B_t^H$ have the same variance, that this variance is proportional to $|t|^{2H}$,
and that the covariance kernel of $B^H$ must be of the form
\begin{align}
   C(s,t)&=\esp\bigl[B_s^HB_t^H\bigr]
   =\frac12\esp\Bigl[(B_s^H)^2+(B_t^H)^2-(B_t^H-B_s^H)^2\Bigr]\notag\\
   &=\frac12\esp\Bigl[(B_s^H)^2+(B_t^H)^2-(B_{t-s}^H)^2\Bigr]\notag\\
   &=\frac\rho2\Bigl(|s|^{2H}+|t|^{2H}-|t-s|^{2H}\Bigr)\label{bhcov}
\end{align}
for a positive parameter $\rho=\esp[(B_1^H)^2]$ (we always assume that $\rho\ne0$). The process
$B^H$ has a continuous modification (we always choose this modification), and its law is
characterised by the two parameters $\rho$ and $H$; however, the important parameter is $H$, and
$\rho$ is easily modified by multiplying $B^H$ by a constant. In this article, it will be
convenient to suppose $\rho=\rho(H)$ given in \eqref{rhoh}; this choice corresponds to the
representation of $B^H$ given in \eqref{mandness}. We also consider the restriction of $B^H$ to
intervals of $\reel$ such as $\reel_+$, $\reel_-$ or $[0,1]$.

Notice that the fractional Brownian motion also exists for $H=1$ and satisfies $B_t^1=t\,B_1^1$;
this is however a very particular process which is excluded from our study (with our choice of
$\rho(H)$ we have $\rho(1)=\infty$).

The standard Brownian motion $W_t=B_t^{1/2}$ is the process corresponding to $H=1/2$ and
$\rho=\rho(1/2)=1$. It is often useful to represent $B^H$ for $0<H<1$ as a linear functional of
$W$; this means that one looks for a kernel $K^H(t,s)$ such that the Wiener-It\^{o} integral
\begin{equation}\label{bhw}
   B_t^H=\int K^H(t,s)dW_s
\end{equation}
is a $H$-fractional Brownian motion. More generally, considering the family $(B^H;\;0<H<1)$
defined by \eqref{bhw}, we would like to find $K^{J,H}$ so that
\begin{equation}\label{bhbj}
   B_t^H=\int K^{J,H}(t,s)dB_s^J.
\end{equation}
In this case however, we have to give a sense to the integral; the process $B^J$ is a Gaussian
process but is not a semimartingale for $J\ne1/2$, so we cannot consider It\^{o} integration. In order
to solve this issue, we approximate $B^J$ with smooth functions for which the Lebesgue-Stieltjes
integral can be defined, and then verify that we can pass to the limit in an adequate functional
space in which $B^J$ lives almost surely. Alternatively, it is also possible to use integration by
parts.

The case where $K^{J,H}$ is a Volterra kernel ($K^{J,H}(t,s)=0$ if $s>t$) is of particular
interest; in this case, the completed filtrations of $B^H$ and of the increments of $B^J$ satisfy
$\tribuf_t(B^H)\subset\tribuf_t(dB^J)$, with the notation
\begin{equation}\label{tribu}
   \tribuf_t(X)=\sigma\bigl(X_s;s\le t\bigr),\quad
   \tribuf_t(dX)=\sigma\bigl(X_s-X_u;u\le s\le t\bigr).
\end{equation}
Notice that when the time interval is $\reel_+$, then $\tribuf_t(dB^J)=\tribuf_t(B^J)$ (because
$B_0^J=0$), but this is false for $t<0$ when the time interval is $\reel$ or $\reel_-$. When
$\tribuf_t(B^H)=\tribuf_t(B^J)$, we say that the representation \eqref{bhbj} is canonical;
actually, we extend here a terminology, introduced by \cite{levy56} (see \cite{hida:hitsuda93}),
which classically describes representations with respect to processes with independent increments
(so here the representation \eqref{bhw}); such a canonical representation is in some sense unique.

Another purpose of this article is to compare $B^H$ with two other families of processes with
similar properties and which are easier to handle in some situations:
\begin{itemize}
\item The so-called Riemann-Liouville processes on $\reel_+$ (they are also sometimes called
    type II fractional Brownian motions, see \cite{marinu:robin99}), are deduced from the
    standard Brownian motion by applying Riemann-Liouville fractional operators, whereas, as we
    shall recall it, the genuine fractional Brownian motion requires a weighted fractional
    operator.

\item We shall also consider here some processes defined by means of a Fourier-Wiener series on
    a finite time interval; they are easy to handle in Fourier analysis, whereas the Fourier
    coefficients of the genuine fractional Brownian motion do not satisfy good independence
    properties.
\end{itemize}
We shall prove that the Cameron-Martin spaces of all these processes are equivalent, and we shall
compare their laws; more precisely, it is known from \cite{feldman58,hajek58,hida:hitsuda93} that
two Gaussian measures are either equivalent, or mutually singular, and we shall decide between
these two possibilities.

Let us now describe the contents of this article. Notations and definitions which are used
throughout the article are given in Section \ref{defnot}; we also give in this section a short
review of fractional calculus, in particular Riemann-Liouville operators and some of their
modifications which are important for our study; we introduce some functional spaces of H\"older
continuous functions; much more results can be found in \cite{samko:kil:mari93}. In Section
\ref{timeinvers}, we give some resuts concerning the time inversion ($t\mapsto1/t$) of Gaussian
self-similar processes.

We enter the main topic in Section \ref{represent}. Our first aim is to explore the relationship
between two classical representations of $B^H$ with respect to $W$, namely the representation
of \cite{mandel:ness68},
\begin{equation}\label{mandness}
    B_t^H=\frac1{\Gamma(H+1/2)}\int_\reel\Bigl((t-s)_+^{H-1/2}-(-s)_+^{H-1/2}\Bigr)dW_s
\end{equation}
on $\reel$ (with the notation $u_+^\lambda=u^\lambda1_{\{u>0\}}$), and the canonical representation
on $\reel_+$ obtained in \cite{molchan:golo69,molchan03}, see also
\cite{decreuse:ustu99,norros:val:vir99} (this is a representation of type \eqref{bhw} for a
Volterra kernel $K^H$, and such that $W$ and $B^H$ generate the same filtration). Let us explain
the idea by means of which this relationship can be obtained; in the canonical representation on
$\reel_+$, we want $B_t^H$ to depend on past values $W_s$, $s\le t$, or equivalently, we want the
infinitesimal increment $dB_t^H$ to depend on past increments $dW_s$, $s\le t$. In
\eqref{mandness}, values of $B_t^H$ for $t\ge0$ involve values of $W_s$ for all $-\infty\le s\le
t$, so this is not convenient for a study on $\reel_+$. However, we can reverse the time
($t\mapsto-t$) and use the backward representation
\begin{equation*}
    B_t^H=\frac1{\Gamma(H+1/2)}\int_0^{+\infty}\Bigl(s^{H-1/2}-(s-t)_+^{H-1/2}\Bigr)dW_s
\end{equation*}
on $\reel_+$. Now the value of $B_t^H$ involves the whole path of $W$ on $\reel_+$, but we
can notice that the infinitesimal increment $dB_t^H$ only involves future increments $dW_s$,
$s\ge t$. Thus $dB^H(1/t)$ depends on past increments $dW(1/s)$, $s\le t$. We can then
conclude by applying the invariance of fractional Brownian motions by time inversion which has
been proved in Section \ref{timeinvers}. This argument is justified in \cite{molchan03} by using
the generalised processes $dB_t^H/dt$, but we shall avoid the explicit use of these processes
here. This technique can be used to work out a general relationship of type \eqref{bhbj} between
$B^H$ and $B^J$ for any $0<J,H<1$, see Theorem \ref{molgolo} (such a relation was obtained
by \cite{cjost06}).

Application of time inversion techniques also enables us to deduce in Theorem \ref{rmoins} a
canonical representation on $\reel_-$, and to obtain in Theorem \ref{noncan} some non canonical
representations of $B^H$ with respect to itself, extending the classical case $H=1/2$; these
representations are also considered by \cite{cjost07}.

Representations of type \eqref{bhw} or \eqref{bhbj} can be applied to descriptions of the
Cameron-Martin spaces $\hilbert_H$ of the fractional Brownian motions $B^H$; these spaces are
Hilbert spaces which characterise the laws of centred Gaussian processes (see Appendix
\ref{equiv}). The space $\hilbert_{1/2}$ is the classical space of absolutely continuous functions
$h$ such that $h(0)=0$ and the derivative $D^1h$ is square integrable, and \eqref{bhw} implies that
$\hilbert_H$ is the space of functions of the form
\begin{equation*}
    t\mapsto\frac1{\Gamma(H+1/2)}\int_\reel\Bigl((t-s)_+^{H-1/2}-(-s)_+^{H-1/2}\Bigr)f(s)ds
\end{equation*}
for square integrable functions $f$.

Sections \ref{riemliou} and \ref{series} are devoted to the comparison of $B^H$ with two
processes. One of them is self-similar but has only asymptotically stationary increments in large
time, and the other one has stationary increments, but is only asymptotically self-similar in small
time.

In Section \ref{riemliou}, we consider on $\reel_+$ the so-called Riemann-Liouville process
defined for $H>0$ by
\begin{equation*}
    X_t^H=\frac1{\Gamma(H+1/2)}\int_0^t(t-s)^{H-1/2}dW_s.
\end{equation*}
This process is $H$-self-similar but does not have stationary increments; contrary to $B^H$, the
parameter $H$ can be larger than 1. The Cameron-Martin space $\hilbert_H'$ of $X^H$ is the
space of functions
\begin{equation*}
    t\mapsto\frac1{\Gamma(H+1/2)}\int_0^t(t-s)^{H-1/2}f(s)ds
\end{equation*}
for square integrable functions $f$. We explain in Theorem \ref{rlcamar} a result of
\cite{samko:kil:mari93}, see \cite{decreuse:ustu99}, stating that $\hilbert_H$ and $\hilbert_H'$
are equivalent for $0<H<1$ (they are the same set with equivalent norms). We also compare the paths
of $B^H$ and $X^H$, and in particular study the equivalence or mutual singularity of the laws of
these processes (Theorem \ref{rlabs}); it appears that these two processes can be discriminated by
looking at their behaviour in small (or large) time. As an application, we also estimate the mutual
information of increments of $B^H$ on disjoint time intervals (more results of this type can be
found in \cite{norros:sak09}).

Another classical representation of the fractional Brownian motion on $\reel$ is its spectral
representation which can be written in the form
\begin{equation}\label{spectral}
   B_t^H=\frac1{\sqrt\pi}
   \int_0^{+\infty}s^{-1/2-H}\Bigl(\bigl(\cos(st)-1\bigr)dW_s^1+\sin(st)dW_s^2\Bigr),
\end{equation}
where $W_t^1$ and $W_t^2$, $t\ge0$, are two independent standard Brownian motions; it is indeed not
difficult to check that the right-hand side is Gaussian, centred, $H$-self-similar with stationary
increments, and $1/\sqrt\pi$ is the constant for which this process has the same variance as
\eqref{mandness} (see Appendix \ref{variance}). If now we are interested in $B^H$ on a bounded
interval, say $[0,1]$, we look for its Fourier coefficients. Thus the aim of Section \ref{series}
is to study the relationship between $B^H$ on $[0,1]$ and some series of trigonometric functions
with independent Gaussian coefficients. More precisely, the standard Brownian motion can be defined
on $[0,1]$ by series such as
\begin{equation}\label{wtrigo}
   W_t=\xi_0t+\sqrt2\,\sum_{n\ge1}\Bigl(\xi_n\frac{\cos(2n\pi t)-1}{2n\pi}
   +\xi_n'\frac{\sin(2n\pi t)}{2n\pi}\Bigr),
\end{equation}
\begin{equation}\label{wanti}
   W_t=\sqrt2\,\sum_{n\ge0}\Bigl(\xi_n\frac{\cos((2n+1)\pi t)-1}{(2n+1)\pi}
   +\xi_n'\frac{\sin((2n+1)\pi t)}{(2n+1)\pi}\Bigr),
\end{equation}
or
\begin{equation}\label{karh}
   W_t=\sqrt2\,\sum_{n\ge0}\xi_n\frac{\sin\bigl((n+1/2)\pi t\bigr)}{(n+1/2)\pi},
\end{equation}
where $\xi_n$, $\xi_n'$ are independent standard Gaussian variables. The form \eqref{karh} is the
Karhunen-Lo\`{e}ve expansion; it provides the orthonormal basis $\sqrt2\sin\bigl((n+1/2)\pi
t\bigr)$ of $L^2([0,1])$, such that the expansion of $W_t$ on this basis consists of independent
terms; it is a consequence of \eqref{wanti} which can be written on $[-1/2,1/2]$, and of the
property
\begin{equation*}
   W_t\simeq\sqrt2\,W_{t/2}\simeq W_{t/2}-W_{-t/2}.
\end{equation*}
It is not possible to write on $[0,1]$ the analogues of these formulae for $B^H$, $H\ne1/2$, but
it is possible (Theorem \ref{bhtrigo}) to write $B^H$ on $[0,1]$ as
\begin{equation}\label{bhanh}
     B_t^H=a_0^H\xi_0t
     +\sum_{n\ge1}a_n^H\Bigl(\bigl(\cos(\pi nt)-1\bigr)\xi_n+\sin(\pi nt)\xi_n'\Bigr)
\end{equation}
with $\sum(a_n^H)^2<\infty$. This result was proved in \cite{istas05} when $H\le1/2$, and the case
$H>1/2$ was studied in \cite{igloi05} with an approximation method. Formula \eqref{bhanh} is not
completely analogous to \eqref{wtrigo}, \eqref{wanti} or \eqref{karh}; contrary to these expansions
of $W$, the $\sigma$-algebra generated by $B^H$ in \eqref{bhanh} is strictly smaller than the
$\sigma$-algebra of the sequence $(\xi_n,\xi_n')$; in other words, the right hand side of
\eqref{bhanh} involves an extra information not contained in $B^H$, and this is a drawback for some
questions. This is why we define for $H>0$ a process
\begin{equation*}
   \Bhat_t^H=\xi_0t+\sqrt2\,\sum_{n\ge1}\Bigl(\xi_n
   \frac{\cos(2\pi nt)-1}{(2\pi n)^{H+1/2}}
   +\xi_n'\frac{\sin(2\pi nt)}{(2\pi n)^{H+1/2}}\Bigr)
\end{equation*}
which is a direct generalisation of \eqref{wtrigo}, and a similar process $\Bbar_t^H$ which
generalises \eqref{wanti}. It appears that for $0<H<1$, these processes have local properties
similar to $B^H$, and we can prove that their Cameron-Martin spaces are equivalent to $\hilbert_H$
(Theorem \ref{percamar}). As an application, we obtain Riesz bases of $\hilbert_H$, and show that
functions of $\hilbert_H$ can be characterised on $[0,1]$ by means of their Fourier coefficients.
We then study the equivalence or mutual singularity of the laws of $B^H$ and $\Bhat^H$, $\Bbar^H$
(Theorem \ref{perequiv}). We also discuss the extension of \eqref{karh} which has been proposed in
\cite{feyel:prad99}. In Theorem \ref{cherid}, we recover a result of \cite{cheridito01,vanzanten07}
which solves the following question: if we observe a path of a process, can we say whether it is a
pure fractional Brownian motion $B^J$, or whether this process $B^J$ has been corrupted by an
independent fractional Brownian motion of different index $H$?

Technical results which are required in our study are given in the three appendices:
\begin{itemize}
\item a lemma about some continuous endomorphisms of the standard Cameron-Martin space
    (Appendix \ref{analytical});
\item the computation of the variance of fractional Brownian motions (Appendix
    \ref{variance});
\item results about the equivalence and mutual singularity of laws of Gaussian processes, and
    about their relative entropies, with in particular a short review of Cameron-Martin spaces
    (Appendix \ref{equiv}).
\end{itemize}

Notice that many aspects concerning the fractional Brownian motion $B^H$ are not considered in this
work. Concerning the representations, it is possible to expand $B^H$ on a wavelet basis; we do not
consider this question to which several works have been devoted, see for instance
\cite{ymeyer:s:taqqu99}. We also do not study stochastic differential equations driven by $B^H$
(which can be solved by means of the theory of rough paths, see \cite{cou:qian02}), or the
simulation of fractional Brownian paths. On the other hand, fractional Brownian motions have
applications in many scientific fields, and we do not describe any of them.

\section{Fractional calculus}
\label{defnot}

Let us first give some notations. All random variables and processes are supposed to be defined on
a probability space $(\Omega,\tribuf,\proba)$ and the expectation is denoted by $\esp$; processes
are always supposed to be measurable functions $\Xi:(t,\omega)\mapsto\Xi_t(\omega)$, where $t$ is
in a subset of $\reel$ endowed with its Borel $\sigma$-algebra; the $\sigma$-algebra generated by
$\Xi$ is denoted by $\sigma(\Xi)$, and for the filtrations we use the notation \eqref{tribu}. The
derivative of order $n$ of $f$ is denoted by $D^nf$; the function is said to be smooth if it is
$C^\infty$. The function $f_1$ is said to be dominated by $f_2$ if $|f_1|\le Cf_2$. The notation
$u_n\asymp v_n$ means that $v_n/u_n$ is between two positive constants. We say that two Hilbert
spaces $\hilbert$ and $\hilbert'$ are equivalent (and write $\hilbert\sim\hilbert'$) if they are
the same set and
\begin{equation}\label{equivdef}
    C_1\|h\|_\hilbert\le\|h\|_{\hilbert'}\le C_2\|h\|_\hilbert
\end{equation}
for some positive $C_1$ and $C_2$; this means that the two spaces are continuously embedded into
each other. We often use the classical function $\Gamma$ defined on $\complex\setminus\relatif_-$,
and in particular the property $\Gamma(z+1)=z\,\Gamma(z)$.

We now describe the functional spaces, fractional integrals and derivatives which are used in this
work; see \cite{samko:kil:mari93} for a much more complete study of the fractional calculus. These
functional spaces are weighted H\"{o}lder spaces which are convenient for the study of the
fractional Brownian motion. The results are certainly not stated in their full generality, but are
adapted to our future needs.

\subsection{Functional spaces}
\label{functional}

The main property which is involved in our study is the H\"{o}lder continuity, but functions will
often exhibit a different behaviour near time 0 and for large times. More precisely, on the time
interval $\reel_+^\star$, let $\holder^{\beta,\gamma,\delta}$ for $0<\beta<1$ and $\gamma$,
$\delta$ real, be the Banach space of real functions $f$ such that
\begin{equation}\label{fbgd}
   \|f\|_{\beta,\gamma,\delta}=\sup_t\frac{|f(t)|}{t^\beta t^{\gamma,\delta}}
   +\sup_{s<t}\frac{\bigl|f(t)-f(s)\bigr|}{(t-s)^\beta\sup_{s\le u\le t}u^{\gamma,\delta}}
\end{equation}
is finite, with the notation
\begin{equation}\label{tgd}
   t^{\gamma,\delta}=t^\gamma 1_{\{t\le1\}}+t^\delta 1_{\{t>1\}}.
\end{equation}
Thus functions of this space are locally H\"{o}lder continuous with index $\beta$, and parameters
$\gamma$ and $\delta$ make more precise the behaviour at 0 and at infinity. If $\beta+\gamma>0$,
the function $f$ can be extended by continuity at 0 by $f(0)=\lim_0f=0$. If $\gamma\ge0$ and
$\delta\ge0$ and if we consider functions $f$ such that $\lim_0f=0$, then the second term of
\eqref{fbgd} dominates the first one (let $s$ decrease to 0).

\begin{remark}\label{dyadic}
Define
\begin{equation*}
   \|f\|_{\beta,\gamma,\delta}'=\sup\biggl\{\frac{\bigl|f(t)-f(s)\bigr|}{\bigl(2^n\bigr)^{\gamma,\delta}(t-s)^\beta},\;
   2^n\le s\le t\le2^{n+1},\;n\in\relatif\biggr\}.
\end{equation*}
Then this semi-norm is equivalent to the second term in \eqref{fbgd}; in particular, if
$\gamma\ge0$ and $\delta\ge0$, then $\|.\|_{\beta,\gamma,\delta}$ and
$\|.\|_{\beta,\gamma,\delta}'$ are equivalent on the space of functions $f$ such that $\lim_0f=0$.
It is indeed easy to see that $\|.\|_{\beta,\gamma,\delta}'$ is dominated by the second term of
\eqref{fbgd}. For the inverse estimation, notice that upper bounds for $|f(t)-f(s)|$ can be
obtained by adding the increments of $f$ on the dyadic intervals $[2^n,2^{n+1}]$ intersecting
$[s,t]$. More precisely, if $2^{k-1}\le s\le2^k\le2^n\le t\le2^{n+1}$, then
\begin{align*}
    \bigl|f(t)-f(s)\bigr|
    &\le\|f\|_{\beta,\gamma,\delta}'\sup_{k-1\le j\le n}\bigl(2^j\bigr)^{\gamma,\delta}
    \Bigl(\sum_{j=k}^{n-1}2^{j\beta}+(2^k-s)^\beta+(t-2^n)^\beta\Bigr)\\
    &\le C\,\|f\|_{\beta,\gamma,\delta}'\sup_{s\le u\le t}u^{\gamma,\delta}
    \Bigl(2^{n\beta}-2^{k\beta}+(2^k-s)^\beta+(t-2^n)^\beta\Bigr)\\
    &\le3C\,\|f\|_{\beta,\gamma,\delta}'\sup_{s\le u\le t}u^{\gamma,\delta}(t-s)^\beta
\end{align*}
because $2^{n\beta}-2^{k\beta}\le(2^n-2^k)^\beta\le(t-s)^\beta$.
\end{remark}

In particular, one can deduce from Remark \ref{dyadic} that $\holder^{\beta,\gamma,\delta}$ is
continuously embedded into $\holder^{\beta-\eps,\gamma+\eps,\delta+\eps}$ for $0<\eps<\beta$.

\begin{theorem}\label{prod}
The map $(f_1,f_2)\mapsto f_1f_2$ is continuous from
$\holder^{\beta,\gamma_1,\delta_1}\times\holder^{\beta,\gamma_2,\delta_2}$ into
$\holder^{\beta,\beta+\gamma_1+\gamma_2,\beta+\delta_1+\delta_2}$.
\end{theorem}

\begin{proof}
This is a bilinear map, so it is sufficient to prove that the image of a bounded subset is bounded.
If $f_1$ and $f_2$ are bounded in their respective H\"older spaces, it is easy to deduce that
$f_1(t)f_2(t)$ is dominated by $t^{2\beta}t^{\gamma_1+\gamma_2,\delta_1+\delta_2}$. On the other
hand, following Remark \ref{dyadic}, we verify that for $2^n\le s\le t\le2^{n+1}$,
\begin{align*}
   \bigl|f_1(t)f_2(t)-f_1(s)f_2(s)\bigr|
   &\le\bigl|f_1(s)\bigr|\>\bigl|f_2(t)-f_2(s)\bigr|+\bigl|f_2(t)\bigr|\>\bigl|f_1(t)-f_1(s)\bigr|\\
   &\le C\Bigl(s^\beta s^{\gamma_1,\delta_1}(2^n)^{\gamma_2,\delta_2}(t-s)^\beta
   +t^\beta t^{\gamma_2,\delta_2}(2^n)^{\gamma_1,\delta_1}(t-s)^\beta\Bigr)\\
   &\le C'(2^n)^\beta(2^n)^{\gamma_1,\delta_1}(2^n)^{\gamma_2,\delta_2}(t-s)^\beta.
\end{align*}
The theorem is therefore proved. \qed
\end{proof}

Let us define
\begin{equation*}
   \holder^{\beta,\gamma}=\holder^{\beta,\gamma,0},\qquad
   \holder^\beta=\holder^{\beta,0,0}.
\end{equation*}
These spaces can be used for functions defined on a finite time interval $[0,T]$, since in this
case the parameter $\delta$ is unimportant. For functions defined on $\reel_-^\star$, we say that
$f$ is in $\holder^{\beta,\gamma,\delta}$ if $t\mapsto f(-t)$ is in it, and for functions defined
on a general interval of $\reel$, we assume that the restrictions to $\reel_+^\star$ and
$\reel_-^\star$ are in $\holder^{\beta,\gamma,\delta}$. For $\gamma=0$, the regularity at time 0 is
similar to other times, so spaces $\holder^{\beta,0,\delta}$ are invariant by the time shifts
$f\mapsto f(.+t_0)-f(t_0)$. If we consider a time interval of type $[1,+\infty)$, then the
parameter $\gamma$ can be omitted and we denote the space by $\holder^{\beta,.,\delta}$.

We use the notations
\begin{equation}\label{hbmoins}
   \holder^{\beta-,\gamma,\delta+}=\bigcap_{\eps>0}\holder^{\beta-\eps,\gamma,\delta+2\eps},
   \quad\holder^{\beta-,\gamma}=\bigcap_{\eps>0}\holder^{\beta-\eps,\gamma},\quad
   \holder^{\beta-}=\bigcap_{\eps>0}\holder^{\beta-\eps}.
\end{equation}
They are Fr\'{e}chet spaces.

\begin{example}
If $B^H$ is a $H$-fractional Brownian motion on the time interval $[0,1]$, the probability of the
event $\{B^H\in\holder^\beta\}$ is 1 if $\beta<H$ (this follows from the Kolmogorov continuity
theorem). In particular, $B^H$ lives almost surely in $\holder^{H-}$. We shall see in Remark
\ref{hhmz} that this implies that on the time interval $\reel_+$, the process $B^H$ lives in
$\holder^{H-,0,0+}$.
\end{example}

The parameters $\gamma$ and $\delta$ can be modified by means of some multiplication operators.
More precisely, on $\reel_+^\star$, define
\begin{equation}\label{pialphadef}
   \Pi^\alpha f(t)=t^\alpha f(t),\qquad\Pi^{\alpha_1,\alpha_2}f(t)=t^{\alpha_1}(1+t)^{\alpha_2-\alpha_1}f(t).
\end{equation}

\begin{theorem}\label{pialphath}
The operator $\Pi^{\alpha_1,\alpha_2}$ maps continuously $\holder^{\beta,\gamma,\delta}$ into
$\holder^{\beta,\gamma+\alpha_1,\delta+\alpha_2}$. In particular, on the time interval $(0,1]$, the
operator $\Pi^\alpha$ maps continuously $\holder^{\beta,\gamma}$ into
$\holder^{\beta,\gamma+\alpha}$.
\end{theorem}

\begin{proof}
The quantity $|t^\alpha-s^\alpha|(t-s)^{-\beta}t^{\beta-\alpha}$ is bounded for $2^n\le s\le
t\le2^{n+1}$, and the bound does not depend on $n$ (use the scaling). Thus it follows from Remark
\ref{dyadic} that the function $t\mapsto t^\alpha$ is in
$\holder^{\beta,\alpha-\beta,\alpha-\beta}$. The same property implies that
$(1+t)^\alpha-(1+s)^\alpha$ is dominated by $(1+t)^{\alpha-\beta}(t-s)^\beta$ (with the same
assumptions on $s$ and $t$), and we can deduce that $t\mapsto(1+t)^\alpha$ is in
$\holder^{\beta,-\beta,\alpha-\beta}$ (the coefficient $-\beta$ is due to the fact that the
function tends to 1 at 0). We deduce from Theorem \ref{prod} that the function
$t^{\alpha_1}(1+t)^{\alpha_2-\alpha_1}$ is in $\holder^{\beta,\alpha_1-\beta,\alpha_2-\beta}$. The
operator $\Pi^{\alpha_1,\alpha_2}$ is the multiplication by this function, and the result follows
by again applying Theorem \ref{prod}. \qed
\end{proof}

It is then possible to deduce a density result for the spaces of \eqref{hbmoins} (the result is
false with $\beta$ instead of $\beta-$). Fractional polynomials are linear combinations of
monomials $t^\alpha$, $\alpha\in\reel$, and these monomials are in $\holder^{\beta,\gamma}$ on
$(0,1]$ if $\alpha\ge\beta+\gamma$.

\begin{theorem}\label{dense}
Let $0<\beta<1$.
\begin{itemize}
\item On $(0,1]$, fractional polynomials (belonging to $\holder^{\beta-,\gamma}$) are dense
    in $\holder^{\beta-,\gamma}$.
\item On $\reel_+^\star$, smooth functions with compact support are dense in
    $\holder^{\beta-,\gamma,\delta+}$.
\end{itemize}
\end{theorem}

\begin{proof}
Let us consider separately the two statements.

\medskip\noindent\emph{Study on $(0,1]$.} The problem can be reduced to the case $\gamma=0$ with Theorem
\ref{pialphath}, and functions $f$ of $\holder^{\beta-}$ are continuous on the closed interval
$[0,1]$ with $f(0)=0$. If $f$ is in $\holder^{\beta-\eps}$ (for $\eps$ small), it can be
approximated by classical polynomials $f_n$ by means of the Stone-Weierstrass theorem; more
precisely, if we choose the Bernstein approximations $\esp f\bigl(\frac1n\sum_{j=1}^n1_{\{U_j\le
x\}}\bigr)$ for independent uniformly distributed variables $U_j$ in $[0,1]$, then $f_n$ is bounded
in $\holder^{\beta-\eps}$ and converges uniformly to $f$. Thus
\begin{align}
   \bigl|f_n(t)&-f_n(s)-f(t)+f(s)\bigr|\notag\\
   &\le C\bigl(|f_n(t)-f_n(s)|^{(\beta-2\eps)/(\beta-\eps)}+|f(t)-f(s)|^{(\beta-2\eps)/(\beta-\eps)}\bigr)
   \notag\\&\qquad\sup_u|f_n(u)-f(u)|^{\eps/(\beta-\eps)}\notag\\
   &\le C'(t-s)^{\beta-2\eps}\sup_u|f_n(u)-f(u)|^{\eps/(\beta-\eps)}.\label{fntfns}
\end{align}
These inequalities can also be written for $s=0$ to estimate $|f_n(t)-f(t)|$, so $f_n$ converges to
$f$ in $\holder^{\beta-2\eps}$.

\medskip\noindent\emph{Study on $\reel_+^\star$.} The technique is similar. By means of
$\Pi^{\alpha_1,\alpha_2}$, we can reduce the study to the case $\gamma=0$ and
$-2\beta<\delta<-\beta$. Let $f$ be in  $\holder^{\beta-,0,\delta+}$ and let us fix a small
$\eps>0$; then $f$ is in $\holder^{\beta-\eps,0,\delta+2\eps}$; in particular, it tends to 0 at 0
and at infinity. A standard procedure enables to approximate it uniformly by smooth functions $f_n$
with compact support, such that $f_n$ is bounded in $\holder^{\beta-\eps,0,\delta+2\eps}$; to this
end, we first multiply $f$ by the function $\phi_n$ supported by $[2^{-n-1},2^{n+1}]$, taking the
value 1 on $[2^{-n},2^n]$, and which is affine on $[2^{-n-1},2^{-n}]$ and on $[2^n,2^{n+1}]$; then
we take the convolution of $f\,\phi_n$ with $2^{n+2}\,\psi(2^{n+2}t)$ for a smooth function $\psi$
supported by $[-1,1]$ and with integral 1. By proceeding as in \eqref{fntfns}, we can see that
\begin{align*}
   \bigl|f_n(t)&-f_n(s)-f(t)+f(s)\bigr|\\
   &\le C(t-s)^{\beta-2\eps}\sup_{s\le u\le t}\bigl(u^{0,\delta+2\eps}\bigr)^{(\beta-2\eps)/(\beta-\eps)}
   \sup_u|f_n(u)-f(u)|^{\eps/(\beta-\eps)}
\end{align*}
so $f_n$ converges to $f$ in $\holder^{\beta-2\eps,0,\delta+4\eps}$ because
$(\delta+2\eps)(\beta-2\eps)/(\beta-\eps)\le\delta+4\eps$ for $\eps$ small enough. \qed
\end{proof}

\subsection{Riemann-Liouville operators}

An important tool for the stochastic calculus of fractional Brownian motions is the fractional
calculus obtained from the study of Riemann-Liouville operators $I_\pm^\alpha$. These operators can
be defined for any real index $\alpha$ (and even for complex indices), but we will mainly focus on
the case $|\alpha|<1$.

\subsubsection{Operators with finite horizon}

The fractional integral operators $I_{\tau\pm}^\alpha$ (Riemann-Liouville operators) are defined
for $\tau\in\reel$ and $\alpha>0$ by
\begin{equation}\label{riemdef}
   I_{\tau+}^\alpha f(t)=\frac1{\Gamma(\alpha)}\int_\tau^t(t-s)^{\alpha-1}f(s)ds,\quad
   I_{\tau-}^\alpha f(t)=\frac1{\Gamma(\alpha)}\int_t^\tau(s-t)^{\alpha-1}f(s)ds,
\end{equation}
respectively for $t>\tau$ and $t<\tau$. These integrals are well defined for instance if $f$ is
locally bounded on $(\tau,+\infty)$ or $(-\infty,\tau)$, and is integrable near $\tau$. If $f$ is
integrable, they are defined almost everywhere, and $I_{\tau\pm}^\alpha$ is a continuous
endomorphism of $L^1([\tau,T])$ or $L^1([T,\tau])$. These operators satisfy the semigroup property
\begin{equation}\label{semig}
   I_{\tau\pm}^{\alpha_2}I_{\tau\pm}^{\alpha_1}=I_{\tau\pm}^{\alpha_1+\alpha_2}
\end{equation}
which can be proved from the relation between Beta and Gamma functions recalled in \eqref{beta}. If
$\alpha$ is an integer, we get iterated integrals; in particular, $I_{\tau\pm}^1f$ is $\pm$ the
primitive of $f$ taking value 0 at $\tau$. Notice that relations \eqref{riemdef} can also be
written as
\begin{equation}\label{riemvar}
  \begin{split}
   I_{\tau+}^\alpha f(t)&=\frac1{\Gamma(\alpha)}\int_\tau^t(t-s)^{\alpha-1}\bigl(f(s)-f(t)\bigr)ds
   +\frac{(t-\tau)^\alpha}{\Gamma(\alpha+1)}f(t),\\
   I_{\tau-}^\alpha f(t)&=\frac1{\Gamma(\alpha)}\int_t^\tau(s-t)^{\alpha-1}\bigl(f(s)-f(t)\bigr)ds
   +\frac{(\tau-t)^\alpha}{\Gamma(\alpha+1)}f(t).
  \end{split}
\end{equation}
If $f$ is Lipschitz with $f(\tau)=0$, an integration by parts shows that
\begin{equation}\label{riemlip}
   I_{\tau+}^\alpha f(t)=\frac1{\Gamma(\alpha+1)}\int_\tau^t(t-s)^\alpha df(s),\quad
   I_{\tau-}^\alpha f(t)=\frac{-1}{\Gamma(\alpha+1)}\int_t^\tau(s-t)^\alpha df(s).
\end{equation}
For $\alpha=0$, the operators $I_{\tau\pm}^0$ are by definition the identity (this is coherent with
\eqref{riemlip}). The study of the operators $I_{\tau\pm}^\alpha$ can be reduced to the study of
$I_{0+}^\alpha$, since the other cases can be deduced by means of an affine change of time.

\begin{example}
The value of $I_{0+}^\alpha$ on fractional polynomials can be obtained from
\begin{equation}\label{itbeta}
   I_{0+}^\alpha\Bigl(\frac{t^\beta}{\Gamma(\beta+1)}\Bigr)=\frac{t^{\alpha+\beta}}{\Gamma(\alpha+\beta+1)}
\end{equation}
which is valid for $\beta>-1$.
\end{example}

Riemann-Liouville operators can also be defined for negative exponents, and are called fractional
derivatives. Here we restrict ourselves to $-1<\alpha<0$, and in this case the derivative of order
$-\alpha$ is defined by
\begin{equation}\label{riemneg}
   I_{\tau+}^\alpha f=D^1I_{\tau+}^{1+\alpha}f,\qquad I_{\tau-}^\alpha f=-D^1I_{\tau-}^{1+\alpha}f
\end{equation}
if $I_{\tau\pm}^{1+\alpha}f$ is absolutely continuous, for the differentiation operator $D^1$. The
relation \eqref{itbeta} is easily extended to negative $\alpha$ (with result 0 if
$\alpha+\beta+1=0$). Fractional derivatives operate on smooth functions, and we have the following
result.

\begin{theorem}\label{smooth}
Suppose that $f$ is smooth and integrable on $(0,1]$. Then, for any $\alpha>-1$, $I_{0+}^\alpha f$
is well defined, is smooth on $(0,1]$, and
\begin{equation}\label{d1i}
\begin{split}
   &\bigl|D^1I_{0+}^\alpha f(t)\bigr|\\&\le C_\alpha\Bigl(t^{\alpha-2}\int_0^{t/2}|f(s)|ds
   +t^{\alpha-1}\sup_{[t/2,t]}|f|+t^\alpha\sup_{[t/2,t]}|D^1f|+t^{\alpha+1}\sup_{[t/2,t]}|D^2f|\Bigr).
\end{split}
\end{equation}
If $D^1f$ is integrable and $\lim_0f=0$, then $D^1I_{0+}^\alpha f=I_{0+}^\alpha D^1f$.
\end{theorem}

\begin{proof}
First suppose $\alpha>0$. Then, for $t>u>0$, we can write \eqref{riemdef} in the form
\begin{equation}\label{dia}
   I_{0+}^\alpha f(t)=\Gamma(\alpha)^{-1}\Bigl(\int_0^u(t-s)^{\alpha-1}f(s)ds
   +\int_0^{t-u}s^{\alpha-1}f(t-s)ds\Bigr).
\end{equation}
This expression is smooth, and
\begin{equation}\label{dialpha}
\begin{split}
   D^1I_{0+}^\alpha f(t)=\Gamma(\alpha)^{-1}\Bigl(&(\alpha-1)\int_0^u(t-s)^{\alpha-2}f(s)ds\\
   &+\int_0^{t-u}s^{\alpha-1}D^1f(t-s)ds+(t-u)^{\alpha-1}f(u)\Bigr).
\end{split}
\end{equation}
In particular, by letting $u=t/2$, we obtain \eqref{d1i} without the $D^2f$ term. Moreover, if
$D^1f$ is integrable and $\lim_0f=0$, we see by writing
\begin{equation*}
   (t-u)^{\alpha-1}f(u)=-(\alpha-1)\int_0^u(t-s)^{\alpha-2}f(s)ds+\int_0^u(t-s)^{\alpha-1}D^1f(s)ds
\end{equation*}
that
\begin{align*}
   D^1I_{0+}^\alpha f(t)&=\Gamma(\alpha)^{-1}\Bigl(\int_0^{t-u}s^{\alpha-1}D^1f(t-s)ds
   +\int_0^u(t-s)^{\alpha-1}D^1f(s)ds\Bigr)\\
   &=I_{0+}^\alpha D^1f(t)
\end{align*}
(apply \eqref{dia} with $f$ replaced by $D^1f$). Let us now consider the case $-1<\alpha<0$; we use
the definition \eqref{riemneg} of the fractional derivative, and in particular deduce that
$I_{0+}^\alpha f$ is again smooth. Moreover, from \eqref{dialpha},
\begin{align*}
   D^1I_{0+}^\alpha f(t)&=D^2I_{0+}^{\alpha+1}f(t)\\
   &=\Gamma(\alpha+1)^{-1}\Bigl(\alpha(\alpha-1)\int_0^u(t-s)^{\alpha-2}f(s)ds
   +(t-u)^\alpha D^1f(u)\\&\quad+\int_0^{t-u}s^\alpha D^2f(t-s)ds+\alpha(t-u)^{\alpha-1}f(u)\Bigr).
\end{align*}
We deduce \eqref{d1i} by letting again $u=t/2$. If $\lim_0f=0$ and $D^1f$ is integrable, then
\begin{equation*}
 D^1I_{0+}^\alpha f=D^2I_{0+}^{\alpha+1}f=D^1I_{0+}^{\alpha+1}D^1f=I_{0+}^\alpha D^1f
\end{equation*}
from the definition \eqref{riemneg} and the property for $\alpha+1$ which has already been proved.
\qed
\end{proof}

For $-1<\alpha<0$, a study of \eqref{riemvar} shows that $I_{\tau\pm}^\alpha f$ is defined as soon
as $f$ is H\"{o}lder continuous with index greater than $-\alpha$, and that \eqref{riemvar} again holds
true. If $f$ is Lipschitz and $f(\tau)=0$, then we can write
\begin{equation*}
   I_{\tau\pm}^\alpha f=\pm D^1I_{\tau\pm}^{1+\alpha}f
   =D^1I_{\tau\pm}^{1+\alpha}I_{\tau\pm}^1D^1f=D^1I_{\tau\pm}^1I_{\tau\pm}^{1+\alpha}D^1f
   =\pm I_{\tau\pm}^{1+\alpha}D^1f
\end{equation*}
where we have used \eqref{semig} in the third equality, so \eqref{riemlip} again holds true. Thus
relations \eqref{riemvar} and \eqref{riemlip} can be used for any $\alpha>-1$ ($\alpha\ne0$ for
\eqref{riemvar}). By using the multiplication operators $\Pi^\alpha$ defined in \eqref{pialphadef},
we can deduce from \eqref{riemvar} a formula for weighted fractional operators; if $f$ is smooth
with compact support in $\reel_+^\star$, then
\begin{equation}\label{pimgam}
    \Pi^{-\gamma}I_{0+}^\alpha\Pi^\gamma f(t)=I_{0+}^\alpha f(t)
    +\frac1{\Gamma(\alpha)}\int_0^t(t-s)^{\alpha-1}\Bigl(\bigl(\frac st\bigr)^\gamma-1\Bigr)f(s)ds
\end{equation}
for $\alpha>-1$, $\alpha\ne0$.

Here are some results about $I_{0+}^\alpha$ related to the functional spaces of Subsection
\ref{functional}. They can easily be translated into properties of $I_{\tau\pm}^\alpha$, see also
\cite{samko:kil:mari93,norros:val:vir99}.

\begin{theorem}\label{riemth}
Consider the time interval $(0,1]$ and let $\gamma>-1$.
\begin{itemize}
\item If $\beta$ and $\beta+\alpha$ are in $(0,1)$, then the operator $I_{0+}^\alpha$ maps
    continuously $\holder^{\beta,\gamma}$ into $\holder^{\beta+\alpha,\gamma}$.
\item The composition rule $I_{0+}^{\alpha_2}I_{0+}^{\alpha_1}=I_{0+}^{\alpha_1+\alpha_2}$
    holds on $\holder^{\beta,\gamma}$ provided $\beta$, $\beta+\alpha_1$ and
    $\beta+\alpha_1+\alpha_2$ are in $(0,1)$.
\end{itemize}
\end{theorem}

\begin{proof}
Let us prove the first statement. Let $f$ be in $\holder^{\beta,\gamma}$. The property
$I_{0+}^\alpha f(t)=O(t^{\alpha+\beta+\gamma})$ can be deduced from \eqref{riemvar} and
\eqref{itbeta}. By applying Remark \ref{dyadic}, it is then sufficient to compare $I_{0+}^\alpha f$
at times $s$ and $t$ for $2^n\le s\le t\le2^{n+1}$, $n<0$. Consider the time $v=(3s-t)/2$, so that
$2^{n-1}\le s/2\le v\le s\le2^{n+1}$. By again applying \eqref{riemvar}, we have
\begin{align*}
   I_{0+}^\alpha f(t)-I_{0+}^\alpha f(s)
   =&\frac{t^\alpha f(t)-s^\alpha f(s)}{\Gamma(\alpha+1)}
   +\frac{A_{v,t}-A_{v,s}}{\Gamma(\alpha)}
   +\frac{f(s)-f(t)}{\Gamma(\alpha)}\int_0^v(t-u)^{\alpha-1}du\\
   &+\frac1{\Gamma(\alpha)}\int_0^v\bigl((t-u)^{\alpha-1}-(s-u)^{\alpha-1}\bigr)\bigl(f(u)-f(s)\bigr)du
\end{align*}
with
\begin{equation*}
   A_{v,w}=\int_v^w(w-u)^{\alpha-1}\bigl(f(u)-f(w)\bigr)du=O\bigl((v^\gamma+w^\gamma)(w-v)^{\alpha+\beta}\bigr).
\end{equation*}
We deduce that
\begin{equation}\label{iaf}
\begin{split}
   I_{0+}^\alpha f(t)-I_{0+}^\alpha f(s)
   =&\frac{\bigl(t^\alpha-s^\alpha\bigr)f(s)}{\Gamma(\alpha+1)}
   +\frac{A_{v,t}-A_{v,s}}{\Gamma(\alpha)}
   -\frac{f(s)-f(t)}{\Gamma(\alpha+1)}(t-v)^\alpha\\
   &+\frac1{\Gamma(\alpha)}\int_0^v\bigl((t-u)^{\alpha-1}-(s-u)^{\alpha-1}\bigr)\bigl(f(u)-f(s)\bigr)du.
\end{split}
\end{equation}
The second and third terms are easily shown to be dominated by $2^{n\gamma}(t-s)^{\alpha+\beta}$.
The first term is dominated by
\begin{equation*}
   \sup_{s\le u\le t}u^{\alpha-1}(t-s)s^{\beta+\gamma}\le C\,2^{n\gamma}(t-s)^{\alpha+\beta}.
\end{equation*}
The last term is dominated by
\begin{align*}
   \int_0^v&\bigl((s-u)^{\alpha-1}-(t-u)^{\alpha-1}\bigr)(s-u)^\beta\bigl(u^\gamma+s^\gamma\bigr)du\\
   &\le(1-\alpha)(t-s)\int_0^v(s-u)^{\alpha+\beta-2}(u^\gamma+s^\gamma\bigr)du\displaybreak[0]\\
   &\le C(t-s)\Bigl(2^{n\gamma}(s-v)^{\alpha+\beta-1}+\int_0^{s/2}(s-u)^{\alpha+\beta-2}(u^\gamma+s^\gamma\bigr)du\Bigr)\\
   &\le C'2^{n\gamma}(t-s)^{\alpha+\beta}
\end{align*}
because $s-v=(t-s)/2$ and the integral on $[0,s/2]$ is proportional to
$s^{\alpha+\beta+\gamma-1}\le c\,2^{n(\alpha+\beta+\gamma-1)}\le
c\,2^{n\gamma}(t-s)^{\alpha+\beta-1}$. Thus the continuity of $I_{0+}^\alpha$ is proved. For the
composition rule, it is easily verified for monomials $f(t)=t^\beta$ (apply \eqref{itbeta}), and is
then extended by density to the space $\holder^{\beta-,\gamma}$ from Theorem \ref{dense}. By
applying this property to a slightly larger value of $\beta$, it appears that the composition rule
actually holds on $\holder^{\beta,\gamma}$. \qed
\end{proof}

Notice that fractional monomials $t^\kappa$ are eigenfunctions of $\Pi^{-\alpha}I_{0+}^\alpha$
and $I_{0+}^\alpha\Pi^{-\alpha}$ when they are in the domains of definitions of these
operators, so when $\kappa$ is large enough. This implies that these operators commute on
fractional polynomials. This property is then extended to other functions by density. In particular,
\begin{equation}\label{ipii}
   I_{0+}^{\alpha_2}\Pi^{-\alpha_1-\alpha_2}I_{0+}^{\alpha_1}
   =\Pi^{-\alpha_1}I_{0+}^{\alpha_1+\alpha_2}\Pi^{-\alpha_2},
\end{equation}
see (10.6) in \cite{samko:kil:mari93}.

\subsubsection{Operators with infinite horizon}
\label{infhor}

The operators $I_{\pm}^\alpha$ are defined by letting $\tau\to\mp\infty$ in
$I_{\tau\pm}^\alpha$. However, we will be more interested in the modified operators
\begin{equation*}
    \Itilde_\pm^\alpha f(t)=I_\pm^\alpha f(t)-I_\pm^\alpha f(0)
    =\lim_{\tau\to\mp\infty}\bigl(I_{\tau\pm}^\alpha f(t)-I_{\tau\pm}^\alpha f(0)\bigr)
\end{equation*}
when the limit exists. For $\alpha>0$, we can write
\begin{equation}\label{itildedef}
 \begin{split}
   \Itilde_+^\alpha f(t)&=\frac1{\Gamma(\alpha)}\int\bigl((t-s)_+^{\alpha-1}
   -(-s)_+^{\alpha-1}\bigr)f(s)ds,\\
   \Itilde_-^\alpha f(t)&=\frac1{\Gamma(\alpha)}\int\bigl((s-t)_+^{\alpha-1}
   -s_+^{\alpha-1}\bigr)f(s)ds
 \end{split}
\end{equation}
where we use the notation $u_+^\lambda=u^\lambda 1_{\{u>0\}}$. These integrals are well defined if
$f(t)$ is dominated by $(1+|t|)^\delta$ for $\delta<1-\alpha$ (there are also cases where the
integrals are only semi-convergent). In particular, the fractional integrals are generally not
defined for large values of $\alpha$, as it was the case for $I_{0+}^\alpha$. We are going to study
$\Itilde_\pm^\alpha$ on the functional spaces $\holder^{\beta,0,\delta}$.

\begin{remark}
The operator $\Itilde_\pm^\alpha$ is a normalisation of $I_\pm^\alpha$ in the sense that it can be
defined in more cases than $I_\pm^\alpha f$. For instance, for $\alpha>0$, if we compare
$I_-^\alpha f$ and $\Itilde_-^\alpha f$ on $\reel_+^\star$ for $f(s)=s^\delta$, we see that the
former one is defined for $\delta<-\alpha$, whereas the latter one is defined for
$\delta<1-\alpha$.
\end{remark}

Let us now consider the case $-1<\alpha<0$; we can let $\tau$ tend to infinity in \eqref{riemvar}
and obtain
\begin{equation}\label{martilde}
  \begin{split}
   \Itilde_+^\alpha f(t)&=\frac1{\Gamma(\alpha)}\int\Bigl((t-s)_+^{\alpha-1}\bigl(f(s)-f(t)\bigr)
   -(-s)_+^{\alpha-1}\bigl(f(s)-f(0)\bigr)\Bigr)ds,\\
   \Itilde_-^\alpha f(t)&=\frac1{\Gamma(\alpha)}\int\Bigl((s-t)_+^{\alpha-1}\bigl(f(s)-f(t)\bigr)
   -s_+^{\alpha-1}\bigl(f(s)-f(0)\bigr)\Bigr)ds.
  \end{split}
\end{equation}
This expression is defined on $\holder^{\beta,0,\delta}$ provided $\beta+\alpha>0$ and
$\beta+\alpha+\delta<1$.

Let $\alpha>-1$. Suppose that $f$ is Lipschitz and has compact support, so that $f$ is 0 on
$(-\infty,\tau]$, respectively $[\tau,+\infty)$. Then $I_\pm^\alpha f=I_{\tau\pm}^\alpha f$ on
$[\tau,+\infty)$, respectively $(-\infty,\tau]$, so $\Itilde_{\pm}^\alpha f(t)$ is equal to
$I_{\tau\pm}^\alpha f(t)-I_{\tau\pm}^\alpha f(0)$, which can be expressed by means of
\eqref{riemlip}. Thus
\begin{equation}\label{rieminf}
\begin{split}
   \Itilde_+^\alpha f(t)&=\frac1{\Gamma(\alpha+1)}\int\bigl((t-s)_+^\alpha-(-s)_+^\alpha\bigr)df(s),\\
   \Itilde_-^\alpha f(t)&=\frac{1}{\Gamma(\alpha+1)}\int\bigl(s_+^\alpha-(s-t)_+^\alpha\bigr)df(s).
\end{split}
\end{equation}
By applying Theorem \ref{smooth}, we see that if $f$ is smooth with compact support, then
$\Itilde_\pm^\alpha f$ is smooth and
\begin{equation}\label{ditil}
 D^1\Itilde_\pm^\alpha f=D^1I_\pm^\alpha f=I_\pm^\alpha D^1f.
\end{equation}

\begin{remark}\label{fzeroplus}
If $f=0$ on $\reel_+$ and if we look for $\Itilde_+^\alpha f$ on $\reel_+^\star$, we see when
$\alpha<0$ that $f(0)$ and $f(t)$ disappear in \eqref{martilde}, so \eqref{itildedef} can be used
on $\reel_+^\star$ for both positive and negative $\alpha$, and $\Itilde_+^\alpha f$ is $C^\infty$
on $\reel_+^\star$.
\end{remark}

\begin{theorem}\label{itildeth}
Consider the operators $\Itilde_+^\alpha$ and $\Itilde_-^\alpha$ on the respective time intervals
$(-\infty,T]$ for $T\ge0$, and $[T,+\infty)$ for $T\le0$. Let $\delta>0$.
\begin{itemize}
\item The operator $\Itilde_\pm^\alpha$ maps continuously $\holder^{\beta,0,\delta}$ into
    $\holder^{\beta+\alpha,0,\delta}$ provided $\beta$, $\beta+\alpha$ and
    $\beta+\alpha+\delta$ are in $(0,1)$.
\item The composition rule
    $\Itilde_\pm^{\alpha_2}\Itilde_\pm^{\alpha_1}=\Itilde_\pm^{\alpha_1+\alpha_2}$ holds on
    $\holder^{\beta,0,\delta}$ provided $\beta$, $\beta+\alpha_1$, $\beta+\alpha_1+\alpha_2$,
    $\beta+\alpha_1+\delta$ and $\beta+\alpha_1+\alpha_2+\delta$ are in $(0,1)$.
\end{itemize}
\end{theorem}

\begin{proof}
It is of course sufficient to study $\Itilde_+^\alpha$. We prove separately the two statements.

\medskip\noindent\emph{Continuity of $\Itilde_+^\alpha$.}
We want to study the continuity on the time interval $(-\infty,T]$; by means of a time shift, let
us consider the time interval $(-\infty,-1]$, and let us prove that if $f$ is in
$\holder^{\beta,.,\delta}$, then the function $\lim_{\tau\to-\infty}(I_{\tau+}^\alpha
f(t)-I_{\tau+}^\alpha f(-1))$ is in $\holder^{\beta+\alpha,.,\delta}$. From Remark \ref{dyadic}, it
is sufficient to estimate the increments of this function on intervals
$[s,t]\subset[-2^{n+1},-2^n]$ for $n\ge0$. Consider the proof of Theorem \ref{riemth} where
$I_{0+}^\alpha$ is replaced by $I_{\tau+}^\alpha$, and let us estimate $I_{\tau+}^\alpha
f(t)-I_{\tau+}^\alpha f(s)$ for $\tau\to-\infty$. We can write a formula similar to \eqref{iaf}.
The first term involves $(t-\tau)^\alpha-(s-\tau)^\alpha$ which tends to 0 as $\tau\to-\infty$, so
this first term vanishes. The second and third terms are dealt with similarly to Theorem
\ref{riemth}; the only difference is that the weight $2^{n\gamma}$ now becomes $2^{n\delta}$. The
last term is an integral on $(-\infty,v)$ and is dominated by
\begin{align*}
   (t-s)\int_{-\infty}^v(s-u)^{\alpha+\beta-2}|u|^\delta du
   &=(t-s)\int_{(t-s)/2}^{+\infty}u^{\alpha+\beta-2}(u-s)^\delta du\\
   &\le(t-s)\int_{(t-s)/2}^{+\infty}\bigl(u^{\alpha+\beta+\delta-2}+u^{\alpha+\beta-2}|s|^\delta\bigr)du\\
   &\le C(t-s)\Bigl((t-s)^{\alpha+\beta+\delta-1}+(t-s)^{\alpha+\beta-1}|s|^\delta\Bigr)\\
   &\le2C(t-s)^{\alpha+\beta}|s|^\delta.
\end{align*}

\medskip\noindent\emph{Composition rule.} If $f$ is 0 before some time $\tau_0$, then
$\Itilde_+^{\alpha_1}f(t)=I_{\tau+}^{\alpha_1}f(t)-I_{\tau+}^{\alpha_1}f(0)$ for
$\tau\le\tau_0\wedge t$. Thus
\begin{equation*}
   \Itilde_+^{\alpha_2}\Itilde_+^{\alpha_1}f(t)=\lim_{\tau\to-\infty}\bigl(I_{\tau+}^{\alpha_2}\Itilde_+^{\alpha_1}f(t)
   -I_{\tau+}^{\alpha_2}\Itilde_+^{\alpha_1}f(0)\bigr)
\end{equation*}
with
\begin{equation*}
   I_{\tau+}^{\alpha_2}\Itilde_+^{\alpha_1}f(t)=I_{\tau+}^{\alpha_2}I_{\tau+}^{\alpha_1}f(t)
   -\frac{(t-\tau)^{\alpha_2}}{\Gamma(\alpha_2+1)}I_{\tau+}^{\alpha_1}f(0)
   =I_{\tau+}^{\alpha_1+\alpha_2}f(t)
   -\frac{(t-\tau)^{\alpha_2}}{\Gamma(\alpha_2+1)}I_+^{\alpha_1}f(0)
\end{equation*}
from Theorem \ref{riemth}. Thus
\begin{equation*}
   \Itilde_+^{\alpha_2}\Itilde_+^{\alpha_1}f(t)=\Itilde_+^{\alpha_1+\alpha_2}f(t)
   -\lim_{\tau\to-\infty}\frac{(t-\tau)^{\alpha_2}-(-\tau)^{\alpha_2}}{\Gamma(\alpha_2+1)}I_+^{\alpha_1}f(0)
   =\Itilde_+^{\alpha_1+\alpha_2}f(t).
\end{equation*}
The case of general functions is then deduced from the density of functions with compact support in
$\holder^{\beta-,0,\delta+}$ (Theorem \ref{dense}); the proof on $\holder^{\beta,0,\delta}$ is
obtained as in Theorem~\ref{riemth} by increasing $\beta$ and decreasing $\delta$ slightly. \qed
\end{proof}

In particular, we deduce from Theorem \ref{itildeth} that $\Itilde_\pm^\alpha$ is a homeomorphism
from $\holder^{\beta-,0,0+}$ onto $\holder^{(\alpha+\beta)-,0,0+}$ if $\beta$ and $\alpha+\beta$
are in $(0,1)$, and $\Itilde_\pm^{-\alpha}$ is its inverse map.

\subsubsection{Operators for periodic functions}
\label{periodop}

Consider a bounded 1-periodic function $f$. Let $|\alpha|<1$; if $\alpha<0$, suppose moreover that
$f$ is in $\holder^\beta$ for some $\beta>-\alpha$. Then $\Itilde_+^\alpha f$ is well defined and
is given by \eqref{itildedef} or \eqref{martilde}; moreover, this function is also 1-periodic, and
is 0 at time 0; this follows from
\begin{equation*}
    I_{\tau+}^\alpha f(t+1)=I_{(\tau-1)+}^\alpha f(t)
\end{equation*}
so that
\begin{equation*}
    I_{\tau+}^\alpha f(t+1)- I_{\tau+}^\alpha f(0)=
    \bigl(I_{(\tau-1)+}^\alpha f(t)-I_{(\tau-1)+}^\alpha f(0)\bigr)
    +\bigl(I_{(\tau-1)+}^\alpha f(0)- I_{\tau+}^\alpha f(0)\bigr).
\end{equation*}
By letting $\tau\to-\infty$, one easily checks that the second part tends to 0, so $\Itilde_+^\alpha
f(t+1)=\Itilde_+^\alpha f(t)$.

The following example explains the action of $\Itilde_+^\alpha$ on trigonometric functions.

\begin{example}
Let us compute $\Itilde_+^\alpha$ on the family of complex functions $\phi_r(t)=e^{irt}-1$ for
$r>0$. Suppose $0<\alpha<1$. The formula
\begin{equation*}
   \Gamma(\alpha)=\int_0^\infty s^{\alpha-1}e^{-s}ds=u^\alpha\int_0^\infty s^{\alpha-1}e^{-us}ds
\end{equation*}
is valid for $u>0$ and can be extended to complex numbers with positive real part. One can also
write it for $u=\mp ir$, $r>0$, and we obtain
\begin{equation}\label{trigo}
    \int_0^\infty s^{\alpha-1}e^{\pm irs}ds=e^{\pm i\alpha\pi/2}r^{-\alpha}\Gamma(\alpha)
\end{equation}
where the integral is only semi-convergent. Thus we obtain the classical formula (see Section 7 of
\cite{samko:kil:mari93})
\begin{align*}
    I_+^\alpha e^{irt}
    &=\frac1{\Gamma(\alpha)}\int_{-\infty}^t(t-s)^{\alpha-1}e^{irs}ds\\
    &=\frac{e^{irt}}{\Gamma(\alpha)}\int_0^\infty s^{\alpha-1}e^{-irs}ds
    =r^{-\alpha}e^{-i\alpha\pi/2}e^{irt}.
\end{align*}
We deduce that $\Itilde_+^\alpha\phi_r=r^{-\alpha}e^{-i\alpha\pi/2}\phi_r$, and this relation is
extended to negative $\alpha$ since the operators of exponents $\alpha$ and $-\alpha$ are the
inverse of each other (Theorem \ref{itildeth}). In particular,
\begin{equation}\label{ihattrigo}
\begin{split}
   \Itilde_+^\alpha\bigl(1-\cos(rt)\bigr)&=
   r^{-\alpha}\Bigl(\cos(\alpha\pi/2)-\cos(rt-\alpha\pi/2)\Bigr)\\
   \Itilde_+^\alpha\sin(rt)&=
   r^{-\alpha}\Bigl(\sin(\alpha\pi/2)+\sin(rt-\alpha\pi/2)\Bigr).
\end{split}
\end{equation}
\end{example}

\begin{remark}
We can similarly study $\Itilde_-^\alpha$ which multiplies $\phi_r$ by
$r^{-\alpha}e^{i\alpha\pi/2}$; consequently, the two-sided operator
$(\Itilde_+^\alpha+\Itilde_-^\alpha)/(2\cos(\alpha\pi/2))$ multiplies $\phi_r$ by $r^{-\alpha}$.
\end{remark}

Let us now define two modifications $\Ihatp^\alpha$ and $\Ibarp^\alpha$ of $\Itilde_+^\alpha$
which will be useful for the study of the fractional Brownian motion on $[0,1]$. Consider a
bounded function $f$ defined on the time interval $[0,1]$ and such that $f(0)=0$. If $\alpha<0$,
suppose again that $f$ is in $\holder^\beta$ for some $\beta>-\alpha$. Let $g(t)$ be the
1-periodic function coinciding on $[0,1]$ with $f(t)-t\,f(1)$. We now define on $[0,1]$
\begin{equation}\label{ihatdef}
    \Ihatp^\alpha f(t)=t\,f(1)+\Itilde_+^\alpha g(t).
\end{equation}
Thus $\Ihatp^\alpha$ satisfies the formulae \eqref{ihattrigo} for $r=2n\pi$, and we decide
arbitrarily that $\Ihatp^\alpha t=t$. On the other hand, let $h$ be the function with
1-antiperiodic increments, so that
\begin{equation*}
  h(1+t)-h(1+s)=-h(t)+h(s),
\end{equation*}
and coinciding with $f$ on $[0,1]$. We define
\begin{equation}\label{iflatdef}
  \Ibarp^\alpha f(t)=\Itilde_+^\alpha h(t).
\end{equation}
Then $\Ibarp^\alpha$ satisfies \eqref{ihattrigo} for $r=(2n+1)\pi$.

It is clear that $\Ihatp^{\alpha_2}\Ihatp^{\alpha_1}=\Ihatp^{\alpha_1+\alpha_2}$ is satisfied on
$\holder^\beta$ as soon as $\beta$, $\beta+\alpha_1$ and $\beta+\alpha_1+\alpha_2$ are in $(0,1)$,
and the same property is valid for $\Ibarp^\alpha$ (apply Theorem \ref{itildeth}). Actually, these
composition rules can be used to extend the two operators to arbitrarily large values of $\alpha$.
Moreover, $\Ihatp^\alpha$ and $\Ibarp^\alpha$ are homeomorphisms from $\holder^\beta$ onto
$\holder^{\beta+\alpha}$ if $\beta$ and $\beta+\alpha$ are in $(0,1)$, and their inverse maps are
$\Ihatp^{-\alpha}$ and $\Ibarp^{-\alpha}$.

\subsection{Some other operators}

Let us describe the other operators which are used in this work. The multiplication operator
$\Pi^\alpha$, $\alpha\in\reel$, has already been defined in \eqref{pialphadef} on $\reel_+^\star$,
and let us complement it with
\begin{equation}\label{pitildedef}
   \Pitilde^\alpha f(t)=I_{0+}^1\Pi^\alpha D^1 f(t)
   =\int_0^ts^\alpha df(s)
   =t^\alpha f(t)-\alpha\int_0^ts^{\alpha-1}f(s)ds
\end{equation}
for $f$ smooth with compact support. In the last form, we see that $\Pitilde^\alpha f$ can be
defined as soon as $t^{\alpha-1}f(t)$ is integrable on any $[0,T]$, so on
$\holder^{\beta,\gamma,\delta}$ if $\alpha+\beta+\gamma>0$.

On the other hand, let us define for $\alpha\in\reel$ the time inversion operators $T_\alpha$ and
$T_\alpha'$ on $\reel_+^\star$ by
\begin{equation}\label{talpha}
   T_\alpha f(t)=t^{2\alpha}f(1/t)
\end{equation}
and
\begin{align}
   T_\alpha'f(t)&=-I_{0+}^1T_{\alpha-1}D^1f(t)=-\int_0^ts^{2\alpha-2}D^1f(1/s)ds
   =-\int_{1/t}^\infty s^{-2\alpha}df(s)\notag\\
   &=t^{2\alpha}f(1/t)-2\alpha\int_{1/t}^\infty s^{-2\alpha-1}f(s)ds\label{tpalpha}
\end{align}
and the last form can be used if $t^{-2\alpha-1}f(t)$ is integrable on any $[T,\infty)$, so in
particular on $\holder^{\beta,\gamma,\delta}$ if $2\alpha>\beta+\delta$. Actually, the form of
$T_\alpha$ and a comparison of \eqref{tpalpha} and \eqref{pitildedef} show that
\begin{equation}\label{tppit}
   T_\alpha=\Pi^{2\alpha}T_0=T_0\Pi^{-2\alpha},\qquad T_\alpha'=\Pitilde^{2\alpha}T_0.
\end{equation}
Notice that $T_\alpha$ and $T_\alpha'$ are involutions, so that
\begin{equation}\label{taltpal}
  T_\alpha T_\alpha'f(t)=f(t)-2\alpha\,t^{2\alpha}\int_t^\infty s^{-2\alpha-1}f(s)ds
\end{equation}
and
\begin{equation}\label{tpaltal}
   T_\alpha'T_\alpha f(t)=\Pitilde^{2\alpha}\Pi^{-2\alpha}f(t)=f(t)-2\alpha\int_0^tf(s)\frac{ds}{s}
\end{equation}
are the inverse transformation of each other.

\begin{theorem}\label{otherth}
Let $0<\beta<1$ and consider the time interval $\reel_+^\star$.
\begin{itemize}
\item The operator $\Pitilde^\alpha$ maps continuously $\holder^{\beta,\gamma,\delta}$ into
    $\holder^{\beta,\gamma+\alpha,\delta+\alpha}$ if $\beta+\gamma+\alpha>0$ and
    $\beta+\delta+\alpha>0$. It satisfies the composition rule
    $\Pitilde^{\alpha_2}\Pitilde^{\alpha_1}=\Pitilde^{\alpha_1+\alpha_2}$ on
    $\holder^{\beta,\gamma,\delta}$ if $\beta+\gamma+\alpha_1>0$ and $\beta+\gamma+\alpha_1+\alpha_2>0$.
\item The operator $T_\alpha$ maps continuously $\holder^{\beta,\gamma,\delta}$ into
    $\holder^{\beta,-\delta+2(\alpha-\beta),-\gamma+2(\alpha-\beta)}$. If moreover
    $2\alpha>\beta+\delta$ and $2\alpha>\beta+\gamma$, the operator $T_\alpha'$ satisfies the
    same property.
\end{itemize}
\end{theorem}

\begin{proof}
We prove separately the two parts.

\medskip\noindent\emph{Study of $\Pitilde^\alpha$.}
The continuity on $\holder^{\beta,\gamma,\delta}$ is proved by noticing
\begin{equation*}
   \bigl|\Pitilde^\alpha f(t)\bigr|\le\bigl|\Pi^\alpha f(t)\bigr|+C\int_0^ts^{\alpha-1+\beta}
   s^{\gamma,\delta}ds\le C't^{\alpha+\beta}t^{\gamma,\delta},
\end{equation*}
\begin{align*}
   \bigl|\Pitilde^\alpha f(t)-\Pitilde^\alpha f(s)\bigr|
   &\le\bigl|\Pi^\alpha f(t)-\Pi^\alpha f(s)\bigr|+C\int_s^tu^{\alpha+\beta-1}u^{\gamma,\delta}du\\
   &\le\bigl|\Pi^\alpha f(t)-\Pi^\alpha f(s)\bigr|+C'(t-s)^\beta
   \sup_{s\le u\le t}u^\alpha u^{\gamma,\delta},
\end{align*}
and by applying Theorem \ref{pialphath}. The composition rule is evident for smooth functions (use
the first equality of \eqref{pitildedef}), and can be extended by density (the parameter $\delta$
is unimportant since we only need the functions on bounded time intervals).

\medskip\noindent\emph{Study of $T_\alpha$ and $T_\alpha'$}.
If $f$ is in $\holder^{\beta,\gamma,\delta}$, then $f(1/t)$ is dominated by
$t^{-\beta}t^{-\delta,-\gamma}$, and if $2^n\le s\le t\le2^{n+1}$,
\begin{align*}
   \bigl|f(1/t)-f(1/s)\bigr|&\le C\sup_{s\le u\le t}u^{-\delta,-\gamma}\bigl(1/s-1/t\bigr)^\beta
   \le C'(2^n)^{-\delta,-\gamma}s^{-\beta}t^{-\beta}(t-s)^\beta\\
   &\le C''(2^n)^{-\delta-2\beta,-\gamma-2\beta}(t-s)^\beta,
\end{align*}
so $T_0f:t\mapsto f(1/t)$ is in $\holder^{\beta,-2\beta-\delta,-2\beta-\gamma}$. The continuity of
$T_\alpha$ and $T_\alpha'$ is then a consequence of \eqref{tppit} and of the continuity of
$\Pi^{2\alpha}$ and $\Pitilde^{2\alpha}$. \qed
\end{proof}

\begin{remark}\label{ttpcont}
We deduce in particular from Theorem \ref{otherth} that $T_\alpha$ and $T_\alpha'$ are
homeomorphisms from $\holder^{\alpha-,0,0+}$ into itself for $0<\alpha<1$. We also deduce that
$T_\alpha T_\alpha'$, respectively $T_\alpha'T_\alpha$, is a continuous endomorphism of
$\holder^{\beta,\gamma,\delta}$ when $2\alpha>\beta+\gamma$ and $2\alpha>\beta+\delta$,
respectively when $\beta+\gamma>0$ and $\beta+\delta>0$; when the four conditions are satisfied,
they are the inverse of each other. The form $\Pitilde^{2\alpha}\Pi^{-2\alpha}$ of
$T_\alpha'T_\alpha$ can be used on a bounded time interval $[0,T]$, and in this case we only need
$\beta+\gamma>0$.
\end{remark}

The time inversion operator $T_0$ enables to write the relationship between $I_-^\alpha$ and
$I_{0+}^\alpha$ on $\reel_+^\star$. If $\alpha>0$ and if $f$ is a smooth function with compact
support in $\reel_+^\star$, we deduce from the change of variables $s\mapsto1/s$ that
\begin{equation*}
    I_-^\alpha f(1/t)
    =\int_{1/t}^\infty\Bigl(s-\frac1t\Bigr)^{\alpha-1}f(s)ds
    =\int_0^t\Bigl(\frac1s-\frac1t\Bigr)^{\alpha-1}f(1/s)\frac{ds}{s^2}
\end{equation*}
so that
\begin{equation}\label{t0ialpha}
    T_0I_-^\alpha T_0=\Pi^{1-\alpha}I_{0+}^\alpha\Pi^{-1-\alpha}.
\end{equation}

\section{Time inversion for self-similar processes}
\label{timeinvers}

We give here time inversion properties which are valid for any $H$-self-similar centred Gaussian
process $(\Xi_t;t>0)$, and not only for the fractional Brownian motion. Such a process must have a
covariance kernel of the form
\begin{equation}\label{cst}
   C(s,t)=s^Ht^H\rho(s/t)
\end{equation}
where $\rho(u)=\rho(1/u)$ and $|\rho(u)|\le\rho(1)$. It then follows immediately by comparing the
covariance kernels that if $T_H$ is the time inversion operator defined in \eqref{talpha}, then one
has the equality in law $T_H\Xi\simeq\Xi$. Notice that this holds even when $H$ is not positive.

\begin{remark}
The Lamperti transform (see for instance \cite{cherid:kawa:mae03})
\begin{equation}\label{lamper}
\bigl(\Xi(t);\;t>0\bigr)\mapsto\bigl(e^{-Ht}\Xi(e^t);\;t\in\reel\bigr)
\end{equation}
maps $H$-self-similar processes $\Xi_t$ into stationary processes $Z_t$. Then $T_H\Xi\simeq\Xi$ is
equivalent to the property $Z_{-t}\simeq Z_t$ which is valid for stationary Gaussian processes
(invariance by time reversal).
\end{remark}

\begin{remark}\label{hhmz}
We have $T_HB^H\simeq B^H$ and can deduce properties of $B^H$ on $[1,+\infty)$ from its properties
on $[0,1]$. For instance, $B^H$ lives in $\holder^{H-}$ on $[0,1]$, and we can check from Theorem
\ref{otherth} that $T_H$ sends this space on $[0,1]$ into the space $\holder^{H-,.,0+}$ on
$[1,+\infty)$; thus $B^H$ lives in $\holder^{H-,0,0+}$ on $\reel_+$ (notation \eqref{hbmoins}).
\end{remark}

We now prove another time inversion property when $H>0$ (we do not assume $H<1$). Assume
provisionally that the paths of $\Xi$ are absolutely continuous; then its derivative $D^1\Xi$ is
$(H-1)$-self-similar, so $T_{H-1}D^1\Xi\simeq D^1\Xi$ and
\begin{equation*}
   T_H'\Xi=-I_{0+}^1T_{H-1}D^1\Xi\simeq-I_{0+}^1D^1\Xi=-\Xi\simeq\Xi.
\end{equation*}
In the general case (when $\Xi$ is not absolutely continuous), the same property can be proved with
the theory of generalised processes (as said in \cite{molchan03}); we here avoid this theory.

\begin{theorem}\label{tinvth}
For $H>0$, let $(\Xi_t;\;t\ge0)$ be a $H$-self-similar centred Gaussian process, and consider the
time inversion operators $T_H$ and $T_H'$. Then one has the equalities in law $T_H'\Xi\simeq
T_H\Xi\simeq\Xi$.
\end{theorem}

\begin{proof}
As it has already been said in the beginning of this Subsection, $T_H\Xi\simeq\Xi$ is obtained by
comparing the covariance kernels. Since $\Xi$ is $H$-self-similar, the norm of $\Xi_t$ in
$L^1(\Omega)$ is proportional to $t^H$, so the variable $\int_T^\infty|\Xi_t|t^{-2H-1}dt$ is in
$L^1(\Omega)$ for any $T>0$, and is therefore almost surely finite. Thus $T_H'\Xi$ is well defined.
Moreover, $T_H'\Xi=T_H'T_HT_H\Xi\simeq T_H'T_H\Xi$, so let us compare the covariance kernels of
$\Xi$ and $T_H'T_H\Xi=\Pitilde^{2H}\Pi^{-2H}\Xi$ given by \eqref{tpaltal}. We have from \eqref{cst}
that
\begin{equation*}
   \esp\Bigl[\Xi_T\int_0^S\Xi_s\frac{ds}{s}\Bigr]
   =T^H\int_0^Ss^{H-1}\rho(s/T)ds=T^{2H}\int_0^{S/T}u^{H-1}\rho(u)du.
\end{equation*}
Thus
\begin{align*}
   \esp\Bigl[&\Bigl(\int_0^T\Xi_t\frac{dt}{t}\Bigr)\Bigl(\int_0^S\Xi_s\frac{ds}{s}\Bigr)\Bigr]\\
   &=\int_0^Tt^{2H-1}\int_0^{S/t}u^{H-1}\rho(u)du\,dt
   =\frac1{2H}\int_0^\infty\bigl(T\wedge\frac Su\bigr)^{2H}u^{H-1}\rho(u)du\\
   &=\frac1{2H}\Bigl(T^{2H}\int_0^{S/T}u^{H-1}\rho(u)du+S^{2H}\int_{S/T}^\infty u^{-H-1}\rho(u)du\Bigr)\\
   &=\frac1{2H}\Bigl(T^{2H}\int_0^{S/T}u^{H-1}\rho(u)du+S^{2H}\int_0^{T/S}u^{H-1}\rho(u)du\Bigr)
\end{align*}
(we used $\rho(1/u)=\rho(u)$ in the last equality). We deduce from these two equations that
\begin{equation*}
   \esp\Bigl[\Bigl(\Xi_T-2H\int_0^T\Xi_t\frac{dt}{t}\Bigr)\Bigl(\Xi_S-2H\int_0^S\Xi_s\frac{ds}{s}\Bigr)\Bigr]
   =\esp\bigl[\Xi_T\Xi_S\bigr]
\end{equation*}
since the other terms cancel one another, so $T_H'T_H\Xi$ has the same covariance kernel as $\Xi$.
\qed
\end{proof}

\begin{remark}
Theorem \ref{tinvth} can be applied to the fractional Brownian motion $B^H$. Moreover, the
relations $B^H\simeq T_HB^H\simeq T_H'B^H$ can be extended to $\reel^\star$ by defining
\begin{equation*}
   T_Hf(t)=|t|^{2H}f(1/t),\qquad T_H'f=\mp I_{0\pm}^1T_{H-1}D^1f\quad\text{on $\reel_\pm^\star$.}
\end{equation*}
Since $B^H$ also has stationary increments, we can deduce how the law of the generalised process
$D^1B^H$ is transformed under the time transformations $t\mapsto(at+b)/(ct+d)$, see
\cite{molchan03}.
\end{remark}

The law of the $H$-self-similar process $\Xi$ is therefore invariant by the transformations
$T_HT_H'$ and $T_H'T_H=\Pitilde^{2H}\Pi^{-2H}$ given by \eqref{taltpal} and \eqref{tpaltal}. We now
introduce a generalisation $T_{H,L}$ of $T_H'T_H$, which was also studied in \cite{cjost07}.

\begin{theorem}\label{thlth}
On the time interval $\reel_+$,  for $H>0$ and $L>0$, the operator
\begin{equation}\label{thldef}
 T_{H,L}=\Pi^{H-L}T_L'T_L\Pi^{L-H}=\Pi^{H-L}\Pitilde^{2L}\Pi^{-L-H}
\end{equation}
is a continuous endomorphism of $\holder^{\beta,\gamma,\delta}$ when $0<\beta<1$, and
$\beta+\gamma$ and $\beta+\delta$ are greater than $H-L$; in particular, it is a continuous
endomorphism of $\holder^{H-,0,0+}$ if $0<H<1$. It is defined on a function $f$ as soon as
$t^{L-H-1}f(t)$ is integrable on any $[0,T]$, and it satisfies
\begin{equation}\label{thlft}
 T_{H,L}f(t)=f(t)-2L\,t^{H-L}\int_0^tf(s)s^{L-H-1}ds.
\end{equation}
If $\Xi$ is a $H$-self-similar centred Gaussian process, then $T_{H,L}\Xi$ has the same law as
$\Xi$.
\end{theorem}

\begin{proof}
The continuity property of $T_{H,L}$ can be deduced from Theorem \ref{otherth} and Remark
\ref{ttpcont}. The representation \eqref{thlft} follows easily from \eqref{pitildedef} and the
second form of $T_{H,L}$ in \eqref{thldef}. Let $\Xi$ be a centred Gaussian $H$-self-similar
process; then the $L^1$-norm of $\Xi_t$ is proportional to $t^H$, so $\int_0^Tt^{L-H-1}|\Xi_t|dt$
is integrable and therefore almost surely finite for any $T>0$. We deduce that $T_{H,L}\Xi$ is well
defined; we have
\begin{equation*}
   T_L'T_L\Pi^{L-H}\Xi\simeq\Pi^{L-H}\Xi
\end{equation*}
because $\Pi^{L-H}\Xi$ is $L$-self-similar. By applying $\Pi^{H-L}$ to both sides we obtain
$T_{H,L}\Xi\simeq\Xi$. \qed
\end{proof}

\begin{remark}
In the non centred case, we have $T_H\Xi\simeq\Xi$ and $T_H'\Xi\simeq T_{H,L}\Xi\simeq-\Xi$.
\end{remark}

We will resume our study of $T_{H,L}$ for self-similar processes in Subsection \ref{noncanon}.

\section{Representations of fractional Brownian motions}
\label{represent}

Starting from the classical representation of fractional Brownian motions on $\reel$ described in
Subsection \ref{repreel}, we study canonical representations on $\reel_+$ (Subsection
\ref{repplus}) and $\reel_-$ (Subsection \ref{repmoins}). In Subsection \ref{noncanon}, we also
consider the non canonical representations on $\reel_+$ introduced in Theorem \ref{thlth}.

\subsection{A representation on $\reel$}\label{repreel}

For $0<H<1$, the basic representation of a fractional Brownian motion $B^H$ is
\begin{equation}\label{mandkappa}
   B_t^H=\kappa\int_{-\infty}^{+\infty}\Bigl((t-s)_+^{H-1/2}-(-s)_+^{H-1/2}
   \Bigr)dW_s
\end{equation}
for a positive parameter $\kappa$, see \cite{mandel:ness68}. It is not difficult to check that the
integral of the right-hand side is Gaussian, centred, with stationary increments, and
$H$-self-similar. Thus $B_t^H$ is a fractional Brownian motion; its covariance is given by
\eqref{bhcov}, and the variance $\rho$ of $B_1^H$ is proportional to $\kappa^2$; the precise
relationship between $\rho$ and $\kappa$ is given in Theorem \ref{varbh}. Subsequently, we will
consider the fractional Brownian motion corresponding to
\begin{equation}\label{kappah}
   \kappa=\kappa(H)=1/\Gamma(H+1/2),
\end{equation}
so that (following \eqref{rhokcos})
\begin{equation}\label{rhoh}
   \rho=\rho(H)=-2\frac{\cos(\pi H)}{\pi}\Gamma(-2H),\qquad\rho(1/2)=1.
\end{equation}
In particular, $B^{1/2}=W$ is the standard Brownian motion. This choice of $\kappa$ is due to the
following result, where we use the modified Riemann-Liouville operators of Subsection \ref{infhor}.

\begin{theorem}\label{reprth}
The family of processes $(B^H;\;0<H<1)$ defined by \eqref{mandkappa} with \eqref{kappah} can be
written as
\begin{equation}\label{bhitildw}
   B^H=\Itilde_+^{H-1/2}W.
\end{equation}
More generally,
\begin{equation}\label{bhitildbj}
   B^H=\Itilde_+^{H-J}B^J
\end{equation}
for any $0<J,H<1$.
\end{theorem}

\begin{proof}
The formula \eqref{bhitildw} would hold true from \eqref{rieminf} if $W$ were Lipschitz with
compact support; the operator $\Itilde_+^{H-1/2}$ is continuous on $\holder^{1/2-,0,0+}$ (Theorem
\ref{itildeth}) in which $W$ lives, and Lipschitz functions with compact support are dense in this
space (Theorem \ref{dense}); moreover, integration by parts shows that the stochastic integral in
the right-hand side of \eqref{mandkappa} can also be computed by approximating $W$ with smooth
functions with compact support, so \eqref{bhitildw} holds almost surely. Then \eqref{bhitildbj}
follows from the composition rules for Riemann-Liouville operators (Theorem \ref{itildeth}). \qed
\end{proof}

We deduce in particular from \eqref{bhitildbj} that \eqref{bhitildw} can be reversed
($W=B^{1/2}$), and
\begin{equation*}
   W=\Itilde_+^{1/2-H}B^H.
\end{equation*}
Thus the increments of $W$ and $B^H$ generate the same completed filtration, namely
$\tribuf_t(dB^H)=\tribuf_t(dW)$ (with notation \eqref{tribu}).

\begin{remark}
Relation \eqref{bhitildbj} can be written by means of \eqref{itildedef} ($H>J$) or \eqref{martilde}
($H<J$). It can be written more informally as
\begin{equation*}
   B_t^H=\frac1{\Gamma(H-J+1)}\int_{-\infty}^{+\infty}\Bigl((t-s)_+^{H-J}-(-s)_+^{H-J}
   \Bigr)dB_s^J,
\end{equation*}
where the integral is obtained by approximating $B^J$ by Lipschitz functions with compact support,
and passing to the limit.
\end{remark}

Relations \eqref{bhitildw} or \eqref{bhitildbj} can be restricted to the time interval $\reel_-$;
in order to know $B^H$ on $\reel_-$, we only need $W$ on $\reel_-$, and vice-versa. On the other
hand, they cannot be used on $\reel_+$; in order to know $B^H$ on $\reel_+$, we have to know $W$ on
the whole real line $\reel$. If we want a representation on $\reel_+$, we can reverse the time
($t\mapsto-t$) for all the processes, so that the operators $\Itilde_+$ are replaced by
$\Itilde_-$. We obtain on $\reel_+$ the backward representation
\begin{equation}\label{mandback}
    B_t^H=\Itilde_-^{H-1/2}W(t)
    =\frac1{\Gamma(H+1/2)}\int_0^\infty\Bigl(s^{H-1/2}-(s-t)_+^{H-1/2}\Bigr)dW_s.
\end{equation}
However, in this formula, if we want to know $B^H$ at a single time $t$, we need $W$ on the
whole half-line $\reel_+$; next subsection is devoted to a representation formula where we only
need $W$ on $[0,t]$.

\subsection{Canonical representation on $\reel_+$}\label{repplus}

We shall here explain the derivation of the canonical representation of fractional Brownian motions
on $\reel_+$ which was found by \cite{molchan:golo69,molchan03}, and the general relationship
between $B^J$ and $B^H$ which was given in \cite{cjost06}. More precisely, we want the various
processes $(B^H;0<H<1)$ to be deduced from one another, so that all of them generate the same
filtration.

As explained in the introduction, we start from the relation $B^H=\Itilde_-^{H-1/2}W$ of
\eqref{mandback} and apply the time inversion $t\mapsto1/t$ on the increments $dW_t$ and
$dB_t^H$; this time inversion is made by means of the operators $T_{1/2}'$ and $T_H'$ defined
in \eqref{talpha} (they are involutions), which preserve respectively the laws of $W$ and $B^H$
(Theorem \ref{tinvth}). Thus
\begin{equation*}
   B^H\simeq\bigl(T_H'\Itilde_-^{H-1/2}T_{1/2}'\bigr)W.
\end{equation*}
It appears that this is the canonical representation of $B^H$. We now make more explicit this
calculation, and generalise it to the comparison of $B^H$ and $B^J$ for any $J$ and $H$; starting
from $B^H=\Itilde_-^{H-J}B^J$, we can show similarly that
\begin{equation}\label{bhsimeq}
   B^H\simeq\bigl(T_H'\Itilde_-^{H-J}T_J'\bigr)B^J.
\end{equation}

\begin{theorem}\label{molgolo}
On the time interval $\reel_+$, the family of fractional Brownian motions $B^H$, $0<H<1$, can be
defined jointly so that $B^H=G_{0+}^{J,H}B^J$ for
\begin{equation}\label{gpiipi}
    G_{0+}^{J,H}=\Pitilde^{H+J-1}I_{0+}^{H-J}\Pitilde^{1-H-J}
\end{equation}
(see Section \ref{defnot} for the definitions of $I_{0+}^\alpha$ and $\Pitilde^\alpha$). This
family of operators satisfies the composition rule $G_{0+}^{H,L}G_{0+}^{J,H}=G_{0+}^{J,L}$, and all
the processes $B^H$ generate the same completed filtration. Moreover, the operator $G_{0+}^{J,H}$
maps continuously $\holder^{J-,0,0+}$ (where paths of $B^J$ live) into $\holder^{H-,0,0+}$, and can
be defined by the following relation; if we define
\begin{equation}\label{phijhdef}
    \phi^{J,H}(u)=(H-J)\int_1^u\Bigl(v^{H+J-1}-1\Bigr)(v-1)^{H-J-1}dv+(u-1)^{H-J}
\end{equation}
for $0<J,H<1$ and $u>1$, and if
\begin{equation}\label{k0pjh}
   K_{0+}^{J,H}(t,s)=\frac1{\Gamma(H-J+1)}\phi^{J,H}\Bigl(\frac ts\Bigr)s^{H-J},
\end{equation}
then
\begin{equation}\label{gjhkjh}
   G_{0+}^{J,H}f(t)=\int_0^tK_{0+}^{J,H}(t,s)df(s)
\end{equation}
for $f$ Lipschitz with compact support in $\reel_+^\star$. Moreover, $B^H$ is given by the It\^o
integral
\begin{equation}\label{bhkw}
 B_t^H=\int_0^tK_{0+}^{1/2,H}(t,s)dW_s
\end{equation}
for $W=B^{1/2}$.

\end{theorem}

\begin{proof}
Let us divide the proof into four steps.

\medskip\noindent\emph{Step 1: Definition of the families $G_{0+}^{J,H}$ and $B^H$.} Following
\eqref{bhsimeq}, we define
\begin{equation}\label{gpjh}
 G_{0+}^{J,H}=T_H'\Itilde_-^{H-J}T_J',\quad
 B^H=G_{0+}^{1/2,H}W,
\end{equation}
so that $B^H$ is a $H$-fractional Brownian motion. The continuity of $G_{0+}^{J,H}$ from
$\holder^{J-,0,0+}$ into $\holder^{H-,0,0+}$ is then a consequence of Theorems \ref{itildeth} and
\ref{otherth}; it indeed follows from these two theorems that $T_J'$ and $T_H'$ are continuous
endomorphisms of respectively $\holder^{J-,0,0+}$ and $\holder^{H-,0,0+}$, and that
$\Itilde_-^{H-J}$ is continuous from $\holder^{J-,0,0+}$ into $\holder^{H-,0,0+}$. Moreover
\begin{equation*}
    G_{0+}^{H,L}G_{0+}^{J,H}=T_L'\Itilde_-^{L-H}T_H'T_H'\Itilde_-^{H-J}T_J'
    =T_L'\Itilde_-^{L-H}\Itilde_-^{H-J}T_J'=T_L'\Itilde_-^{L-J}T_J'=G_{0+}^{J,L}
\end{equation*}
and consequently
\begin{equation*}
    G_{0+}^{J,H}B^J=G_{0+}^{J,H}G_{0+}^{1/2,J}W=G_{0+}^{1/2,H}W=B^H.
\end{equation*}
The equality between filtrations of $B^H$ also follows from this relation.

\medskip\noindent\emph{Step 2: Proof of \eqref{gpiipi}.} First assume $H>J$, and let us work on smooth
functions with compact support in $\reel_+^\star$. We deduce from \eqref{t0ialpha} and the
relations $T_\alpha=\Pi^{2\alpha}T_0=T_0\Pi^{-2\alpha}$ that
\begin{align}
    T_{H-1}I_-^{H-J}T_{J-1}&=\Pi^{2H-2}T_0I_-^{H-J}T_0\Pi^{2-2J}
    =\Pi^{2H-2}\Pi^{1-H+J}I_{0+}^{H-J}\Pi^{-1-H+J}\Pi^{2-2J}\notag\\
    &=\Pi^{H+J-1}I_{0+}^{H-J}\Pi^{1-H-J}.\label{th1i}
\end{align}
On the other hand, $T_\alpha'$ has been defined as $-I_{0+}^1T_{\alpha-1}D^1$, and
$\Itilde_-^\alpha=I_{0+}^1I_-^\alpha D^1$ from \eqref{ditil}, so the definition \eqref{gpjh} can
be written as
\begin{align}
   G_{0+}^{J,H}&=(I_{0+}^1T_{H-1}D^1)(I_{0+}^1I_-^{H-J}D^1)(I_{0+}^1T_{J-1}D^1)\notag\\
   &=I_{0+}^1T_{H-1}I_-^{H-J}T_{J-1}D^1\notag\\
   &= I_{0+}^1\Pi^{H+J-1}I_{0+}^{H-J}\Pi^{1-H-J}D^1\label{gipiipid}\\
   &=I_{0+}^1\Pi^{H+J-1}I_{0+}^{H-J}D^1I_{0+}^1\Pi^{1-H-J}D^1\notag\\
   &=\bigl(I_{0+}^1\Pi^{H+J-1}D^1\bigr)I_{0+}^{H-J}\bigl(I_{0+}^1\Pi^{1-H-J}D^1\bigr)\notag\\
   &=\Pitilde^{H+J-1}I_{0+}^{H-J}\Pitilde^{1-H-J}\notag
\end{align}
(we used \eqref{th1i} in the third equality and Theorem \ref{smooth} in the fifth one). The
equality can be extended to the functional space $\holder^{J-,0,0+}$, since $G_{0+}^{J,H}$ is
continuous on this space, and the right-hand side is continuous on $\holder^{J-}$ on any interval
$[0,T]$. Moreover, inverting this relation provides $G_{0+}^{H,J}$, so that this expression of
$G_{0+}^{J,H}$ also holds when $H<J$.

\medskip\noindent\emph{Step 3: Proof of \eqref{gjhkjh}.} For smooth functions $f$ with compact support
in $\reel_+^\star$, \eqref{pimgam} yields
\begin{align*}
    \Pi^{H+J-1}&I_{0+}^{H-J}\Pi^{1-H-J}f(t)\\&=I_{0+}^{H-J}f(t)+\frac1{\Gamma(H-J)}\int_0^t
    \Bigl(\bigl(\frac st\bigr)^{1-H-J}-1\Bigr)(t-s)^{H-J-1}f(s)ds,
\end{align*}
so \eqref{gipiipid} implies
\begin{align*}
    G_{0+}^{J,H}f(t)=&I_{0+}^{H-J}f(t)+\frac1{\Gamma(H-J)}
    \int_0^t\biggl(\int_0^v\Bigl(\bigl(\frac sv\bigr)^{1-H-J}-1\Bigr)(v-s)^{H-J-1}df(s)\biggr)dv\\
    =&\frac1{\Gamma(H-J+1)}\int_0^t(t-s)^{H-J}df(s)\\&+\frac{H-J}{\Gamma(H-J+1)}
    \int_0^t\biggl(\int_s^t\Bigl(\bigl(\frac sv\bigr)^{1-H-J}-1\Bigr)(v-s)^{H-J-1}dv\biggr)df(s).
\end{align*}
This expression can be written as \eqref{gjhkjh} for a kernel $K_{0+}^{J,H}$, and a scaling
argument shows that $K_{0+}^{J,H}$ is of the form \eqref{k0pjh} for
$\phi^{J,H}(u)=\Gamma(H-J+1)K_{0+}^{J,H}(u,1)$. Then \eqref{phijhdef} follows from a simple
verification.

\medskip\noindent\emph{Step 4: Proof of \eqref{bhkw}.} By means of an integration by parts, we write
\eqref{gjhkjh} for $J=1/2$ and $H\ne1/2$ in the form
\begin{align}
   G_{0+}^{1/2,H}f(t)
   &=\frac{f(t)}{t}\int_0^tK_{0+}^{1/2,H}(t,s)ds+\int_0^tK_{0+}^{1/2,H}(t,s)\bigl(D^1f(s)-f(t)/t\bigr)ds\notag\\
   &=\frac{f(t)}{t}\int_0^tK_{0+}^{1/2,H}(t,s)ds
   -\int_0^t\bigl(f(s)-\frac stf(t)\bigr)\partial_sK_{0+}^{1/2,H}(t,s)ds.\label{g12jh}
\end{align}
On the other hand,
\begin{equation*}
   (\phi^{J,H})'(u)=(H-J)(u-1)^{H-J-1}u^{H+J-1}
\end{equation*}
so that
\begin{equation*}
   \partial_sK_{0+}^{J,H}(t,s)=\frac1{\Gamma(H-J)}\Bigl(\phi^{J,H}\bigl(\frac ts\bigr)s^{H-J-1}
   -(t-s)^{H-J-1}\bigl(\frac ts\bigr)^{H+J}\Bigr).
\end{equation*}
An asymptotic study of \eqref{phijhdef} shows that $\phi^{1/2,H}(u)$ is $O((u-1)^{H-1/2})$ as
$u\downarrow1$ and $O(u^{2H-1}\vee1)$ as $u\uparrow\infty$; thus $\partial_sK_{0+}^{1/2,H}(t,s)$ is
$O((t-s)^{H-3/2})$ as $s\uparrow t$, and is $O(s^{-H-1/2}\vee s^{H-3/2})$ as $s\downarrow0$. An
approximation by smooth functions shows that \eqref{g12jh} is still valid for $W$, and a stochastic
integration by parts leads to \eqref{bhkw}. \qed
\end{proof}

\begin{remark}
It is also possible to write a representation $B^H=G_{T+}^{J,H}B^J$ on the time interval
$[T,+\infty)$, associated to the kernel $K_{T+}^{J,H}(t,s)=K_{0+}^{J,H}(t-T,s-T)$. In
\cite{cjost08}, it is proved that letting $T$ tend to $-\infty$, we recover at the limit
\eqref{mandkappa}.
\end{remark}

\begin{remark}
If $H>J$, we have
\begin{equation*}
    \phi^{J,H}(u)=(H-J)\int_1^uv^{H+J-1}(v-1)^{H-J-1}dv.
\end{equation*}
If $H<J$, this integral diverges and $\phi^{J,H}(u)$ is its principal value. This function, and
therefore the kernel $K_{0+}^{J,H}(t,s)$ can also be written by means of the Gauss hypergeometric
function, see \cite{decreuse:ustu99,cjost06}.
\end{remark}

\begin{remark}\label{hpj1}
If $H+J=1$, then \eqref{gpiipi} is simply written as $G_{0+}^{J,H}=I_{0+}^{H-J}$. Thus the relation
between $B^H$ and $B^{1-H}$ is particularly simple (as it has already been noticed in
\cite{cjost06}), but we have no intuitive explanation of this fact.
\end{remark}

\begin{remark}
The expression \eqref{gpiipi} for $G_{0+}^{J,H}$ is close to the representation given in
\cite{norros:val:vir99} for $J=1/2$. We define
\begin{equation*}
   Z_t^{J,H}=I_{0+}^{H-J}\Pitilde^{1-J-H}B^J(t)
   =\frac1{\Gamma(H-J+1)}\int_0^t(t-s)^{H-J}s^{1-J-H}dB_s^J
\end{equation*}
which is an It\^{o} integral in the case $J=1/2$, and the fractional Brownian motion $B^H$ is given
by
\begin{equation*}
   B_t^H=\Pitilde^{H+J-1}Z^{J,H}(t)=\int_0^ts^{H+J-1}dZ_s^{J,H}
\end{equation*}
which can be defined by integration by parts.
\end{remark}

\begin{remark}
In the case $J=1/2$, let us compare our result with the decomposition of \cite{decreuse:ustu99}. We
look for a decomposition of $G_{0+}^{1/2,H}$ which would be valid on the classical Cameron-Martin
space $\hilbert_{1/2}=I_{0+}^1L^2$ of $W$. To this end, we start from \eqref{gipiipid}
\begin{equation*}
   G_{0+}^{1/2,H}=I_{0+}^1\Pi^{H-1/2}I_{0+}^{H-1/2}\Pi^{1/2-H}D^1
\end{equation*}
which is valid for smooth functions. When $H>1/2$, this formula is valid on $\hilbert_{1/2}$ for
any finite time interval $[0,T]$ because these five operators satisfy the continuity properties
\begin{equation*}
   \hilbert_{1/2}\to L^2\to L^1\to L^1\to L^1\to L^\infty
\end{equation*}
(use the fact that $I_{0+}^\alpha$ is a continuous endomorphism of $L^1$ for $\alpha>0$). However,
it does not make sense on $\hilbert_{1/2}$ for $H<1/2$ because $I_{0+}^{H-1/2}$ is in this case a
fractional derivative, and is not defined for non continuous functions. Thus let us look for an
alternative definition of the operator $G_{0+}^{1/2,H}$; in order to solve this question, we apply
the property \eqref{ipii} of Riemann-Liouville operators and get
\begin{align*}
   G_{0+}^{1/2,H}
   &=I_{0+}^{2H}\bigl(I_{0+}^{1-2H}\Pi^{H-1/2}I_{0+}^{H-1/2}\bigr)\Pi^{1/2-H}D^1\\
   &=I_{0+}^{2H}\bigl(\Pi^{1/2-H}I_{0+}^{1/2-H}\Pi^{2H-1}\bigr)\Pi^{1/2-H}D^1\\
   &=I_{0+}^{2H}\Pi^{1/2-H}I_{0+}^{1/2-H}\Pi^{H-1/2}D^1
\end{align*}
which makes sense on $\hilbert_{1/2}$ if $H<1/2$. This is the expression of \cite{decreuse:ustu99}.
\end{remark}

\begin{remark}
A consequence of \eqref{bhkw} is that we can write the conditional law of $(B_t^H;t\ge S)$ given
$(B_t^H;0\le t\le S)$. This is the prediction problem, see also \cite{gripen:norros96,molchan03}.
\end{remark}

\begin{remark}
Theorem \ref{molgolo} can also be proved by using the time inversion operators $T_H$ rather than
$T_H'$. If we start again from \eqref{mandback} and consider the process with independent
increments
\begin{equation*}
   V_t^H=\int_0^ts^{H-1/2}dW_s,
\end{equation*}
then it appears that $B_t^H$ depends on future values of $V^H$; consequently, $T_HB^H(t)$ depends
on past values of $T_HV^H$. On the other hand, $T_HB^H\simeq B^H$ and $T_HV^H\simeq V^H$ from
Theorem \ref{tinvth}, so we obtain an adapted representation of $B^H$ with respect to $V^H$, and
therefore with respect to $W$. One can verify that this is the same representation as Theorem
\ref{molgolo}; however, the composition rule for the operators $G_{0+}^{J,H}$ is less direct with
this approach.
\end{remark}

Let us give another application of Theorem \ref{molgolo}. The process $B^H$ has stationary
increments, so a natural question is to know whether it can be written as $B_t^H=A_t^H-A_0^H$ for a
stationary centred Gaussian process $A^H$, and to find $A^H$. This is clearly not possible on an
infinite time interval, since the variance of $B^H$ is unbounded. However, let us check that this
is possible in an explicit way on a finite time interval, and that moreover we do not have to
increase the $\sigma$-algebra of $B^H$. Since we are on a bounded time interval $[0,T]$, the
stationarity means that $(A_{U+t}^H;\;0\le t\le T-U)$ and $(A_t^H;\;0\le t\le T-U)$ have the same
law for any $0<U<T$.

\begin{theorem}\label{statiofrac}
Let $T>0$. There exists a stationary centred Gaussian process $(A_t^H;\;0\le t\le T)$ such that
$B_t^H=A_t^H-A_0^H$ is a $H$-fractional Brownian motion on $[0,T]$, and $B^H$ and $A^H$ generate
the same $\sigma$-algebra.
\end{theorem}

\begin{proof}
Consider $B^H=G_{0+}^{1/2,H}W$. We look for a variable $A_0^H$ such that $A_t^H=B_t^H+A_0^H$ is
stationary; this will hold when
\begin{equation*}
   \esp\bigl[A_t^HA_s^H\bigr]=\frac\rho2\bigl(t^{2H}+s^{2H}-|t-s|^{2H}\bigr)
   +\esp\bigl[B_t^HA_0^H\bigr]+\esp\bigl[B_s^HA_0^H\bigr]+\esp\bigl[(A_0^H)^2\bigr]
\end{equation*}
is a function of $t-s$, so when
\begin{equation*}
   \esp\bigl[B_t^HA_0^H\bigr]=-\rho\,t^{2H}/2.
\end{equation*}
By applying the operator $G_{0+}^{H,1/2}$, this condition is shown to be equivalent to
\begin{equation*}
   \esp\bigl[W_tA_0^H\bigr]=-\frac\rho2\,G_{0+}^{H,1/2}t^{2H}
   =-\frac\rho2\,\frac{2H}{H+1/2}\Gamma(H+1/2)t^{H+1/2}
\end{equation*}
by using the formulae \eqref{gipiipid} and \eqref{itbeta} for computing $G_{0+}^{H,1/2}$, and for
$\rho$ given by \eqref{rhoh}. Thus we can choose
\begin{equation*}
   A_0^H=\int_0^T\frac{d}{dt}\esp\bigl[W_tA_0^H\bigr]dW_t
   =-\rho\,H\,\Gamma(H+1/2)\int_0^Tt^{H-1/2}dW_t.
\end{equation*}
\qed
\end{proof}

In particular we have $A_0^{1/2}=-W_T/2$. Of course we can add to $A_0^H$ any independent variable;
this increases the $\sigma$-algebra, but this explains the mutual compatibility of the variables
$A_0^H$ when $T$ increases. More generally, the technique used in the proof enables to write any
variable $A$ of the Gaussian space of $B^H$, knowing the covariances $\esp[A\,B_t^H]$.

\begin{remark}
We can also try to write $B^H$ on $[0,T]$ as the increments of a process which would be stationary
on $\reel$. We shall address this question in Remark \ref{statioproc}.
\end{remark}

\begin{remark}
Another classical stationary process related to the Brownian motion is the Ornstein-Uhlenbeck
process; actually there are two different fractional extensions of this process, see
\cite{cherid:kawa:mae03}.
\end{remark}

\subsection{Canonical representation on $\reel_-$}\label{repmoins}

In the representation \eqref{mandness}, we have $\tribuf_t(dB^H)=\tribuf_t(dW)$ (with notation
\eqref{tribu}). However, when $t<0$, the filtration $\tribuf_t(dB^H)$ is strictly included into
$\tribuf_t(B^H)$. We now give a representation of $B^H$ on the time interval $\reel_-$ for which
$\tribuf_t(B^H)=\tribuf_t(dW)$; one can then deduce a canonical representation of $B^H$ (see Remark
\ref{rmcan} below). In the particular case $H=1/2$ of a standard Brownian motion, we recover the
classical representation of the Brownian bridge.

We want $B_t^H$, $t<0$, to depend on past increments of $W$; by applying the time reversal
$t\mapsto-t$, this is equivalent to wanting $B_t^H$, $t>0$, to depend on future increments of
$W$. The starting point is the operator $T_\alpha T_\alpha'$ of \eqref{taltpal} which can be
written in the form
\begin{equation*}
    T_HT_H'f(t)=-2Ht^{2H}\int_t^\infty s^{-2H-1}\bigl(f(s)-f(t)\bigr)ds.
\end{equation*}
Thus $T_HT_H'f(t)$ depends on future increments of $f$, and the equality in law $B^H\simeq
T_HT_H'B^H$ enables to write $B^H$ as a process depending on future increments of another
$H$-fractional Brownian motion. On the other hand, in the representation
$B^H\simeq\Itilde_-^{H-1/2}W$ of \eqref{mandback}, future increments of $B^H$ depend on
future increments of $W$. Thus, in $B^H\simeq T_HT_H'\Itilde_-^{H-1/2}W$, the value of
$B_t^H$ depends on future increments of $W$, and this answers our question. The same
method can be used with $W$ replaced by $B^J$.

\begin{theorem}\label{rmoins}
Let $B^J$ be a $J$-fractional Brownian motion on $\reel_-$; consider the function $\phi^{J,H}$ of
\eqref{phijhdef}. On $\reel_-^\star$, the operator
\begin{equation*}
   G_+^{J,H}f(t)=\int_{-\infty}^tK_+^{J,H}(t,s)df(s)
\end{equation*}
for $f$ smooth with compact support, with
\begin{equation*}
    K_+^{J,H}(t,s)=\Gamma(H-J+1)^{-1}\phi^{J,H}(s/t)(-t)^{2H}(-s)^{-H-J},
    \qquad s<t<0,
\end{equation*}
can be extended to a continuous operator from $\holder^{J-,0,0+}$ into $\holder^{H-,0,0+}$, and
$\Btilde^{J,H}=G_+^{J,H}B^J$ is a $H$-fractional Brownian motion on $\reel_-$. Moreover,
$\tribuf_t(\Btilde^{J,H})=\tribuf_t(dB^J)$ (with notation \eqref{tribu}).
\end{theorem}

\begin{proof}
We transform the question on $\reel_-$ into a question on $\reel_+$ by means of the time reversal
$t\mapsto-t$. Following the discussion before the theorem, we introduce on $\reel_+^\star$ the
operator
\begin{equation*}
    G_-^{J,H}=T_HT_H'\Itilde_-^{H-J}.
\end{equation*}
It follows from Theorems \ref{itildeth} and \ref{otherth} that $G_-^{J,H}$ maps continuously
$\holder^{J-,0,0+}$ into $\holder^{H-,0,0+}$; moreover $\Btilde^{J,H}=G_-^{J,H}B^J$ is a
$H$-fractional Brownian motion. If we compare $G_-^{J,H}$ with $G_{0+}^{J,H}$ given in
\eqref{gpjh}, we see that
\begin{equation*}
    G_-^{J,H}=T_HG_{0+}^{J,H}T_J'.
\end{equation*}
For $f$ smooth with compact support in $\reel_+^\star$,
\begin{equation*}
    G_{0+}^{J,H}T_J'f(t)=\int_0^tK_{0+}^{J,H}(t,s)s^{2J-2}D^1f(1/s)ds
    =\int_{1/t}^\infty K_{0+}^{J,H}(t,1/s)s^{-2J}df(s)
\end{equation*}
so
\begin{equation*}
    G_-^{J,H}f(t)=t^{2H}\int_t^\infty K_{0+}^{J,H}(1/t,1/s)s^{-2J}df(s)
    =\int_t^\infty K_-^{J,H}(t,s)df(s)
\end{equation*}
with
\begin{equation*}
    K_-^{J,H}(t,s)=t^{2H}s^{-2J}K_{0+}^{J,H}(1/t,1/s)=\Gamma(H-J+1)^{-1}\phi^{J,H}(s/t)t^{2H}s^{-H-J}
\end{equation*}
(apply \eqref{k0pjh}). We still have to check that
\begin{equation*}
   \sigma\bigl(\Btilde_s^{J,H};\;s\ge t\bigr)=\sigma\bigl(B_s^J-B_u^J;\;s\ge u\ge t\bigr)
\end{equation*}
for $t\ge0$. The inclusion of the left-hand side in the right-hand side follows from the discussion
before the theorem. For the inverse inclusion, notice that $\Btilde^{J,H}=G_-^{J,H}B^J$ can be
reversed and
\begin{equation*}
    B^J=\Itilde_-^{J-H}T_H'T_H\Btilde^{J,H}.
\end{equation*}
Thus future increments of $B^J$ depend on future increments of $T_H'T_H\Btilde^{J,H}$, which depend
on future values of $\Btilde^{J,H}$ from \eqref{tpaltal}. \qed
\end{proof}

\begin{remark}\label{rmcan}
The theorem involves $\tribuf_t(dB^J)$ which is strictly smaller than $\tribuf_t(B^J)$, so the
representation is not really canonical on $\reel_-$; however, $\tribuf_t(dB^J)$ is also the
filtration generated by (for instance) the increments of the process
\begin{equation*}
   \Upsilon_t^J=\int_{-\infty}^t(-s)^{-2J}dB_s^J
   =(-t)^{-2J}B_t^J+2J\int_{-\infty}^t(-s)^{-2J-1}B_s^Jds,
\end{equation*}
and
\begin{equation}\label{bups}
   \Btilde_t^{J,H}=\int_{-\infty}^tK_+^{J,H}(t,s)(-s)^{2J}d\Upsilon_s^J.
\end{equation}
The process $\Upsilon_t^J$ tends to 0 at $-\infty$, so
\begin{equation*}
   \tribuf_t(\Btilde^{J,H})=\tribuf_t(dB^J)=\tribuf_t(d\Upsilon^J)=\tribuf_t(\Upsilon^J)
\end{equation*}
and \eqref{bups} is therefore a canonical representation on $\reel_-$ (notice that $\Upsilon^{1/2}$
has independent increments).
\end{remark}

\begin{remark}
By applying Theorem \ref{rmoins} with $J=1/2$, we can predict on $\reel_-$ future values of $B^H$
knowing previous values; this prediction must take into account the fact $B_0^H=0$; this can be
viewed as a bridge; actually for $H=J=1/2$, we recover the classical Brownian bridge. More
precisely, $\phi^{1/2,1/2}\equiv1$, so $K_+^{1/2,1/2}(t,s)=|t|/|s|$ on $\reel_-$; thus $W=B^{1/2}$
and $\Wbar=\Btilde^{1/2,1/2}$ are Brownian motions on $\reel_-$, and satisfy
\begin{equation*}
   \Wbar_t=|t|\int_{-\infty}^t|s|^{-1}dW_s,\qquad
   d\Wbar_t=-\frac{\Wbar_t}{|t|}dt+dW_t.
\end{equation*}
Notice in the same vein that $B_{t-T}^H\simeq B_{T-t}^H$ on $[0,T]$ for $T>0$, so the study on
$[-T,0]$ is related to the time reversal of $B^H$ on $[0,T]$; some general results for this problem
were obtained in \cite{darses:sauss07}.
\end{remark}

\subsection{Some non canonical representations}\label{noncanon}

Let us come back to general $H$-self-similar centred Gaussian processes $\Xi_t$, $t\ge0$. In
Theorem \ref{thlth}, we have proved the equality in law
\begin{equation*}
   \Xi_t\simeq T_{H,L}\Xi(t)=\Xi_t-2L\,t^{H-L}\int_0^ts^{L-H-1}\Xi_sds
\end{equation*}
for $L>0$. When $\Xi=W$ is a standard Brownian motion so that $H=1/2$, this is the classical L\'evy
family of non canonical representations of $W$ with respect to itself. We now verify that this
property of non canonical representation holds in many cases, in the sense that
$\tribuf_t(T_{H,L}\Xi)$ is strictly included in $\tribuf_t(\Xi)$ for $t>0$ (it is of course
sufficient to consider the case $t=1$). In the following theorem we need some notions about
Cameron-Martin spaces and Wiener integrals (see a short introduction in Appendix \ref{camarsec}).

\begin{theorem}\label{thlnc}
Let $\Xi=(\Xi_t;\;0\le t\le1)$ be the restriction to $[0,1]$ of a $H$-self-similar centred Gaussian
process for $H>0$. Let $\wiener$ be a separable Fr\'echet space of paths in which $\Xi$ lives, and
let $\hilbert$ be its Cameron-Martin space. Suppose that the function $\psi(t)=t^{H+L}$ is in
$\hilbert$, and denote by $\langle\Xi,\psi\rangle_\hilbert$ its Wiener integral. Then
\begin{equation*}
   \sigma(\Xi)=\sigma(T_{H,L}\Xi)\vee\sigma(\langle\Xi,\psi\rangle_\hilbert)
\end{equation*}
where the two $\sigma$-algebras of the right-hand side are independent.
\end{theorem}

\begin{proof}
The operator $T_{H,L}$ operates on $\hilbert$, and it is easy to check that functions proportional
to $\psi$ constitute the kernel of $T_{H,L}$. On the other hand, for any $h$ in $\hilbert$,
$h\ne0$, we can write the decomposition
\begin{equation*}
   \Xi=\langle\Xi,h\rangle_\hilbert\frac{h}{|h|_\hilbert^2}+\Bigl(\Xi-\langle\Xi,h\rangle_\hilbert\frac{h}{|h|_\hilbert^2}\Bigr)
\end{equation*}
where the two terms are independent: this is because independence and orthogonality are equivalent
in Gaussian spaces, and
\begin{equation*}
   \esp\Bigl[\langle\Xi,h\rangle_\hilbert\bigl\langle\Xi-\langle\Xi,h\rangle_\hilbert\frac{h}{|h|_\hilbert^2}
   ,h'\bigr\rangle_\hilbert\Bigr]=0
\end{equation*}
for any $h'$ in $\hilbert$ (apply \eqref{whwh}). Thus
\begin{equation*}
   T_{H,L}\Xi=\langle\Xi,h\rangle_\hilbert\frac{T_{H,L}h}{|h|_\hilbert^2}
   +\mbox{process independent of }\langle\Xi,h\rangle_\hilbert,
\end{equation*}
and $T_{H,L}\Xi$ is independent of $\langle\Xi,h\rangle_\hilbert$ if and only if $h$ is in the
kernel of $T_{H,L}$, so if and only if $h$ is proportional to $\psi$. Thus the Gaussian space of
$\Xi$, which is generated by $\langle\Xi,h\rangle_\hilbert$, $h\in\hilbert$, is the orthogonal sum
of the Gaussian space generated by $T_{H,L}\Xi$ and of the variables proportional to
$\langle\Xi,\psi\rangle_\hilbert$. We deduce the theorem. \qed
\end{proof}

Notice that on the other hand, the transformation $T_{H,L}$ becomes injective on the whole time
interval $\reel_+$, so $\sigma(\Xi)$ and $\sigma(T_{H,L}\Xi)$ coincide; actually, the theorem
cannot be used on $\reel_+$ because $\psi$ is no more in $\hilbert$; this can be viewed from the
fact that $\Xi$ lives in the space of functions $f$ such that $t^{-H-1-\eps}f(t)$ is integrable on
$[1,\infty)$ (for $\eps>0$), so $\hilbert$ is included in this space, whereas $\psi$ does not
belong to it for $\eps\le L$.

In the case where $\Xi$ is the standard Brownian motion $W$, we obtain the well known property
\begin{equation}\label{triwt}
   \tribuf_t(W)=\tribuf_t(T_{1/2,L}W)\vee\sigma\bigl(\Pitilde^{L-1/2}W(t)\bigr).
\end{equation}
Let us prove that this property enables to write Theorem \ref{thlnc} in another form when $\Xi$ has
a canonical representation with respect to $W$, see also \cite{cjost07}.

\begin{theorem}\label{ncvol}
Consider the standard Brownian motion $W$ on $\reel_+$, and let
\begin{equation*}
   \Xi_t=(AW)(t)=\int_0^tK(t,s)dW_s
\end{equation*}
be given by a kernel $K$ satisfying $K(\lambda t,\lambda s)=\lambda^{H-1/2}K(t,s)$ for any
$\lambda>0$ and some $H>0$. Suppose that $\tribuf_t(\Xi)=\tribuf_t(W)$ (the representation is
canonical). Then $\Xi$ is a $H$-self-similar process, and we have
\begin{equation}\label{tribind}
   T_{H,L}\Xi=AT_{1/2,L}W,\quad
   \tribuf_t(\Xi)=\tribuf_t(T_{H,L}\Xi)\vee\sigma\bigl(\Pitilde^{L-1/2}W(t)\bigr)
\end{equation}
where the two $\sigma$-algebras of the right side are independent.
\end{theorem}

\begin{proof}
The scaling condition on $K$ implies that $\Xi$ is $H$-self-similar. It can be viewed for instance
as a random variable in the space of functions $f$ such that $t^{\eps-1-H,-\eps-H-1}f(t)$ is
integrable on $\reel_+^\star$. On the other hand, notice that
\begin{equation}\label{thlpi}
   T_{H,L}=\Pi^{H-1/2}\Pi^{1/2-L}\Pitilde^{2L}\Pi^{-L-1/2}\Pi^{1/2-H}=\Pi^{H-1/2}T_{1/2,L}\Pi^{1/2-H}
\end{equation}
from \eqref{thldef}, and consider the linear functional $\Pi^{1/2-H}A$ mapping $W$ to the
$1/2$-self-similar process $\Pi^{1/2-H}\Xi$. The monomials $\psi_\beta(t)=t^\beta$, $\beta>1/2$,
generate the Cameron-Martin space $\hilbert_{1/2}$ of $W$; we deduce from the scaling condition
that they are eigenfunctions of $\Pi^{1/2-H}A$ and of $T_{1/2,L}$, so the commutativity relation
\begin{equation}\label{piat}
   \Pi^{1/2-H}AT_{1/2,L}=T_{1/2,L}\Pi^{1/2-H}A
\end{equation}
holds on fractional polynomials, and therefore on $\hilbert_{1/2}$ and on the paths of $W$ (a
linear functional of $W$ which is zero on the Cameron-Martin space must be zero on $W$). We deduce
from \eqref{thlpi} and \eqref{piat} that
\begin{equation*}
   T_{H,L}\Xi=\Pi^{H-1/2}T_{1/2,L}\Pi^{1/2-H}AW=\Pi^{H-1/2}\Pi^{1/2-H}AT_{1/2,L}W=AT_{1/2,L}W
\end{equation*}
and the first part of \eqref{tribind} is proved. We have moreover assumed that
$\tribuf_t(AW)=\tribuf_t(W)$; this can be applied to the Brownian motion $T_{1/2,L}W$ so
$\tribuf_t(AT_{1/2,L}W)=\tribuf_t(T_{1/2,L}W)$. Thus, by applying \eqref{triwt},
\begin{align*}
   \tribuf_t(\Xi)&=\tribuf_t(W)=\tribuf_t(T_{1/2,L}W)\vee\sigma(\Pitilde^{L-1/2}W(t))\\
   &=\tribuf_t(AT_{1/2,L}\Xi)\vee\sigma(\Pitilde^{L-1/2}W(t))
   =\tribuf_t(T_{H,L}\Xi)\vee\sigma(\Pitilde^{L-1/2}W(t))
\end{align*}
so the second part of \eqref{tribind} is also proved. \qed
\end{proof}

\begin{remark}
Another proof of the second part of \eqref{tribind} is to use directly Theorem \ref{thlnc}; we
verify that on $[0,1]$
\begin{equation*}
   \Pitilde^{L-1/2}W(1)=\langle W,\phi\rangle_{\hilbert_{1/2}}
   =\langle\Xi,A\phi\rangle_\hilbert
\end{equation*}
for $\phi(t)=t^{L+1/2}/(L+1/2)$, and $A\phi$ is proportional to the function $\psi(t)=t^{L+H}$ from
the scaling condition.
\end{remark}

\begin{theorem}\label{noncan}
Consider on $\reel_+$ the family of fractional Brownian motions $B^H=G_{0+}^{1/2,H}W$, so that
$B^H=G_{0+}^{J,H}B^J$. Then, for any $L>0$, the process $B^{H,L}=T_{H,L}B^H$ is a $H$-fractional
Brownian motion satisfying the relation $B^{H,L}=G_{0+}^{J,H}B^{J,L}$. Moreover, for any $t$,
\begin{equation}
    \tribuf_t(B^H)=\tribuf_t(B^{H,L})\vee\sigma\bigl(\Pitilde^{L-1/2}W(t)\bigr),
\end{equation}
and the two $\sigma$-algebras of the right-hand side are independent.
\end{theorem}

\begin{proof}
This is a direct application of Theorem \ref{ncvol} with $A=G_{0+}^{1/2,H}$. The first part of
\eqref{tribind} implies that
\begin{equation*}
   B^{H,L}=G_{0+}^{1/2,H}T_{1/2,L}W,
\end{equation*}
and the relationship between $B^{J,L}$ and $B^{H,L}$ follows from the composition rule satisfied by
the family $G_{0+}^{J,H}$. \qed
\end{proof}

\section{Riemann-Liouville processes}
\label{riemliou}

In this section, we compare the fractional Brownian motion $B^H$ with the process
$X^H=I_{0+}^{H-1/2}W$.

\subsection{Comparison of processes}

The processes
\begin{equation}\label{xjhdef}
    X_t^H=I_{0+}^{H-1/2}W(t)=\frac1{\Gamma(H-1/2)}\int_0^t(t-s)^{H-1/2}dW_s
\end{equation}
defined on $\reel_+$ are often called Riemann-Liouville processes. Notice that these processes can
be defined for any $H>0$. When $0<H<1$, these processes have paths in $\holder^{H-}$ on bounded
time intervals from Theorem~\ref{riemth}, and can be viewed as good approximations of fractional
Brownian motions $B^H$ for large times, as it is explained in the following result.

\begin{theorem}\label{coupl1}
For $0<H<1$, we can realise jointly the two processes $(X^H,B^H)$ on $\reel_+$, so that $X^H-B^H$
is $C^\infty$ on $\reel_+^\star$. Moreover, for $T>0$, $S>0$ and $1\le p<\infty$,
\begin{equation}\label{bhxhsm}
   \Bigl\|\sup_{0\le t\le T}\bigl|(X_{S+t}^H-X_S^H)-(B_{S+t}^H-B_S^H)\bigr|\Bigr\|_p
   \le C_p\,S^{H-1}T
\end{equation}
(where $\|.\|_p$ denotes the $L^p(\Omega)$-norm for the probability space).
\end{theorem}

\begin{proof}
Let $(B_t^H;\;t\ge0)$ be defined by $B^H=\Itilde_+^{H-1/2}W$ for a standard Brownian motion
$(W_t;\;t\in\reel)$. The process $W$ can be decomposed into the two independent processes
$W_t^+=W_t$ and $W_t^-=W_{-t}$ for $t\ge0$, and consequently, the process $B^H$ is decomposed into
$B^H=X^H+Y^H$ where
\begin{equation*}
   X^H=\Itilde_+^{H-1/2}\bigl(W\,1_{\reel_+}\bigr)=I_{0+}^{H-1/2}W^+
\end{equation*}
is a Riemann-Liouville process, and $Y^H=\Itilde_+^{H-1/2}\bigl(W\,1_{\reel_-}\bigr)$ can be
written by means of Remark \ref{fzeroplus}; more precisely, $Y^H=I_\triangle^{H-1/2}W^-$, where
\begin{equation}\label{itriang}
    I_\triangle^\alpha f(t)=\frac1{\Gamma(\alpha)}\int_0^\infty\bigl((t+s)^{\alpha-1}
    -s^{\alpha-1}\bigr)f(s)ds.
\end{equation}
We deduce from this representation that $Y^H$ is $C^\infty$ on $\reel_+^\star$, so the first
statement is proved. On the other hand, it follows from the scaling property that its derivative is
$(H-1)$-self-similar, and is therefore of order $t^{H-1}$ in $L^p(\Omega)$; thus the left hand side
of \eqref{bhxhsm} is bounded by
\begin{equation*}
  \Bigl\|\int_S^{S+T}\bigl|D^1Y_u^H\bigr|du\Bigr\|_p\le C_p\int_S^{S+T}u^{H-1}du\le C_p\,S^{H-1}T.
\end{equation*}
\qed
\end{proof}

\begin{remark}
Inequality \eqref{bhxhsm} says that the process $X_t^{S,H}=X_{S+t}^H-X_S^H$ is close to a
fractional Brownian motion when $S$ is large; it actually provides an upper bound for the
Wasserstein distance between the laws of these two processes. A result about the total variation
distance will be given later (Theorem \ref{rlabs}).
\end{remark}

Instead of using the representation of $B^H=\Itilde_+^{H-1/2}W$ on $\reel$, we can consider the
coupling based on the canonical representation of $B^H$ on $\reel_+$. It appears that in this case
$X^H-B^H$ is not $C^\infty$ but is still differentiable. In particular, we can deduce that the
estimation \eqref{bhxhsm} also holds for the coupling of Theorem \ref{coupl2}.

\begin{theorem}\label{coupl2}
Consider on $\reel_+$ the family $B^H=G_{0+}^{1/2,H}W$ and the family $X^H$ defined by
\eqref{xjhdef}. Then $X^H-B^H$ is differentiable on $\reel_+^\star$.
\end{theorem}

\begin{proof}
For $f$ smooth with compact support in $\reel_+^\star$, Theorem \ref{smooth} and the
expression \eqref{gipiipid} for $G_{0+}^{J,H}$ shows that $G_{0+}^{J,H}f$ and
$I_{0+}^{H-J}f$ are smooth, and
\begin{equation*}
    D^1\bigl(G_{0+}^{J,H}-I_{0+}^{H-J}\bigr)
    =\bigl(\Pi^{H+J-1}I_{0+}^{H-J}\Pi^{1-H-J}-I_{0+}^{H-J}\bigr)D^1.
\end{equation*}
We therefore deduce from \eqref{pimgam} that
\begin{align*}
    &\frac{d}{dt}\bigl(G_{0+}^{J,H}-I_{0+}^{H-J}\bigr)f(t)\\
    &=\frac1{\Gamma(H-J)}\int_0^t\biggl(\bigl(\frac ts\bigr)^{H+J-1}-1\biggr)(t-s)^{H-J-1}D^1f(s)ds\\
    &=\frac{f(t)}{t}U(t)
    +\frac1{\Gamma(H-J)}\int_0^t\biggl(\bigl(\frac ts\bigr)^{H+J-1}-1\biggr)(t-s)^{H-J-1}\bigl(D^1f(s)-f(t)/t\bigr)ds\\
    &=\frac{f(t)}{t}U(t)
    -\frac1{\Gamma(H-J)}\int_0^t\partial_s\biggl[\biggl(\bigl(\frac ts\bigr)^{H+J-1}-1\biggr)
    (t-s)^{H-J-1}\biggr]\bigl(f(s)-\frac stf(t)\bigr)ds
\end{align*}
with
\begin{equation*}
   U(t)=\frac1{\Gamma(H-J)}\int_0^t\biggl(\bigl(\frac ts\bigr)^{H+J-1}-1\biggr)(t-s)^{H-J-1}ds
\end{equation*}
proportional to $t^{H-J}$. This equality can be extended to any function $f$ of $\holder^{J-}$, so
in particular to $W$ in the case $J=1/2$; we deduce the differentiability announced in the theorem.
\qed
\end{proof}

\subsection{The Riemann-Liouville Cameron-Martin space}

Cameron-Martin spaces are Hilbert spaces which characterise the law of centred Gaussian variables,
so in particular of centred Gaussian processes, see Appendix \ref{camarsec}. The Cameron-Martin
spaces $\hilbert_H$ of $H$-fractional Brownian motions are deduced from each other by means of the
transforms of Theorems \ref{reprth} or \ref{molgolo}, so that
\begin{equation*}
    \hilbert_H=\Itilde_+^{H-J}(\hilbert_J)=\Itilde_-^{H-J}(\hilbert_J),
    \qquad\hilbert_H=G_{0+}^{J,H}(\hilbert_J)=\Itilde_-^{H-J}(\hilbert_J)
\end{equation*}
respectively on $\reel$ and $\reel_+$; the space $\hilbert_{1/2}$ is the classical space of
absolutely continuous functions $h$ such that $h(0)=0$ and $D^1h$ is in $L^2$. Similarly, the
Cameron-Martin space of the Riemann-Liouville process $X^H$ on $\reel_+$ is
\begin{equation*}
   \hilbert_H'=I_{0+}^{H-1/2}\hilbert_{1/2}=I_{0+}^{H+1/2}L^2.
\end{equation*}
In particular, if $f$ is a smooth function on $\reel_+$ such that $f(0)=0$, then, on the time
interval $[0,T]$,
\begin{align}
  |f|_{\hilbert_H'}&=\bigl|D^1I_{0+}^{1/2-H}f\bigr|_{L^2}\nonumber\\
  &\le C\biggl(\sup|D^1f|\Bigl(\int_0^T\bigl(t^{1/2-H}\bigr)^2dt\Bigr)^{1/2}+\sup|D^2f|
  \Bigl(\int_0^T\bigl(t^{3/2-H}\bigr)^2dt\Bigr)^{1/2}\biggr)\nonumber\\
  &\le C'\Bigl(T^{1-H}\sup|D^1f|+T^{2-H}\sup|D^2f|\Bigr)\label{hpsm}
\end{align}
from Theorem \ref{smooth}.

We now explain the proof of a result mentioned in \cite{decreuse:ustu99} (Theorem 2.1) and taken
from \cite{samko:kil:mari93}. We use the equivalence of Hilbert spaces ($\hilbert\sim\hilbert'$)
defined in \eqref{equivdef}. A probabilistic interpretation of this equivalence is given in
Appendix \ref{camarsec}, see \eqref{cmembed}.

\begin{theorem}\label{rlcamar}
For $0<H<1$, the spaces $\hilbert_H$ and $\hilbert_H'$ are equivalent on $\reel_+$.
\end{theorem}

\begin{proof}
The proof is divided into the two inclusions; for the second one, we are going to use an analytical
result proved in Appendix~\ref{analytical}. We can of course omit the case $H=1/2$.

\medskip\noindent\emph{Proof of $\hilbert_H'\subset\hilbert_H$.} We have seen in the proof of Theorem
\ref{coupl1} that $B^H$ can be written as the sum of the Riemann-Liouville process $X^H$ and of an
independent process $Y^H$. If we denote by $\hilbert_H^\triangle$ the Cameron-Martin space of
$Y^H$, then this decomposition implies (see \eqref{contract}) that
\begin{equation}\label{hhh}
    \hilbert_H=\hilbert_H'+\hilbert_H^\triangle\quad\text{with}\quad
    |h|_{\hilbert_H}=\inf\Bigl\{\bigl(|h_1|_{\hilbert_H'}^2+|h_2|_{\hilbert_H^\triangle}^2\bigr)^{1/2};
    \>h=h_1+h_2\Bigr\}.
\end{equation}
In particular, $\hilbert_H'\subset\hilbert_H$ with $|h|_{\hilbert_H}\le|h|_{\hilbert_H'}$.

\medskip\noindent\emph{Proof of $\hilbert_H\subset\hilbert_H'$.} It is sufficient from \eqref{hhh} to prove
that $\hilbert_H^\triangle$ is continuously embedded into $\hilbert_H'$. Let $h$ be in
$\holder^{1/2}$; then $|h(t)|\le|h|_{\holder^{1/2}}\sqrt t$, and we can deduce from \eqref{itriang}
that $I_\triangle^{H-1/2}h$ is $C^\infty$ on $\reel_+^\star$, and that the derivative of order $k$
is dominated by $|h|_{\holder^{1/2}}t^{H-k}$. Theorem \ref{smooth} enables to deduce that
$Ah=I_{0+}^{1/2-H}I_\triangle^{H-1/2}h$ is also smooth, and we have from \eqref{d1i} that
$D^1Ah(t)$ is dominated by $|h|_{\holder^{1/2}}/\sqrt t$. Moreover, the scaling condition
\eqref{scalc} is satisfied, so we deduce from Theorem \ref{ahlambda} that $A$ is a continuous
endomorphism of $\hilbert_{1/2}$. By composing with $I_{0+}^{H-1/2}$, we obtain that
$\bigl|I_\triangle^{H-1/2}g\bigr|_{\hilbert_H'}$ is dominated by $|g|_{\hilbert_{1/2}}$, so
\begin{equation*}
   |h|_{\hilbert_H^\triangle}=\inf\Bigl\{|g|_{\hilbert_{1/2}};\;h=I_\triangle^{H-1/2}g\Bigr\}
   \ge c|h|_{\hilbert_H'}.
\end{equation*}
\qed
\end{proof}

\begin{remark}\label{endom}
Let us give another interpretation of Theorem \ref{rlcamar}. By comparing $\reel$ and $\reel_+$,
the fractional Brownian motion on $\reel_+$ can be obtained as a restriction of the fractional
Brownian motion on $\reel$. This property can be extended to the Cameron-Martin spaces, and
applying \eqref{contract}, we deduce that $\hilbert_H(\reel_+)$ consists of the restrictions to
$\reel_+$ of functions of $\hilbert_H(\reel)$, and
\begin{equation*}
   |h|_{\hilbert_H(\reel_+)}=\inf\Bigl\{|g|_{\hilbert_H(\reel)};\;g=h\text{ on }\reel_+\Bigr\},
\end{equation*}
so $|h|_{\hilbert_H(\reel_+)}\le|h\,1_{\reel_+}|_{\hilbert_H(\reel)}$ for $h$ defined on $\reel_+$.
On the other hand,
\begin{align*}
   |h|_{\hilbert_H'}&=\bigl|I_{0+}^{1/2-H}h\bigr|_{\hilbert_{1/2}(\reel_+)}
   =\bigl|(I_{0+}^{1/2-H}h)1_{\reel_+}\bigr|_{\hilbert_{1/2}(\reel)}\\
   &=\bigl|\Itilde_+^{H-1/2}((I_{0+}^{1/2-H}h)1_{\reel_+})\bigr|_{\hilbert_H(\reel)}
   =\bigl|h\,1_{\reel_+}\bigr|_{\hilbert_H(\reel)}.
\end{align*}
Thus $|h|_{\hilbert_H(\reel_+)}\le|h|_{\hilbert_H'}$, and $\hilbert_H'$ is continuously embedded in
$\hilbert_H(\reel_+)$. The inverse inclusion means that
\begin{equation*}
   \bigl|h\,1_{\reel_+}\bigr|_{\hilbert_H(\reel)}
   \le C\>\inf\Bigl\{|g|_{\hilbert_H(\reel)};\;g=h\text{ on }\reel_+\Bigr\},
\end{equation*}
for $h$ defined on $\reel_+$, and this is equivalent to
\begin{equation*}
   \bigl|g\,1_{\reel_+}\bigr|_{\hilbert_H(\reel)}\le C\,|g|_{\hilbert_H(\reel)}
\end{equation*}
for $g$ defined on $\reel$; thus this means that $g\mapsto g\,1_{\reel_+}$ is a continuous
endomorphism of $\hilbert_H(\reel)$. This is a known analytical result, see also Lemma 1 in
\cite{norros:sak09}.
\end{remark}

\begin{remark}\label{cmodd}
Consider on $\reel_+$ the even and odd parts $B_t^{H\pm}=(B_t^H\pm B_{-t}^H)/2$ of $B^H$. These two
processes are independent (this is easily verified by computing the covariance), and
$B^H1_{\reel_+}=B^{H+}+B^{H-}$, so their Cameron-Martin spaces $\hilbert_{H\pm}$ are continuously
embedded into $\hilbert_H(\reel_+)$. On the other hand
\begin{align*}
   |h|_{\hilbert_{H\pm}}&=\inf\Bigl\{|g|_{\hilbert_H(\reel)};\;h(t)=\frac12(g(t)\pm g(-t))
   \text{~on~}\reel_+\Bigr\}\\
   &\le2\bigl|h\,1_{\reel_+}\bigr|_{\hilbert_H(\reel)}=2|h|_{\hilbert_H'}
   \le C|h|_{\hilbert_H(\reel_+)}
\end{align*}
by means of the result of Remark \ref{endom}, so the three spaces $\hilbert_{H\pm}$ and
$\hilbert_H(\reel_+)$ are equivalent.
\end{remark}

\begin{remark}\label{restrict}
Notice that the endomorphism of Remark \ref{endom} maps the function $h(t)$ to the function
$h(t_+)$; by applying the invariance by time reversal, we deduce that the operator mapping $h(t)$
to $h(1-(1-t)_+)$ is also continuous, so by composing these two operators, we see that the operator
mapping $h(t)$ to the function
\begin{equation}\label{hstar}
   h^\star(t)=
   \begin{cases}
     0&\text{if $t\le0$,}\\
     h(t)&\text{if $0\le t\le1$,}\\
     h(1)&\text{if $t\ge1$,}
   \end{cases}
\end{equation}
is a continuous endomorphism of $\hilbert_H$. On the other hand, we have
\begin{equation*}
   |h|_{\hilbert_H([0,1])}=\inf\Bigl\{|g|_{\hilbert_H(\reel)};\;g=h\text{ on }[0,1]\Bigr\}.
\end{equation*}
Thus $h\mapsto h^\star$ is continuous from $\hilbert_H([0,1])$ into $\hilbert_H(\reel)$.
\end{remark}

\subsection{Equivalence and mutual singularity of laws}

In Theorem \ref{rlcamar}, we have proved that the Cameron-Martin spaces of $B^H$ and $X^H$ are
equivalent. It is known that the laws of two centred Gaussian processes are either equivalent, or
mutually singular, see Appendix \ref{equiv}; the equivalence of Cameron-Martin spaces is necessary
for the equivalence of the laws, but is of course not sufficient (compare for instance a standard
Brownian motion $W_t$ with $2W_t$). In subsequent results, the equivalence or mutual singularity of
laws of processes should be understood by considering these processes as variables with values in
the space of continuous functions.

\begin{theorem}\label{rlabs}
Let $0<H<1$. For any $S>0$, the laws of $B_t^H$ and $X_t^{S,H}=X_{S+t}^H-X_S^H$ are equivalent on
any time interval $[0,T]$; more precisely, the relative entropies of $B^H$ and $X^{S,H}$ with
respect to each other are dominated by $S^{2H-2}$ as $S\uparrow\infty$, and therefore tend to 0; in
particular, the total variation distance between the laws of $X^{S,H}$ and $B^H$ is dominated by
$S^{H-1}$. In the case $S=0$, the two laws are mutually singular as soon as $H\ne1/2$.
\end{theorem}

\begin{proof}
Let us consider separately the cases $S>0$ and $S=0$.

\medskip\noindent\emph{Equivalence for $S>0$.} Consider the coupling and notations of Theorem \ref{coupl1},
so that the process $B_t^H=X_t^H+Y_t^H$ is written as the sum of two independent processes. This
implies that $B^{S,H}=X^{S,H}+Y^{S,H}$, where $B^{S,H}$ and $Y^{S,H}$ are defined similarly to
$X^{S,H}$. Theorem \ref{rlcamar} states that the Cameron-Martin spaces of $X^H$ and $B^H$ are
equivalent; this implies that the Cameron-Martin space of $X^{S,H}$ is equivalent to the
Cameron-Martin space of $B^{S,H}$ which is $\hilbert_H$, and is therefore also equivalent to
$\hilbert_H'=I_{0+}^{H+1/2}L^2(\reel_+)$; thus it contains smooth functions taking value 0 at 0.
But the perturbation $Y^{S,H}$ is smooth, so the equivalence of the laws of $B^{S,H}$ and $X^{S,H}$
follows from the Cameron-Martin theorem for an independent perturbation. Moreover, \eqref{entrop}
yields an estimation of the relative entropies
\begin{align*}
   \max\bigl(\inform(B^H,X^{S,H}),\inform(X^{S,H},B^H)\bigr)&\le\frac12\esp|Y^{S,H}|_{\hilbert_H}^2
   \le C\,\esp|Y^{S,H}|_{\hilbert_H'}^2\\
   &\le C_T\,\esp\Bigl(\sup_{[0,T]}|D^1Y^{S,H}|+\sup_{[0,T]}|D^2Y^{S,H}|\Bigr)^2
\end{align*}
from \eqref{hpsm}. The derivative $D^kY_t^H$ is $O(t^{H-k})$ in $L^2(\Omega)$ from the scaling
property, so
\begin{equation*}
   \sup|D^1Y_t^{S,H}|=\sup|D^1Y_{S+t}^H|\le|D^1Y_S^H|+\int_0^T|D^2Y_{S+t}^H|dt=O(S^{H-1})
\end{equation*}
as $S\uparrow\infty$. The second derivative is even smaller (of order $S^{H-2}$). Thus the relative
entropies are dominated by $S^{2H-2}$. In particular, the total variation distance is estimated
from Pinsker's inequality \eqref{pinsker}.

\medskip\noindent\emph{Mutual singularity for $S=0$.} This is a consequence of Theorem \ref{dichot};
the two processes are self-similar, the initial $\sigma$-algebra $\tribuf_{0+}(B^H)$ is almost
surely trivial (Remark~\ref{blumen}), so it is sufficient to prove that they do not have the same
law. But this is evident since $B^H$ can be written as the sum of $X^H$ and of an independent
process $Y^H$ which is not identically zero. \qed
\end{proof}

\begin{remark}
In the case $S=0$, Theorem \ref{covgauss} provides a criterion to decide whether a process $\Xi$
has the law of $B^H$ or $X^H$. The variances of these two processes differ (they can be computed
from the calculation of Appendix \ref{variance}), so we can decide between them by looking at the
small time behaviour of $\int_t^1s^{-2H-1}(\Xi_s)^2ds$. Actually, by applying the invariance by
time inversion, we can also look at the behaviour in large time.
\end{remark}

For the following result, we recall that the mutual information of two variables $X_1$ and $X_2$ is
defined as the entropy of $(X_1,X_2)$ relative to two independent copies of $X_1$ and $X_2$. We
want to estimate the dependence between the increments of $B^H$ on some interval $[S,S+T]$,
$S\ge0$, and its increments before time 0, and in particular prove that the two processes are
asymptotically independent when $S\uparrow+\infty$. This result and other estimates were proved in
\cite{norros:sak09} with a more analytical method; an asymptotic independence result is also given
in \cite{picard08}.

\begin{theorem}\label{bhindep}
Let $H\ne1/2$. The joint law of the two processes $(B_t^{S,H}=B_{S+t}^H-B_S^H;\;0\le t\le T)$ and
$(B_t^H;\;t\le0)$ is equivalent to the product of laws as soon as $S>0$, and the Shannon mutual
information is $O(S^{2H-2})$ as $S\uparrow\infty$. If $S=0$, the joint law and the product of laws
are mutually singular.
\end{theorem}

\begin{proof}
We consider separately the two cases.

\medskip\noindent\emph{Equivalence for $S>0$.} Let $(W_t;\;t\in\reel)$ and $(\Wbar_t;\;t\in\reel)$ be two
standard Brownian motions such that $\Wbar_t=W_t$ for $t\ge0$ and $(\Wbar_t;\;t\le0)$ is
independent of $W$. We then consider the two fractional Brownian motions $B^H=\Itilde_+^{H-1/2}W$
and $\Lambda^H=\Itilde_+^{H-1/2}\Wbar$. With the notation of Theorem \ref{coupl1}, they can be
written on $\reel_+$ as $B^H=X^H+Y^H$ and $\Lambda^H=X^H+\Ybar^H$, so $\Lambda^H=B^H+\Ybar^H-Y^H$;
by looking at the increments after time $S$, we have $\Lambda^{S,H}=B^{S,H}+\Ybar^{S,H}-Y^{S,H}$.
Conditionally on $\tribuf_0(W,\Wbar)=\tribuf_0(B^H,\Lambda^H)$, the process $\Ybar^{S,H}-Y^{S,H}$
becomes a deterministic process which is almost surely in $\hilbert_H$ (see the proof of
Theorem~\ref{rlabs}), so the conditional laws of
\begin{equation*}
   (B_t^{S,H},\;0\le t\le T;\;B_t^H,\;t\le0)\quad\mbox{and}\quad(\Lambda_t^{S,H},\;0\le t\le T;\;B_t^H,\;t\le0)
\end{equation*}
are equivalent. We deduce that the unconditional laws are also equivalent. Moreover, the two
processes of the right side are independent, and $\Lambda^{S,H}\simeq B^{S,H}$, so the
equivalence of laws stated in the theorem is proved. On the other hand, the relative entropies of
\begin{equation*}
   (B_t^{S,H},\;0\le t\le T;\;B_t^H,\;t\le0;\;\Lambda_t^H,\;t\le0)
\end{equation*}
and
\begin{equation*}
    (\Lambda_t^{S,H},\;0\le t\le T;\;B_t^H,\;t\le0;\;\Lambda_t^H,\;t\le0)
\end{equation*}
with respect to each other are equal to
\begin{equation*}
   \frac12\esp\bigl|\Ybar^{S,H}-Y^{S,H}\bigr|_{\hilbert_H}^2\le2\esp\bigl|Y^{S,H}\bigr|_{\hilbert_H}^2
   =O(S^{2H-2})
\end{equation*}
(proceed as in Theorem \ref{rlabs}). If we project on the two first components, we deduce that the
mutual information that we are looking for is smaller than this quantity.

\medskip\noindent\emph{Mutual singularity for $S=0$.} If we compare the law of
$(B_t^H,B_{-t}^H;\;0\le t\le T)$ with the law of two independent copies of the fractional Brownian
motion, we have two self-similar Gaussian processes with different laws, so the laws are mutually
singular from Theorem \ref{dichot}. \qed
\end{proof}

\begin{remark}\label{odde}
As an application, we can compare $B^H$ with its odd and even parts. Let $B$ and $B'$ be two
independent copies of $B^H$. Let $S>0$. From Theorem \ref{bhindep}, we have on $[0,T]$ the
equivalence of laws
\begin{align*}
   \bigl(B_{S+t}^H-B_S^H\bigr)\pm\bigl(B_{-S-t}^H-B_{-S}^H\bigr)
   \sim\bigl(B_{S+t}-B_S\bigr)\pm\bigl(B_{-S-t}'-B_{-S}'\bigr)&\simeq\sqrt2(B_{S+t}^H-B_S^H\bigr)\\
   &\simeq\sqrt2 B_t^H.
\end{align*}
Thus the law of the increments of $(B_t^H\pm B_{-t}^H)/\sqrt2$ on $[S,S+T]$ have a law equivalent
to the law of $B^H$. For $S=0$, the Cameron-Martin spaces are equivalent (Remark \ref{cmodd}), but
the laws can be proved to be mutually singular from Theorem \ref{dichot}.
\end{remark}

\section{Series expansions}
\label{series}

Let us try to write $B^H$ on $[0,1]$ as some series of type
\begin{equation*}
    B_t^H=\sum_nh_n(t)\xi_n
\end{equation*}
where $h_n$ are deterministic functions and $\xi_n$ are independent standard Gaussian variables.
Such expansions have been described in the standard case $H=1/2$ by \cite{ito:nisio68}, and
actually, an expansion valid for the standard Brownian motion $W$ can be transported to $B^H$ by
means of the operator $G_{0+}^{1/2,H}$, see \cite{gilsing:sott03}.

If we look more precisely for a trigonometric expansion, we can apply \cite{dzhapa:zanten04} where
the functions $h_n$ are trigonometric functions, the coefficients of which are related to some
Bessel function depending on $H$. However, we are here more interested in trigonometric functions
which do not depend on $H$.

\subsection{A trigonometric series}

Suppose that we are interested in the Fourier series of $(B_t^H;\;0\le t\le1)$. The problem is that
the Fourier coefficients are not independent, since this property is already known to be false for
$H=1/2$. What is known for $H=1/2$ is that $W_t$ can be represented by means of \eqref{wtrigo},
\eqref{wanti} or \eqref{karh} for independent standard Gaussian variables $(\xi_n,\xi_n';n\ge1)$;
the series converges in $L^2(\Omega)$, uniformly in $t$, and one easily deduces the Fourier series
of $W$ from \eqref{wtrigo}. Similar representations cannot hold on $[0,1]$ for the fractional
Brownian motion as soon as $H\ne1/2$, but it appears that one can find a representation mixing
\eqref{wtrigo} and \eqref{wanti},
\begin{equation}\label{bhsimsum}
     B_t^H\simeq a_0^H\xi_0t
     +\sum_{n\ge1}a_n^H\Bigl(\bigl(\cos(\pi nt)-1\bigr)\xi_n+\sin(\pi nt)\xi_n'\Bigr)
\end{equation}
on $[0,1]$. This question has been studied in \cite{istas05} and \cite{igloi05} respectively for
the cases $H<1/2$ and $H>1/2$. The sign of $a_n^H$ is of course irrelevant so we will choose
$a_n^H\ge0$. We follow a general technique for finding series expansions of Gaussian processes from
series expansions of their covariance kernels. We are going to find all the possible $a_n^H$ for
which \eqref{bhsimsum} holds; it appears that $a_n^H$, $n\ge1$, is unique as soon as $a_0^H$ has
been chosen in some set of possible values.

\begin{theorem}\label{bhtrigo}
It is possible to find a sequence $(a_n^H;\;n\ge0)$, $a_n^H\ge0$, such that $\sum(a_n^H)^2<\infty$
and \eqref{bhsimsum} holds on $[0,1]$ for independent standard Gaussian variables
$(\xi_0,\xi_n,\xi_n';n\ge1)$. The convergence of the series holds uniformly in $t$, almost surely.
If $H\le1/2$, we have to choose $a_0^H$ in an interval $[0,a(H)]$, $a(H)>0$, and $a_n^H$ is then
uniquely determined; if $H>1/2$ there is only one choice for the sequence. Moreover, except in the
case $H=1/2$, we must have $a_n^H\ne0$ for all large enough $n$. If $H\ne1/2$, then
\eqref{bhsimsum} cannot hold on $[0,T]$ for $T>1$.
\end{theorem}

\begin{proof}
We divide the proof into two parts.

\medskip\noindent\emph{Step 1: Study on $[0,1]$.} It is clear that the convergence of
the series in \eqref{bhsimsum} holds for $t$ fixed (almost surely and in $L^2(\Omega)$); the
uniform convergence comes from the It\^{o}-Nisio theorem \cite{ito:nisio68}. We have to verify that
the right hand side $Z$ has the same covariance kernel as $B^H$ for a good choice of $(a_n^H)$. We
have
\begin{align*}
   \esp[Z_sZ_t]&=(a_0^H)^2st
   +\sum_{n\ge1}(a_n^H)^2\Bigl(\bigl(\cos(\pi nt)-1\bigr)\bigl(\cos(\pi ns)-1\bigr)
   +\sin(\pi nt)\sin(\pi ns)\Bigr)\\
   &=(a_0^H)^2st
   +\sum_{n\ge1}(a_n^H)^2\Bigl(\cos(\pi n(t-s))-\cos(\pi nt)-\cos(\pi ns)+1\Bigr)\\
   &=\bigl(f_H(t)+f_H(s)-f_H(t-s)\bigr)/2
\end{align*}
with
\begin{equation}\label{fht}
   f_H(t)=(a_0^H)^2t^2+2\sum_{n\ge1}(a_n^H)^2\Bigl(1-\cos(\pi nt)\Bigr).
\end{equation}
If we compare this expression with \eqref{bhcov}, it appears that if $f_H$ coincides on $[-1,1]$
with $g_H(t)=\rho\,|t|^{2H}$, then $B^H\simeq Z$ on $[0,1]$; conversely, if $B^H\simeq Z$, then
they have the same variance, so $f_H=g_H$ on $[0,1]$ and therefore on $[-1,1]$ (the two functions
are even). Thus finding an expansion \eqref{bhsimsum} on $[0,1]$ is equivalent to finding
coefficients $a_n^H$ so that $f_H=g_H$ on $[-1,1]$. For any choice of $a_0^H$, one has on $[-1,1]$
the Fourier decomposition
\begin{equation*}
   \rho|t|^{2H}-(a_0^H)^2t^2=b_0^H-2\sum_{n\ge1}b_n^H\cos(\pi nt).
\end{equation*}
Thus the possible expansions correspond to the possible choices of $a_0^H$ such that $b_n^H\ge0$
for $n\ge1$ and $\sum b_n^H<\infty$; then
\begin{equation*}
   \rho|t|^{2H}-(a_0^H)^2t^2=2\sum_{n\ge1}b_n^H(1-\cos(\pi nt))
\end{equation*}
and we take $a_n^H=\sqrt{b_n^H}$ for $n\ge1$. We have
\begin{align}
   b_n^H&=-\rho\int_0^1t^{2H}\cos(\pi nt)dt+(a_0^H)^2\int_0^1t^2\cos(\pi nt)dt\notag\\
   &=\frac{2H}{\pi n}\rho\int_0^1t^{2H-1}\sin(\pi nt)dt-\frac{2(a_0^H)^2}{\pi n}
   \int_0^1t\sin(\pi nt)dt\notag\\
   &=-\frac{2H(2H-1)}{\pi^2n^2}\rho\int_0^1t^{2H-2}\bigl(1-\cos(\pi nt)\bigr)dt\notag\\
   &\qquad\qquad\qquad+\frac{2H}{\pi^2n^2}\rho\bigl(1-(-1)^n\bigr)
   +\frac{2(a_0^H)^2}{\pi^2n^2}(-1)^n.\label{bnh1}
\end{align}
Let us first assume $H<1/2$; then the first term is positive, and the sum of the second and third
terms is nonnegative as soon as $a_0^H\le\sqrt{2\rho H}$. Moreover
\begin{equation}\label{int1cos}
   cn^2\int_0^{1/n}t^{2H}dt\le\int_0^1t^{2H-2}\bigl(1-\cos(\pi nt)\bigr)dt
   \le Cn^2\int_0^{1/n}t^{2H}dt+2\int_{1/n}^\infty t^{2H-2}dt
\end{equation}
so this integral is of order $n^{1-2H}$ (actually a more precise estimate will be proved in Theorem
\ref{btilth}), and we have $b_n^H\asymp n^{-1-2H}$. It is then not difficult to deduce that there
exists a maximal $a(H)\ge\sqrt{2\rho H}$ such that if we choose $a_0^H$ in $[0,a(H)]$, then
$b_n^H\ge0$ for any $n$; the value $a(H)$ is attained when one of the coefficients $b_n^H$ becomes
0. It follows from $b_n^H\asymp n^{-1-2H}$ that $\sum b_n^H<\infty$. Let us now assume $H=1/2$; the
property $b_n^H\ge0$ holds for $a_0^{1/2}\in[0,a(1/2)]=[0,1]$, and $b_n^{1/2}=O(n^{-2})$. Finally,
if $H>1/2$,
\begin{align}
   b_n^H&=\frac{2H(2H-1)}{\pi^2n^2}\rho\int_0^1t^{2H-2}\cos(\pi nt)dt
   +\frac{2(a_0^H)^2-2\rho H}{\pi^2n^2}(-1)^n\label{bnh2}\\
   &=-\frac{2H(2H-1)(2H-2)}{\pi^3n^3}\rho\int_0^1t^{2H-3}\sin(\pi nt)dt
   +\frac{2(a_0^H)^2-2\rho H}{\pi^2n^2}(-1)^n\notag\\
   &=\frac{2H(2H-1)(2H-2)(2H-3)}{\pi^4n^4}\rho\int_0^1t^{2H-4}\bigl(1-\cos(\pi nt)\bigr)dt\notag\\
   &\qquad-\frac{2H(2H-1)(2H-2)}{\pi^4n^4}\rho\bigl(1-(-1)^n\bigr)
   +\frac{2(a_0^H)^2-2\rho H}{\pi^2n^2}(-1)^n.\notag
\end{align}
The integral of the last equality is studied like \eqref{int1cos}, and is of order $n^{3-2H}$, so
the first term of this last equality is positive and of order $n^{-1-2H}$. The second term is
nonnegative and smaller. If we choose $a_0^H\ne\sqrt{\rho H}$, then the third term has an
alternating sign and is the dominant term, so $b_n^H$ is not always positive. Thus we must choose
$a_0^H=\sqrt{\rho H}$, and $b_n^H>0$ for any $n$; we again have $b_n^H\asymp n^{-1-2H}$ so that
$\sum b_n^H<\infty$. Moreover, in the two cases $H<1/2$ and $H>1/2$, we have $a_n^H\asymp
n^{-H-1/2}$, so $a_n^H\ne0$ for all large enough $n$.

\medskip\noindent\emph{Step 2: Study on larger intervals.} Suppose now that \eqref{bhsimsum} holds on
$[0,T]$ for some $T>1$. Then, as in previous step, we should have $f_H(t)=g_H(t)=\rho|t|^{2H}$ on
$[-T,T]$. But $f_H(t)-(a_0^H)^2t^2$ is even and 2-periodic, so
\begin{equation*}
   f_H(1-t)-(a_0^H)^2(1-t)^2=f_H(1+t)-(a_0^H)^2(1+t)^2.
\end{equation*}
Thus
\begin{equation*}
   \rho(1-t)^{2H}-(a_0^H)^2(1-t)^2=\rho(1+t)^{2H}-(a_0^H)^2(1+t)^2
\end{equation*}
for $|t|\le\min(T-1,1)$. By differentiating twice, it appears that this relation is false if
$H\ne1/2$. \qed
\end{proof}

\begin{remark}
For $H=1/2$, we can choose $a_0^{1/2}$ in $[0,1]$, and the expansion \eqref{bhsimsum} is an
interpolation between the decompositions containing respectively only odd terms ($a_0^{1/2}=0$) and
only even terms ($a_0^{1/2}=1$), which are respectively \eqref{wanti} and \eqref{wtrigo}.
\end{remark}

\begin{remark}
Suppose that $H\le1/2$ with $a_0^H=0$; the formula \eqref{bhsimsum} defines a Gaussian process on
the torus $\reel/2\relatif$ with covariance kernel
\begin{equation}\label{covcirc}
   \esp\bigl[B_t^HB_s^H\bigr]=\frac\rho2\bigl(\delta(0,t)^{2H}+\delta(0,s)^{2H}-\delta(s,t)^{2H}\bigr)
\end{equation}
for the distance $\delta$ on the torus. This is the fractional Brownian motion of \cite{istas05}
indexed by the torus. For $H>1/2$, we cannot take $a_0^H=0$; this is related to the fact proved in
\cite{istas05}, that the fractional Brownian motion on the torus does not exist; when indeed such a
process exists, we deduce from \eqref{covcirc} that
\begin{equation*}
   \esp\bigl[B_t^H(B_{1+t}^H-B_1^H)\bigr]
   =\rho\bigl((1-t)^{2H}-1\bigr)\sim-2\rho Ht
\end{equation*}
as $t\downarrow0$ (use the fact $\delta(1+t,0)=1-t$ on the torus), whereas this covariance should
be dominated by $t^{2H}$.
\end{remark}

\begin{remark}\label{statioproc}
When $H\le1/2$ and $a_0^H=0$, we can write $B_t^H$ on $[0,1]$ as $\Abar_t^H-\Abar_0^H$ for the
stationary process $\Abar_t^H=\sum a_n^H(\cos(\pi nt)\xi_n+\sin(\pi nt)\xi_n')$. In the case
$H=1/2$, it generates the same $\sigma$-algebra as $B^{1/2}$, and this process coincides with the
process $A^{1/2}$ of Theorem \ref{statiofrac}. However, a comparison of the variances of the two
processes show that they are generally different when $H<1/2$.
\end{remark}

\begin{remark}
Since the two sides of \eqref{bhsimsum} have stationary increments, we can replace the time
intervals $[0,1]$ and $[0,T]$ of Theorem \ref{bhtrigo} by other intervals of length 1 and $T$
containing 0.
\end{remark}

We now study the asymptotic behaviour of the coefficients $a_n^H$ of Theorem \ref{bhtrigo}.

\begin{theorem}\label{btilth}
The expansion of Theorem \ref{bhtrigo} can be written with $a_0^H=\sqrt{\rho H}$. In this
case, $a_n^H>0$ for any $n$ and
\begin{equation}\label{anhequiv}
   a_n^H=(\pi n)^{-H-1/2}(1+O(n^{2H-3}))
\end{equation}
for $n$ large.
\end{theorem}

\begin{proof}
The only part which has still to be proved is \eqref{anhequiv}. This will be accomplished through
an asymptotic analysis of the integrals in \eqref{bnh1} and \eqref{bnh2}. For $H=1/2$ we have
$a_n^H=(\pi n)^{-1}$ so this is trivial. If $H<1/2$, we have
\begin{align}
   &(1-2H)\int_0^1t^{2H-2}(1-\cos(\pi nt))dt\notag\\
   &=(1-2H)\int_0^\infty t^{2H-2}(1-\cos(\pi nt))dt-1+(1-2H)\int_1^\infty t^{2H-2}\cos(\pi nt)dt\notag\\
   &=(1-2H)(\pi n)^{1-2H}\int_0^\infty t^{2H-2}(1-\cos t)dt-1\notag\\&\qquad\qquad\qquad
   +(1-2H)(\pi n)^{1-2H}\int_{\pi n}^\infty t^{2H-2}\cos t\,dt\notag\\
   &=(\pi n)^{1-2H}\int_0^\infty t^{2H-1}\sin t\,dt-1+O(n^{-2})\label{1m2h}
\end{align}
where we have used in the last equality
\begin{align}
    \Bigl|\int_{\pi n}^\infty t^{2H-2}\cos t\,dt\Bigr|
    &=(2-2H)\Bigl|\int_{\pi n}^\infty t^{2H-3}\sin t\,dt\Bigr|\notag\\
    &=(2-2H)\Bigl|\sum_{k\ge n}\int_{\pi k}^{\pi(k+1)} t^{2H-3}\sin t\,dt\Bigr|\notag\\
    &\le(2-2H)\Bigl|\int_{\pi n}^{\pi(n+1)} t^{2H-3}\sin t\,dt\Bigr|=O(n^{2H-3})\label{2H3}
\end{align}
(this is an alternating series). By applying \eqref{trigo}, we deduce that
\begin{equation*}
    (1-2H)\int_0^1t^{2H-2}(1-\cos(\pi nt))dt=(\pi n)^{1-2H}\Gamma(2H)\sin(\pi H)-1+O(n^{-2}),
\end{equation*}
so \eqref{bnh1} with $a_0^H=\sqrt{\rho H}$ implies
\begin{equation}\label{bnh4}
   b_n^H=\rho(\pi n)^{-1-2H}\Gamma(2H+1)\sin(\pi H)+O(n^{-4}).
\end{equation}
Similarly, if $H>1/2$, then \eqref{2H3} again holds true and
\begin{align*}
   \int_0^1t^{2H-2}\cos(\pi nt)dt
   &=(\pi n)^{1-2H}\int_0^\infty t^{2H-2}\cos t\,dt+O(n^{-2})\\
   &=(\pi n)^{1-2H}\Gamma(2H-1)\sin(\pi H)+O(n^{-2})
\end{align*}
and we deduce from \eqref{bnh2} that we again have \eqref{bnh4}. By using our choice of $\rho$
given in \eqref{rhoh}, we obtain in both cases
\begin{align*}
    b_n^H&=-2\frac{\Gamma(-2H)\Gamma(2H+1)}{\pi^{2H+2}n^{2H+1}}\cos(\pi H)\sin(\pi H)(1+O(n^{2H-3}))\\
    &=(\pi n)^{-2H-1}(1+O(n^{2H-3}))
\end{align*}
from \eqref{beta}. We deduce \eqref{anhequiv} by taking the square root. \qed
\end{proof}

\begin{remark}\label{bch}
Considering the expansion \eqref{bhsimsum} for $a_0^H=\sqrt{\rho H}$, replacing $B^H$ by the
process
\begin{equation*}
   \Bcheck_t^H=c\xi_0t+\sum_{n\ge1}(\pi n)^{-H-1/2}\bigl((\cos(\pi nt)-1)\xi_n+\sin(\pi nt)\xi_n'\bigr)
\end{equation*}
for $c>0$ is equivalent to multiplying $\xi_0$ by $c/a_0^H$ and $(\xi_n,\xi_n')$ by some
$(1+O(n^{2H-3}))$ which remains strictly positive. We can compare the laws of these two sequences
of independent Gaussian variables by means of Kakutani's criterion (Theorem \ref{kakth}), and it
appears that the laws of these two sequences are equivalent ($\sum n^{4H-6}<\infty$). Thus the laws
of $B^H$ and $\Bcheck^H$ are equivalent on $[0,1]$. This implies that the law of
$2^{-H}\Bcheck_{2t}^H$ is equivalent on $[0,1/2]$ to the law of $B_t^H$; actually, we will prove in
Theorem \ref{perequiv} that these two laws are equivalent on $[0,T]$ for any $T<1$.
\end{remark}

\subsection{Approximate expansions}

We now consider the processes
\begin{equation}\label{bhatt}
 \begin{split}
   \Bhat_t^H&=\xi_0t+\sqrt2\,\sum_{n\ge1}\Bigl(\xi_n
   \frac{\cos(2n\pi t)-1}{(2n\pi)^{H+1/2}}
   +\xi_n'\frac{\sin(2n\pi t)}{(2n\pi)^{H+1/2}}\Bigr),\\
   \Bbar_t^H&=\sqrt2\,\sum_{n\ge0}\Bigl(\xi_n
   \frac{\cos((2n+1)\pi t)-1}{((2n+1)\pi)^{H+1/2}}
   +\xi_n'\frac{\sin((2n+1)\pi t)}{((2n+1)\pi)^{H+1/2}}\Bigr)
 \end{split}
\end{equation}
on $[0,1]$. Notice that $\Bhat^{1/2}\simeq\Bbar^{1/2}\simeq W$ from \eqref{wtrigo} and
\eqref{wanti}. On the other hand, it follows from Theorem \ref{bhtrigo} that $\Bhat^H\not\simeq
B^H$ and $\Bbar^H\not\simeq B^H$ for $H\ne1/2$ (because one should have $a_n^H\ne0$ in the
expansion \eqref{bhsimsum} of $B^H$ for all large enough $n$), but we are going to check that these
two processes have a local behaviour similar to $B^H$. The advantage with respect to the exact
expansion \eqref{bhsimsum} is that the sequence of random coefficients and the process will
generate the same $\sigma$-algebra. Then we will apply these approximations to some properties of
the Cameron-Martin space $\hilbert_H$ (Subsection \ref{approxcm}), and to some equivalence of laws
(Subsection \ref{approxlaw}).  As it was the case for Riemann-Liouville processes, $\Bhat^H$ and
$\Bbar^H$ are not only defined for $0<H<1$, but also for any $H>0$.

Let us compare $\Bhat^H$ and $\Bbar^H$ with $B^H$ for $0<H<1$. We use the operators $\Ihatp^\alpha$
and $\Ibarp^\alpha$ defined in \eqref{ihatdef} and \eqref{iflatdef}. By projecting on the Gaussian
spaces generated by $\xi_n$ and $\xi_n'$ and by applying \eqref{ihattrigo}, we can write
\begin{align}
   \Ihatp^{1/2-H}\Bhat_t^H=\xi_0t+\sqrt2\,\sum_{n\ge1}&\Bigl(\xi_n
   \frac{\cos(2\pi nt+(H-1/2)\pi/2)-\cos((H-1/2)\pi/2)}{2\pi n}\notag\\
   &+\xi_n'\frac{\sin(2\pi nt+(H-1/2)\pi/2)-\sin((H-1/2)\pi/2)}{2\pi n}\Bigr).\label{zth}
\end{align}
The two expressions \eqref{wtrigo} and \eqref{zth} are related to each other by applying a rotation
on the vectors $(\xi_n,\xi_n')$, so $\Ihatp^{1/2-H}\Bhat^H$ and $W$ have the same law. A similar
property holds for $\Ibarp^{1/2-H}\Bbar^H$, and we can therefore write
\begin{equation}\label{bhhzh}
    \Bhat^H\simeq\Ihatp^{H-J}\Bhat^J,\quad\Bbar^H\simeq\Ibarp^{H-J}\Bbar^J,
    \quad\Bhat^{1/2}\simeq\Bbar^{1/2}\simeq W.
\end{equation}
We can give an extension of Theorem \ref{coupl1}.

\begin{theorem}\label{coupl3}
It is possible to realise jointly the processes $B^H$, $X^H$, $\Bbar^H$ and $\Bhat^H$ so that the
differences $B^H-X^H$, $\Bbar^H-B^H$ and $\Bhat^H-B^H$ are $C^\infty$ on $(0,1]$; moreover, the
derivatives of order $k$ of these differences are $O(t^{H-k})$ in $L^2(\Omega)$ as $t\downarrow0$.
\end{theorem}

\begin{proof}
We consider the coupling $B^H=\Itilde_+^{H-1/2}W$, $X^H=I_{0+}^{H-1/2}W$, $\Bbar^H=\Ibarp^{H-1/2}W$
and $\Bhat^H=\Ihatp^{H-1/2}W$ for the same $W$ on $\reel$. The smoothness of $B^H-X^H$ is proved in
Theorem \ref{coupl1}, and the estimation of the derivatives follows by a scaling argument. On the
other hand, let $W_t^1$ be equal to $W_t-W_1t$ on $[0,1]$, extend it to $\reel$ by periodicity, and
define $W_t^2=W_{-t}^1$ for $t\ge0$. Then, with the notation \eqref{itriang},
\begin{align*}
   \Bhat_t^H&=W_1t+I_{0+}^{H-1/2}(W_t-W_1t)+I_\triangle^{H-1/2}W_t^2\\
   &=X_t^H+W_1\bigl(t-\Gamma(H+3/2)^{-1}t^{H+1/2}\bigr)+I_\triangle^{H-1/2}W_t^2
\end{align*}
The smoothness of $\Bhat^H-X^H$ follows; the process $W_t^2$ is dominated in $L^2(\Omega)$ by
$\min(\sqrt t,1)$, so we deduce from \eqref{itriang} that
\begin{equation*}
   \bigl\|D^kI_\triangle^{H-1/2}W_t^2\bigr\|_2\le C\int_0^\infty(t+s)^{H-k-3/2}\sqrt s\,ds
   =C't^{H-k}
\end{equation*}
for $k\ge1$. The study of $\Bbar^H$ is similar; let $W^3$ be the process $W$ on $[0,1]$ extended to
$\reel$ so that the increments are 1-antiperiodic, and let $W_t^4=W_{-t}^3$; then $\Bbar^H$ is
equal to $X^H+I_\triangle^{H-1/2}W^4$; the end of the proof is identical. \qed
\end{proof}

\subsection{Application to the Cameron-Martin space}
\label{approxcm}

Let $\hilberthat_H$ and $\hilbertbar_H$ be the Cameron-Martin spaces of $\Bhat^H$ and $\Bbar^H$ on
the time interval $[0,1]$. It follows from \eqref{bhhzh} that
$\hilberthat_{1/2}=\hilbertbar_{1/2}=\hilbert_{1/2}$, and $\hilberthat_H=\Ihatp^{H-J}\hilberthat_J$
as well as $\hilbertbar_H=\Ibarp^{H-J}\hilbertbar_J$.

\begin{theorem}\label{percamar}
For $0<H<1$, the spaces $\hilberthat_H$, $\hilbertbar_H$ and $\hilbert_H$ are equivalent on
$[0,1]$.
\end{theorem}

\begin{proof}
We compare successively $\hilberthat_H$ and $\hilbertbar_H$ with $\hilbert_H([0,1])$, and use the
properties of this last space described in Remark \ref{restrict}.

\medskip\noindent\emph{Proof of $\hilberthat_H\sim\hilbert_H$.}
We know that $\hilberthat_H=\Ihatp^{H-1/2}\hilbert_{1/2}$, so it is sufficient to establish that
$\Ihatp^{H-1/2}$ is a homeomorphism from $\hilbert_{1/2}([0,1])$ onto $\hilbert_H([0,1])$. To this
end, we are going to prove that $\Ihatp^{H-J}$ is continuous from $\hilbert_J([0,1])$ into
$\hilbert_H([0,1])$ for $0<J,H<1$. Consider a function $h$ of $\hilbert_J([0,1])$, consider
$h_0(t)=h(t)-h(1)t$, and extend it by periodicity. Then $h_0$ is generally not in
$\hilbert_J(\reel)$, but the operator $h\mapsto h_1=h_01_{(-1,1]}$ is continuous from
$\hilbert_J([0,1])$ into $\hilbert_J(\reel)$. Moreover, the operator $h\mapsto
h_2=h_01_{(-\infty,-1]}$ is continuous from $\hilbert_J([0,1])$ into the space
$L^\infty((-\infty,-1])$ of bounded functions supported by $(-\infty,-1]$. On the other hand, it is
known that $\hilbert_H=\Itilde_+^{H-J}\hilbert_J$ on $\reel$, and $\Itilde_+^{H-J}$ also maps
continuously $L^\infty((-\infty,-1])$ into the space of smooth functions on $[0,1]$, and therefore
into $\hilbert_H([0,1])$. Thus $h\mapsto\Itilde_+^{H-J}h_0=\Itilde_+^{H-J}h_1+\Itilde_+^{H-J}h_2$
is continuous from $\hilbert_J([0,1])$ into $\hilbert_H([0,1])$. If we add the operator
$h\mapsto(h(1)t)$ which is also continuous, we can conclude.

\medskip\noindent\emph{Proof of $\hilbertbar_H\sim\hilbert_H$.}
In this case, we let $h_0$ be the function $h$ on $[0,1]$, extended to $\reel$ so that the
increments are 1-antiperiodic. We then consider $h_1=h_01_{(-2,1]}$ and $h_2=h_01_{(-\infty,-2]}$.
The proof is then similar, except that we do not have the term $h(1)t$ in this case. \qed
\end{proof}

\begin{remark}
In view of \eqref{spectral}, a function $h$ is in the space $\hilbert_H(\reel)$ if its derivative
$D^1h$ (in distribution sense if $H<1/2$) is in the homogeneous Sobolev space of order $H-1/2$ (see
for instance \cite{norros:sak09}); similarly, it follows from \eqref{bhatt} that $h$ is in
$\hilberthat_H$ is $D^1h$ is in the Sobolev space of order $H-1/2$ of the torus $\reel/\relatif$.
Thus the equivalence $\hilberthat_H\sim\hilbert_H$ of Theorem \ref{percamar} means that the Sobolev
space on the torus is equivalent to the restriction to $[0,1]$ of the Sobolev space on $\reel$.
This classical result is true because we deal with Sobolev spaces of order in $(-1/2,1/2)$.
\end{remark}

\begin{remark}
We have from Theorems \ref{rlcamar} and \ref{percamar} that
$\hilbert_H\sim\hilbert_H'\sim\hilberthat_H\sim\hilbertbar_H$ for any $0<H<1$. Notice however that
the comparison for instance of $\hilberthat_H$ and $\hilbert_H'$ cannot be extended to the case
$H>1$; in this case indeed, functions of $\hilbert_H'$ satisfy $D^1h(0)=0$, contrary to functions
of $\hilberthat_H$.
\end{remark}

Let us now give an immediate corollary of Theorem \ref{percamar}.

\begin{theorem}\label{riesz}
The sets of functions on $[0,1]$
\begin{equation*}
   t,\quad n^{-H-1/2}\bigl(1-\cos(2n\pi t)\bigr),\quad n^{-H-1/2}\sin(2n\pi t),
\end{equation*}
and
\begin{equation*}
   n^{-H-1/2}\bigl(1-\cos((2n+1)\pi t)\bigr),\quad n^{-H-1/2}\sin((2n+1)\pi t),
\end{equation*}
form two Riesz bases of $\hilbert_H$. A function $h$ is in $\hilbert_H$ is and only if it has the
Fourier expansion
\begin{equation*}
   h(t)-h(1)t=\sum_{n\ge0}\alpha_n\cos(2\pi nt)+\sum_{n\ge1}\beta_n\sin(2\pi nt)
\end{equation*}
with
\begin{equation*}
   \sum n^{2H+1}\bigl(\alpha_n^2+\beta_n^2\bigr)<\infty.
\end{equation*}
\end{theorem}

\subsection{Equivalence and mutual singularity of laws}
\label{approxlaw}

We now compare the laws of $B^H$, $\Bhat^H$ and $\Bbar^H$ viewed as variables with values
in the space of continuous functions.

\begin{theorem}\label{perequiv}
Let $H\ne1/2$. The laws of the processes $\Bhat^H$, $\Bbar^H$ and $B^H$ are equivalent on
the time interval $[0,T]$ if $T<1$, and are mutually singular if $T=1$.
\end{theorem}

\begin{proof}
We compare the laws of $B^H$ and $\Bhat^H$. The study of $\Bbar^H$ is similar.

\medskip\noindent\emph{Proof of the equivalence for $0<T<1$.} The increments of both processes are
stationary, so let us study the equivalence of $\Bhat_t^{S,H}=\Bhat_{S+t}^H-\Bhat_S^H$ and
$B_t^{S,H}=B_{S+t}^H-B_S^H$ on $[0,T]$ for $S=1-T$. From Theorem \ref{coupl3}, we can couple $B^H$
and $\Bhat^H$ so that the difference is smooth on $\reel_+^\star$. Consequently,
$\Bhat^{S,H}-B^{S,H}$ is smooth on $[0,T]$, so it lives in $\hilbert_H$. Moreover, we have proved
in Theorem \ref{percamar} that the Cameron-Martin spaces of $\Bhat^H$ and $B^H$ are equivalent, so
the same is true for the Cameron-Martin spaces of $\Bhat^{S,H}$ and $B^{S,H}$. The equivalence of
laws then follows from Theorem \ref{feld}.

\medskip\noindent\emph{Proof of the mutual singularity for $T=1$.}
Consider $\Bhat^H$ on $\reel$. Our aim is to prove that the laws of the two processes
\begin{equation*}
   (B_t^H,B_1^H-B_{1-t}^H)\quad\text{and}
   \quad(\Bhat_t^H,\Bhat_1^H-\Bhat_{1-t}^H)=(\Bhat_t^H,-\Bhat_{-t}^H)\simeq(\Bhat_{2t}^H-\Bhat_t^H,\Bhat_t^H)
\end{equation*}
are mutually singular on the time interval $[0,1/4]$. The law of the first process is equivalent to
a couple $(B_t^{H,1},B_t^{H,2})$ of two independent fractional Brownian motions (see Theorem
\ref{bhindep}), and $\tribuf_{0+}(B^{H,1},B^{H,2})$ is almost surely trivial. On the other hand,
from the first part of this proof, the law of the second process is equivalent to the law of
$(B_{2t}^H-B_t^H,B_t^H)$. We therefore obtain two self-similar processes which do not have the same
law, so we deduce from Theorem \ref{dichot} that the laws are mutually singular. \qed
\end{proof}

\begin{remark}
It follows from Remark \ref{bch} that the law of $B^H$ is equivalent on $[0,1]$ to the law of
$(\Bhat^H+\Bbar^H)/\sqrt2$, where $\Bhat^H$ and $\Bbar^H$ are independent. We have now proved that
this law is equivalent separately to the laws of $\Bhat^H$ and $\Bbar^H$, but only on $[0,T]$ for
$T<1$.
\end{remark}

\begin{theorem}\label{bhsmall}
Let $T>0$. The distance in total variation between the laws of the processes $(\eps^{-H}\Bhat_{\eps
t}^H;\;0\le t\le T)$ and $(B_t^H;\;0\le t\le T)$ is $O(\eps^{1-H})$ as $\eps\downarrow0$. The
process $\Bbar^H$ satisfies the same property.
\end{theorem}

\begin{proof}
As in Theorem \ref{perequiv}, let us compare the laws of $\Bhat^{1/2,H}$ and $B^{1/2,H}$ on
$[0,\eps T]$ for $0<\eps\le1/(2T)$. It follows from Theorem \ref{feld} that the entropy $\inform$
of the former process relative to the latter one satisfies
\begin{equation*}
   \inform\le C\,\esp\bigl|\Bhat^{1/2,H}-B^{1/2,H}\bigr|_{\hilbert_H([0,\eps T])}^2.
\end{equation*}
More precisely it is stated in Theorem \ref{feld} that the constant $C$ involved in this domination
property depends only on the constants involved in the injections of the Cameron-Martin spaces of
$\Bhat^{1/2,H}$ and $B^{1/2,H}$ on $[0,\eps T]$ into each other; but if we choose a constant which
is valid for $\Bhat^H$ and $B^H$ the time interval $[0,1]$ (Theorem \ref{percamar}), then it is
also valid for $\Bhat^{1/2,H}$ and $B^{1/2,H}$ on $[0,1/2]$, and therefore on the subintervals
$[0,\eps T]$, $0<\eps\le1/(2T)$, so we can choose $C$ not depending on $\eps$. Thus
\begin{equation*}
   \inform\le C\,\esp\bigl|\Bhat^{1/2,H}-B^{1/2,H}\bigr|_{\hilbert_H'([0,\eps T])}^2=O(\eps^{2-2H})
\end{equation*}
from \eqref{hpsm}. The convergence in total variation and the speed of convergence are deduced from
\eqref{pinsker}. The proof for $\Bbar^H$ is similar. \qed
\end{proof}

\begin{remark}
We can say that the processes $\Bbar^H$ and $\Bhat^H$ are asymptotically fractional Brownian
motions near time 0. The processes $\Bbar^H$, $\Bhat^H$ and $B^H$ have stationary increments, so
the same local property holds at any time.
\end{remark}

As an application, we recover a result of \cite{cheridito01}, see also
\cite{baudoin:nualart03,vanzanten07} for more general results. Notice that the equivalence stated
in the following theorem may hold even when the paths of $B_2^H$ are not in $\hilbert_J$.

\begin{theorem}\label{cherid}
Let $B_1^J$ and $B_2^H$ be two independent fractional Brownian motions with indices $J<H$, and let
$T>0$. Then the laws of $(B_1^J+\lambda\,B_2^H;\;\lambda\ge0)$ are pairwise equivalent on $[0,T]$
if $H>J+1/4$. Otherwise, they are pairwise mutually singular.
\end{theorem}

\begin{proof}
It is sufficient to prove the result for $T=1$.

\medskip\noindent\emph{Equivalence for $H-J>1/4$.} Let us prove that the laws of $B_1^J$ and
$B_1^J+\lambda B_2^H$, are equivalent. From Theorems \ref{bhtrigo} and \ref{btilth}, the process
$B_1^J$ can be written as \eqref{bhsimsum} for independent standard Gaussian variables
$(\xi_n,\xi_n')$ and coefficients $a_n^J$ such that $a_n^J\ne0$ for any $n$. The process $B_2^H$
can be written similarly with coefficients $a_n^H$ and variables $(\eta_n,\eta_n')$. Thus
$B_1^J+\lambda B_2^H$ is the image by some functional of the sequence
\begin{equation*}
   U_n^\lambda=a_n^J(\xi_n,\xi_n')+\lambda\,a_n^H(\eta_n,\eta_n'),
\end{equation*}
and it is sufficient to prove that the laws of $U_n^\lambda$ and $U_n^0$ are equivalent. This can
be done by means of Kakutani's criterion (Theorem \ref{kakth}) with $\sigma_n^2=(a_n^J)^2$ and
$\bar\sigma_n^2=(a_n^J)^2+\lambda^2(a_n^H)^2$. But
\begin{equation*}
   \sum_{n\ge1}\Bigl(\frac{\lambda^2(a_n^H)^2}{(a_n^J)^2}\Bigr)^2\le C\sum_{n\ge1}n^{4(J-H)}<\infty
\end{equation*}
from Theorem \ref{btilth}.

\medskip\noindent\emph{Mutual singularity for $0<H-J\le1/4$.} Let us use the coupling
\begin{equation*}
   B_1^J=G_{0+}^{1/2,J}W_1,\quad B_2^H=\Itilde_+^{H-1/2}W_2,\quad X_2^K=I_{0+}^{K-1/2}W_2,\quad
   \Bhat_2^K=\Ihatp^{K-1/2}W_2
\end{equation*}
($0<K<1$), for independent $W_1$ on $\reel_+$ and $W_2$ on $\reel$. By applying the
operator $G_{0+}^{J,1/2}$, we can write
\begin{align}
   &G_{0+}^{J,1/2}\bigl(B_1^J+\lambda\,B_2^H\bigr)\label{g0j12}\\
   &=W_1+\lambda\,G_{0+}^{J,1/2}B_2^H\notag\\
   &=W_1+\lambda\Bigl((G_{0+}^{J,1/2}-I_{0+}^{1/2-J})B_2^H+I_{0+}^{1/2-J}(B_2^H-X_2^H)
   +X_2^{1/2+H-J}-\Bhat_2^{1/2+H-J}\Bigr)\notag\\&\qquad+\lambda\,\Bhat_2^{1/2+H-J}.\notag
\end{align}
Let us now prove that the process inside the big parentheses lives in $\hilbert_{1/2}$. We have
checked in the proof of Theorem \ref{coupl2} that $(G_{0+}^{J,1/2}-I_{0+}^{1/2-J})f$ is
differentiable on $\reel_+^\star$ for any $f$ in $\holder^{J-}$, so in particular for $f=B_2^H$;
the scaling property then enables to prove that the derivative is $O(t^{H-J-1/2})$, so
$(G_{0+}^{J,1/2}-I_{0+}^{1/2-J})B_2^H$ is in $\hilbert_{1/2}$. Similarly, $B_2^H-X_2^H$ is smooth,
so $I_{0+}^{1/2-J}(B_2^H-X_2^H)$ is also smooth, and we deduce from the same scaling property that
it is in $\hilbert_{1/2}$. Finally $X_2^{1/2+H-J}-\Bhat_2^{1/2+H-J}$ is also in $\hilbert_{1/2}$
from Theorem \ref{coupl3}. Thus we deduce that the process of \eqref{g0j12} is obtained from
$W_1+\lambda\,\Bhat_2^{1/2+H-J}$ by means of a perturbation which lives in $\hilbert_{1/2}$ and is
independent of $W_1$, so the two laws are equivalent. It is then sufficient to prove that the laws
of $W_1+\lambda_i\,\Bhat_2^{1/2+H-J}$ for $\lambda_1\ne\lambda_2$ are mutually singular. But these
two processes can be expanded on the basis $(t,1-\cos(2\pi nt),\sin(2\pi nt))$; the coefficients
are independent with positive variance; the variance of the coefficients on $1-\cos(2\pi nt)$ and
$\sin(2\pi nt)$ is equal to $2(2\pi n)^{-2}+2\lambda_i^2(2\pi n)^{-2(H-J+1)}$. As in the first
step, we can apply Kakutani's criterion (Theorem \ref{kakth}) and notice that
\begin{equation*}
   \sum_{n\ge1}\Bigl(\frac{(\lambda_2^2-\lambda_1^2)(2\pi n)^{-2(H-J+1)}}{(2\pi n)^{-2}+\lambda_1^2(2\pi n)^{-2(H-J+1)}}
   \Bigr)^2=\infty
\end{equation*}
so that the two laws are mutually singular. \qed
\end{proof}

\begin{remark}
For $H>J$ and $\lambda>0$, the process $B^J+\lambda\,B^H$ exhibits different scaling properties in
finite and large time. It is locally asymptotically $J$-self-similar, whereas it is asymptotically
$H$-self-similar in large time.
\end{remark}

Another application is the comparison with $B^H$ of a fractional analogue of the Karhunen-Lo\`{e}ve
process \eqref{karh} proposed in \cite{feyel:prad99}.

\begin{theorem}
Consider the process
\begin{equation*}
   L_t^H=\sqrt2\,\sum_{n\ge0}\xi_n\frac{\sin\bigl((n+1/2)\pi t\bigr)}{((n+1/2)\pi)^{H+1/2}}
\end{equation*}
for independent standard Gaussian variables $\xi_n$. Then the laws of $L_{S+t}^H-L_S^H$ and $B^H$
are equivalent on $[0,T-S]$ for $0<S<T<1$. On the other hand, these laws are mutually singular if
$S=0$ or $T=1$.
\end{theorem}

\begin{proof}
We deduce from Theorem \ref{perequiv} that the laws of $B_{t/2}^H$ and $\Bbar_{t/2}^H$ are
equivalent on $[0,2T]$ for $T<1$, and therefore on $[-T,T]$ (the two processes have stationary
increments). Thus $(B_t^H-B_{-t}^H)/\sqrt2$, which has the same law as
$2^{H-1/2}(B_{t/2}^H-B_{-t/2}^H)$, has a law equivalent on $[0,T]$ to the law of
\begin{equation*}
   2^{H-1/2}\bigl(\Bbar_{t/2}^H-\Bbar_{-t/2}^H\bigr)
   =2^{H+1}\sum_{n\ge0}\xi_n\frac{\sin((n+1/2)\pi t)}{((2n+1)\pi)^{H+1/2}}=L_t^H,
\end{equation*}
so we have the equivalence of laws
\begin{equation}\label{lthsim}
   L_t^H\sim(B_t^H-B_{-t}^H)/\sqrt2
\end{equation}
on $[0,T]$. Moreover, we deduce from Remark \ref{odde} that the increments of the right hand side
of \eqref{lthsim} on $[S,T]$ are equivalent to the increments of $B^H$, and this proves the first
statement of the theorem. For the case $S=0$, we have also noticed in Remark \ref{odde} that the
laws of the right hand side of \eqref{lthsim} and of $B^H$ are mutually singular. For the case
$T=1$, we have to check that the laws of $L_1^H-L_{1-t}^H$ and of $B^H$ are mutually singular on
$[0,1-S]$. We have
\begin{align*}
   L_1^H-L_{1-t}^H&=2^{H-1/2}\bigl(\Bbar_{1/2}^H-\Bbar_{(1-t)/2}^H-\Bbar_{-1/2}^H+\Bbar_{(t-1)/2}^H\bigr)\\
   &=2^{H-1/2}\bigl(2\Bbar_{1/2}^H-\Bbar_{(1-t)/2}^H-\Bbar_{(1+t)/2}^H\bigr)\\
   &\simeq2^{H-1/2}\bigl(\Bbar_{-t/2}^H+\Bbar_{t/2}^H\bigr)\sim(B_t^H+B_{-t}^H)/\sqrt2
\end{align*}
where we have used the fact that the increments of $\Bbar^H$ are 1-antiperiodic and stationary. But
the law of this process is mutually singular with the law of $B^H$ by again applying Remark
\ref{odde}. \qed
\end{proof}

\setcounter{section}{0}
\renewcommand\thesection{\Alph{section}}

\section*{Appendix}

We now explain some technical results which were used throughout this article.

\section{An analytical lemma}
\label{analytical}

The basic result of this appendix is the following classical lemma, see Theorem 1.5 of
\cite{samko:kil:mari93}.

\begin{theorem}\label{analth}
Consider a kernel $K(t,s)$ on $\reel_+\times\reel_+$ such that
\begin{equation}\label{scalk}
   K(\lambda t,\lambda s)=K(t,s)/\lambda
\end{equation}
for $\lambda>0$, and
\begin{equation*}
   \int_0^\infty\frac{|K(1,s)|}{\sqrt s}ds<\infty.
\end{equation*}
Then $K:f\mapsto\int K(.,s)f(s)ds$ defines a continuous endomorphism of $L^2$.
\end{theorem}

\begin{proof}
For $f$ nonnegative, let us study
\begin{align*}
   E(f)&=\int_0^\infty\Bigl(\int_0^\infty|K(t,s)|
   f(s)ds\Bigr)^2dt
   =\int_0^\infty\Bigl(\int_0^\infty|K(1,s)|
   f(ts)ds\Bigr)^2dt\\
   &=\iiint|K(1,s)|\,|K(1,u)|f(ts)f(tu)ds\,du\,dt
\end{align*}
from the scaling property \eqref{scalk} written as $K(t,s)=K(1,s/t)/t$. We have
\begin{equation*}
  \int f(ts)f(tu)dt\le\|f\|_{L^2}^2/\sqrt{su},
\end{equation*}
so
\begin{equation*}
  E(f)\le\|f\|_{L^2}^2\Bigl(\int\frac{|K(1,s)|}{\sqrt s}ds\Bigr)^2.
\end{equation*}
If now $f$ is a real square integrable function, then $Kf(t)$ is well defined for almost any $t$,
and
\begin{equation*}
   \int_0^\infty Kf(t)^2dt\le E(|f|)\le C\|f\|_{L^2}^2.
\end{equation*}
\qed
\end{proof}

\begin{theorem}\label{ahlambda}
On the time interval $\reel_+$, let
\begin{equation*}
   A:(h(t);\;t\ge0)\mapsto(Ah(t);\;t\ge0)
\end{equation*}
be a linear operator defined on $\holder^{1/2}$ (the space of $1/2$-H\"{o}lder continuous functions
taking the value 0 at 0) such that $Ah(0)=0$. We suppose that
\begin{equation}\label{scalc}
   A(h_\lambda)=(Ah)_\lambda\quad\text{for $h_\lambda(t)=h(\lambda t)$.}
\end{equation}
We also suppose that $Ah$ is differentiable on $\reel_+^\star$ and that $h\mapsto D^1Ah(1)$ is
continuous on $\holder^{1/2}$. Then $A$ is a continuous endomorphism of the standard Cameron-Martin
space $\hilbert_{1/2}=I_{0+}^1L^2$.
\end{theorem}

\begin{proof}
On $\hilbert_{1/2}$, the linear form $h\mapsto D^1Ah(1)$ takes the form $D^1Ah(1)=\langle
a,h\rangle_{\hilbert_{1/2}}$ for some $a$ in $\hilbert_{1/2}$, so
\begin{align*}
   D^1Ah(t)&=\frac1tD^1(Ah)_t(1)=\frac1tD^1Ah_t(1)=\frac1t\langle a,h_t\rangle_{\hilbert_{1/2}}
   =\frac1t\int D^1a(s)\,D^1h_t(s)ds\\&=\int D^1a(s)\,D^1h(ts)ds=\int K(t,s)D^1h(s)ds
\end{align*}
for
\begin{equation*}
   K(t,s)=D^1a(s/t)/t.
\end{equation*}
Then $K$ satisfies the scaling condition \eqref{scalk}, and
\begin{align*}
   \int\frac{|D^1a(s)|}{\sqrt s}ds&\le\sup\Bigl\{\langle a,h\rangle_{\hilbert_{1/2}};\;h\in\hilbert_{1/2},\;
   |D^1h(s)|\le1/\sqrt s\Bigr\}\\
   &\le\sup\Bigl\{D^1Ah(1);\;h(0)=0,\;
   |h(t)-h(s)|\le2\sqrt{t-s}\Bigr\}<\infty
\end{align*}
since $h\mapsto D^1Ah(1)$ is continuous on $\holder^{1/2}$. Thus we can apply Theorem \ref{analth}
and deduce that $D^1AI_{0+}^1$ is a continuous endomorphism of $L^2$, or, equivalently, that $A$ is
a continuous endomorphism of $\hilbert_{1/2}$. \qed
\end{proof}

\section{Variance of fractional Brownian motions}
\label{variance}

We prove here a result stated in Subsection \ref{repreel}, more precisely that if $B^H$ is given
by the representation \eqref{mandkappa} with $\kappa$ given by \eqref{kappah}. then the
variance $\rho$ of $B_1^H$ satisfies \eqref{rhoh}. We also prove that the variance of $B_1^H$
given by the spectral representation \eqref{spectral} is the same.

\begin{theorem}\label{varbh}
The variance of $B_1^H$ defined by \eqref{mandkappa} is given by
\begin{equation}\label{rhokappa2}
   \rho=\kappa^2\frac{3/2-H}{2H}\,B(2-2H,H+1/2)
\end{equation}
for the Beta function
\begin{equation*}
   B(\alpha,\beta)=\int_0^1t^{\alpha-1}(1-t)^{\beta-1}dt,\qquad\alpha>0,\>\beta>0.
\end{equation*}
\end{theorem}

\begin{proof}
For $t>0$, by decomposing the right-hand side of \eqref{mandkappa} into integrals on $[0,t]$
and on $\reel_-$, we obtain
\begin{equation*}
   \esp[(B_t^H)^2]=\kappa^2\Bigl(\frac{t^{2H}}{2H}+\phi(t)\Bigr)
\end{equation*}
with
\begin{equation*}
   \phi(t)=\int_0^\infty\Bigl((t+x)^{H-1/2}-x^{H-1/2}
   \Bigr)^2dx.
\end{equation*}
We can differentiate twice this integral and get
\begin{equation*}
   \phi'(t)=(2H-1)\int_0^\infty\Bigl((t+x)^{2H-2}-(t+x)^{H-3/2}x^{H-1/2}\Bigr)dx,
\end{equation*}
\begin{align*}
   \phi''(t)&=(2H-1)(2H-2)\int_0^\infty(t+x)^{2H-3}dx\\
   &\quad-(2H-1)(H-3/2)\int_0^\infty (t+x)^{H-5/2}x^{H-1/2}dx\\
   &=-(2H-1)t^{2H-2}-(2H-1)(H-3/2)t^{2H-2}\int_1^\infty y^{H-5/2}(y-1)^{H-1/2}dy\\
   &=-(2H-1)t^{2H-2}-(2H-1)(H-3/2)t^{2H-2}\int_0^1\Bigl(\frac{1-z}{z^2}\Bigr)^{H-1/2}dz
\end{align*}
by means of the changes of variables $x=t(y-1)$ and $y=1/z$. Thus
\begin{equation*}
   \phi''(t)=(2H-1)t^{2H-2}\Bigl(-1+(3/2-H)\,B(2-2H,H+1/2)\Bigr).
\end{equation*}
We integrate twice this formula, and since $\phi(t)$ and $\phi'(t)$ are respectively proportional
to $t^{2H}$ and $t^{2H-1}$, we obtain \eqref{rhokappa2} by writing
$\kappa^2\bigl(\phi(1)+1/(2H)\bigr)$. \qed
\end{proof}

By applying properties of Beta and Gamma functions
\begin{equation}\label{beta}
\begin{split}
   &B(\alpha,\beta)=\Gamma(\alpha)\Gamma(\beta)/\Gamma(\alpha+\beta),\\
   &\Gamma(z+1)=z\,\Gamma(z),\qquad\Gamma(z)\Gamma(1-z)=\pi/\sin(\pi z),
\end{split}
\end{equation}
where $\Gamma$ is defined on $\complex\setminus\relatif_-$, we can write equivalent forms
which are used in the literature,
\begin{align}
   \rho
   &=\kappa^2\frac{3/2-H}{2H}\,\frac{\Gamma(2-2H)\Gamma(H+1/2)}{\Gamma(5/2-H)}\nonumber\\
   &=\kappa^2\frac{1}{2H(1/2-H)}\,\frac{\Gamma(2-2H)\Gamma(H+1/2)}{\Gamma(1/2-H)}\nonumber\\
   &=\kappa^2\frac{\cos(\pi H)}{\pi H(1-2H)}\Gamma(2-2H)\Gamma(H+1/2)^2\nonumber\\
   &=-2\kappa^2\frac{\cos(\pi H)}{\pi}\Gamma(-2H)\Gamma(H+1/2)^2\label{rhokcos}
\end{align}
where, except in the first line, we have to assume $H\ne1/2$. Thus if we choose
$\kappa=\kappa(H)=\Gamma(H+1/2)^{-1}$ as this is done in this article, then $\rho$ is given by
\eqref{rhoh}.

If now we consider the spectral representation \eqref{spectral}, then
\begin{align*}
    \esp\bigl[(B_1^H)^2\bigr]
    &=\frac1\pi\int_0^\infty s^{-1-2H}\Bigl(\bigl(\cos s-1\bigr)^2+\sin^2s\Bigr)ds\\
    &=\frac2\pi\int_0^\infty s^{-1-2H}\bigl(1-\cos s\bigr)ds
    =\frac{1}{\pi H}\int_0^\infty s^{-2H}\sin s\,ds
\end{align*}
by integration by parts. If $H<1/2$, an application of \eqref{trigo} shows that this variance is
again given by \eqref{rhoh}; if $H>1/2$, the same property can be proved by using another
integration by parts, and the case $H=1/2$ can be deduced from the continuity of the variance
with respect to $H$.

\begin{remark}
The variance of the spectral decomposition can also be obtained as follows. The process $B^H$
given by \eqref{spectral} can be written as the real part of
\begin{align*}
    B_t^{H,\complex}&=\frac1{\sqrt\pi}\int_0^{+\infty}s^{-H-1/2}\bigl(e^{ist}-1\bigr)
    \bigl(dW_s^1+i\,dW_s^2\bigr)\\
    &\simeq\frac1{\sqrt{2\pi}}\int_{-\infty}^{+\infty}|s|^{1/2-H}\frac{e^{ist}-1}{s}
    \bigl(dW_s^1+i\,dW_s^2\bigr).
\end{align*}
The isometry property of the Fourier transform on $L^2$ enables to check that
$B^{1/2,\complex}$ has the same law as $W^1+i\,W^2$, so in particular $B^{1/2}$ is a
standard Brownian motion. Following Theorem \ref{reprth}, the general case $H\ne1/2$ is
obtained by applying $\Itilde_+^{H-1/2}$ to $B^{1/2,\complex}$ (use \eqref{trigo}).
\end{remark}

\section{Equivalence of laws of Gaussian processes}
\label{equiv}

Our aim is to compare the laws of two centred Gaussian processes. It is known from
\cite{feldman58,hajek58,hida:hitsuda93} that their laws are either equivalent, or mutually singular
(actually this is also true in the non centred case), and we want to decide between these two
possibilities. In Subsection \ref{camarsec}, after a brief review of infinite dimensional Gaussian
variables, we explain how the Cameron-Martin space (or reproducing kernel Hilbert space) can be
used to study this question. In particular, we prove a sufficient condition for the equivalence.
Then, in Subsection \ref{covself}, we describe a more computational method which can be used for
self-similar processes to decide between the equivalence and mutual singularity.

\subsection{Cameron-Martin spaces}\label{camarsec}

A Gaussian process can be viewed as a Gaussian variable $W$ taking its values in an
infinite-dimensional vector space $\wiener$, but the choice of $\wiener$ is not unique; in order to
facilitate the study of $W$, it is better for $\wiener$ to have a good topological structure. This
is with this purpose that the notion of abstract Wiener space was introduced by \cite{gross67}; in
this framework, $\wiener$ is a separable Banach space. However, more general topological vector
spaces can also be considered, see for instance \cite{bogachev98}. Here, we assume that $\wiener$
is a separable Fr\'echet space and we let $\wiener^\star$ be its topological dual. The space
$\wiener$ is endowed with its Borel $\sigma$-algebra, which coincides with the cylindrical
$\sigma$-algebra generated by the maps $w\mapsto \ell(w)$, $\ell\in\wiener^\star$. A
$\wiener$-valued variable $W$ is said to be centred Gaussian if $\ell(W)$ is centred Gaussian for
any $\ell\in\wiener^\star$; the closed subspace of $L^2(\Omega)$ generated by the variables
$\ell(W)$ is the Gaussian space of $W$. The Fernique theorem (see Theorem 2.8.5 in
\cite{bogachev98}) states that if $|.|$ is a measurable seminorm on $\wiener$ (which may take
infinite values) and if $|W|$ is almost surely finite, then $\exp(\lambda|W|^2)$ is integrable for
small enough positive $\lambda$.

For $h$ in $\wiener$, define
\begin{equation}\label{hhsup}
   |h|_\hilbert=\sup\Bigl\{\frac{\ell(h)}{\bigl\|\ell(W)\bigr\|_2};\;\ell\in\wiener^\star\Bigr\}
\end{equation}
with the usual convention $0/0=0$. Then $\hilbert=\bigl\{h;\;|h|_\hilbert<\infty\bigr\}$ is a
separable Hilbert space which is continuously embedded in $\wiener$ and which is called the
Cameron-Martin space of $W$; it is dense in $\wiener$ if the topological support of the law of $W$
is $\wiener$. It can be identified to its dual, and the adjoint of the inclusion
$i:\hilbert\to\wiener$ is a map $i^\star:\wiener^\star\to\hilbert$ with dense image such that
\begin{equation}\label{istarell}
   \langle i^\star(\ell),h\rangle_\hilbert=\ell(h),\quad
   \langle i^\star(\ell_1),i^\star(\ell_2)\rangle_\hilbert=\esp\bigl[\ell_1(W)\ell_2(W)\bigr].
\end{equation}
Consequently, the map $\ell\mapsto \ell(W)$ can be extended to an isometry between $\hilbert$ and
the Gaussian space of $W$, that we denote by $\langle W,h\rangle_\hilbert$ (though $W$ does not
live in $\hilbert$); thus $\ell(W)=\langle W,i^\star(\ell)\rangle_\hilbert$ and
\begin{equation}\label{whwh}
   \esp\bigl[\langle W,h\rangle_\hilbert\>\langle W,h'\rangle_\hilbert\bigr]=\langle h,h'\rangle_\hilbert.
\end{equation}
The variable $\langle W,h\rangle_\hilbert$ is called the Wiener integral of $h$.

\begin{example}
When considering real continuous Gaussian processes, the space $\wiener$ can be taken to be the
space of real-valued continuous functions with the topology of uniform convergence on compact
subsets. The most known example is the standard Brownian motion; its Cameron-Martin space
$\hilbert_{1/2}$ is the space of absolutely continuous functions $h$ such that $h(0)=0$ and $D^1h$
is in $L^2$.
\end{example}

\begin{remark}
Let $\wiener$ be the space of real-valued continuous functions. The coordinate maps
$\ell_t(\omega)=\omega(t)$ are in $\wiener^\star$ and the linear subspace generated by the
variables $\ell_t(W)=W_t$ is dense in the Gaussian space of $W$; equivalently, the space $\hilbert$
is generated by the elements $i^\star(\ell_t)$. On the other hand, we deduce from \eqref{istarell}
that
\begin{equation*}
   i^\star(\ell_t):s\mapsto\ell_s\bigl(i^\star(\ell_t)\bigr)
   =\langle i^\star(\ell_s),i^\star(\ell_t)\rangle_\hilbert=\esp[W_sW_t].
\end{equation*}
Thus, if we denote by $C(s,t)=\esp[W_sW_t]$ the covariance kernel, then $\hilbert$ is the closure
of the linear span of the functions $i^\star(\ell_t)=C(t,.)$ for the inner product
\begin{equation*}
   \langle C(s,.),C(t,.)\rangle_\hilbert=C(s,t).
\end{equation*}
This relation is called the reproducing property, and $\hilbert$ is the reproducing kernel Hilbert
space of $C(.,.)$. This technique can also be used for non continuous processes, see for instance
\cite{vanvaart:zanten08}.
\end{remark}

\begin{remark}
Another viewpoint for the Wiener integrals when $W=(W_t)$ is a continuous Gaussian process is to
consider the integrals $\int f(t)dW_t$ for deterministic functions $f$. This integral is easily
defined when $f$ is an elementary (or step) process, and we can extend by continuity this
definition to more general functions. With this method, we obtain variables which are in the
Gaussian space of $W$, but we do not necessarily obtain the whole space, see the case of the
fractional Brownian motion $B^H$ when $H>1/2$ in \cite{pipiras:taqqu01}.
\end{remark}

Let $W_1$ and $W_2$ be two centred Gaussian variables with values in the same space $\wiener$, with
Cameron-Martin spaces $\hilbert_1$ and $\hilbert_2$. It follows from \eqref{hhsup} that
$\hilbert_1$ is continuously embedded in $\hilbert_2$ if and only if
\begin{equation}\label{cmembed}
   \bigl\|\ell(W_1)\bigr\|_2\le C\bigl\|\ell(W_2)\bigr\|_2
\end{equation}
for any $\ell\in\wiener^\star$.

Let $\wiener^1$ and $\wiener^2$ be separable Fr\'echet spaces, let $W$ be a $\wiener^1$-valued
centred Gaussian variable with Cameron-Martin space $\hilbert^1$, and let $A:\wiener^1\to\wiener^2$
be a measurable linear transformation which is defined on a measurable linear subspace of
$\wiener^1$ supporting the law of $W$. Then $AW$ is a centred Gaussian variable. If $A$ is
injective on $\hilbert^1$, then the Cameron-Martin space of $AW$ is $\hilbert^2=A(\hilbert^1)$.
This explains how the Cameron-Martin space $\hilbert_H$ of the fractional Brownian motion $B^H$ can
be deduced from $\hilbert_{1/2}$; one applies the transformations $\Itilde_+^{H-1/2}$ (Theorem
\ref{reprth}) or $G_{0+}^{1/2,H}$ (Theorem \ref{molgolo}). On the other hand, if $A$ is non
injective, one still has $\hilbert^2=A(\hilbert^1)$ and the norm is now given by
\begin{equation}\label{contract}
   |h_2|_{\hilbert^2}=\inf\bigl\{|h_1|_{\hilbert^1};\;A(h_1)=h_2\bigr\}.
\end{equation}
In particular $|Ah|_{\hilbert_2}\le|h|_{\hilbert_1}$. If $A=0$ on $\hilbert_1$, then $AW=0$.

We now consider the absolute continuity of Gaussian measures with respect to one another. This
notion can be studied by means of the relative entropy, or Kullback-Leibler divergence, defined for
probability measures $\mu_1$ and $\mu_2$ by
\begin{equation*}
   \inform(\mu_2,\mu_1)=\int\ln\bigl(d\mu_2/d\mu_1\bigr)d\mu_2
\end{equation*}
if $\mu_2$ is absolutely continuous with respect to $\mu_1$, and by $+\infty$ otherwise. This
quantity is related to the total variation of $\mu_2-\mu_1$ by the Pinsker inequality
\begin{equation}\label{pinsker}
   \Bigl(\int\bigl|d\mu_2-d\mu_1\bigr|\Bigr)^2\le2\inform(\mu_2,\mu_1).
\end{equation}

The Cameron-Martin theorem enables to characterise elements of $\hilbert$ amongst elements of
$\wiener$. More precisely, $h$ is in $\hilbert$ if and only if the law of $W+h$ is absolutely
continuous with respect to the law of $W$. Moreover, in this case, the density is
$\exp\bigl(\langle W,h\rangle_\hilbert-|h|_\hilbert^2/2\bigr)$. Thus
\begin{equation*}
   \inform(\mu',\mu)=\inform(\mu,\mu')=|h|_\hilbert^2/2
\end{equation*}
when $\mu$ and $\mu'$ are the laws of $W$ and $W+h$.

The transformation $W\mapsto W+h$ of the Cameron-Martin space can be generalised to random $h$. If
we add to $W$ an independent process $X$ taking its values in $\hilbert$, it is easily seen by
working conditionally on $X$ that the laws of $W$ and $W+X$ are again equivalent. Moreover, the law
of $(W+X,X)$ is absolutely continuous with respect to the law of $(W,X)$, with a density equal to
$\exp\bigl(\langle W,X\rangle_\hilbert-|X|_\hilbert^2/2\bigr)$, and relative entropies of the two
variables with respect to each other are equal to $\frac12\esp|X|_\hilbert^2$. By projecting on the
first component, it follows from the Jensen inequality that the relative entropy cannot increase,
so
\begin{equation}\label{entrop}
   \max\bigl(\inform(\mu',\mu),\inform(\mu,\mu')\bigr)\le\esp|X|_\hilbert^2/2
\end{equation}
when $\mu$ and $\mu'$ are the laws of $W$ and $W+X$.

When $W=(W_n)$ and $\Wbar=(\Wbar_n)$ are two sequences consisting of independent centred Gaussian
variables with positive variances, then the equivalence or mutual singularity of their laws can be
decided by means of Kakutani's criterion \cite{kakutani48}. This criterion is actually intended to
general non Gaussian variables; when specialised to the Gaussian case, it leads to the following
result.

\begin{theorem}\label{kakth}
Let $W=(W_n)$ and $\Wbar=(\Wbar_n)$ be two sequences of independent centred Gaussian variables with
variances $\sigma_n^2>0$ and $\bar\sigma_n^2>0$. Then the laws of $W$ and $\Wbar$ are equivalent if
and only if
\begin{equation}\label{kaku}
   \sum_n\Bigl(\frac{\bar\sigma_n^2}{\sigma_n^2}-1\Bigr)^2<\infty.
\end{equation}
\end{theorem}

Returning to general Gaussian variables, we now give a sufficient condition for the equivalence of
$W$ and $W+X$ where $W$ and $X$ are not required to be independent. This result has been used in
the proof of Theorem \ref{perequiv}; it can be deduced from the proof of \cite{feldman58}, but we
explain its proof for completeness.

\begin{theorem}\label{feld}
Let $(W,X)$ be a centred Gaussian variable with values in $\wiener\times\hilbert$, where $\wiener$
is a separable Fr\'echet space, and $\hilbert$ is the Cameron-Martin space of $W$; thus $W+X$ is a
Gaussian variable taking its values in $\wiener$; let $\hilbert'$ be its Cameron-Martin space.
\begin{itemize}
\item The space $\hilbert'$ is continuously embedded in $\hilbert$.
\item If moreover $\hilbert$ is continuously embedded in $\hilbert'$ (so that
    $\hilbert\sim\hilbert'$), then the laws of $W$ and $W+X$ are equivalent. Moreover, the
    entropy of the law of $W+X$ relative to the law of $W$ is bounded by
    $C\,\esp|X|_\hilbert^2$, where $C$ depends only on the norms of the injections of
    $\hilbert$ and $\hilbert'$ into each other.
\end{itemize}
\end{theorem}

\begin{proof}
We have to compare the laws of $(\ell(W);\ell\in\wiener^\star)$ and
$(\ell(W+X),\ell\in\wiener^\star)$. Since $|X|_\hilbert$ is almost surely finite, it follows from
the Fernique theorem that $|X|_\hilbert^2$ has an exponential moment and is in particular
integrable, so $\ell(X)=\langle i^\star(\ell),X\rangle_\hilbert$ is square integrable. Thus
\begin{equation*}
   \bigl\|\ell(W+X)\bigr\|_2\le\bigl\|\ell(W)\bigr\|_2+C\bigl|i^\star(\ell)\bigr|_\hilbert
   \le(C+1)\bigl\|\ell(W)\bigr\|_2
\end{equation*}
and the inclusion $\hilbert'\subset\hilbert$ follows from \eqref{cmembed}. Let us now suppose
$\hilbert\sim\hilbert'$, so that, by again applying \eqref{cmembed},
\begin{equation}\label{ellwx}
  C_1\bigl\|\ell(W)\bigr\|_2\le\bigl\|\ell(W+X)\bigr\|_2\le C_2\bigl\|\ell(W)\bigr\|_2
\end{equation}
for positive $C_1$ and $C_2$. Let us first compare the laws of the families
$(\ell(W+X);\;\ell\in\wiener_1^\star)$ and $(\ell(W);\;\ell\in\wiener_1^\star)$ for a
finite-dimensional subspace $\wiener_1^\star$ of $\wiener^\star$. We have
\begin{equation*}
  \wiener_0^\star=\bigl\{\ell\in\wiener^\star;\;\bigl\|\ell(W)\bigr\|_2=0\bigr\}
  =\bigl\{\ell\in\wiener^\star;\;\bigl\|\ell(W+X)\bigr\|_2=0\bigr\}
\end{equation*}
and it is sufficient to consider the case where $\wiener_1^\star\cap\wiener_0^\star=\{0\}$. Then
$|\ell|=\|\ell(W)\bigr\|_2$ and $|\ell|'=\bigl\|\ell(W+X)\bigr\|_2$ define two Euclidean structures
on $\wiener_1^\star$, and it is possible to find a basis $(\ell_n;\;1\le n\le N)$ which is
orthonormal for the former norm, and orthogonal for the latter norm. We have to compare the laws
$\mu_N$ and $\mu_N'$ of $U_N=(\ell_n(W);1\le n\le N)$ and $U_N'=(\ell_n(W+X);1\le n\le N)$. The
vectors $U_N$ and $U_N'$ consist of independent centred Gaussian variables; moreover, $U_n$ has
variance 1, and it follows from \eqref{ellwx} that $U'_n$ has a variance $\sigma_n^2$ satisfying
$C_1\le\sigma_n^2\le C_2$. We deduce that
\begin{equation*}
   \inform(\mu_N',\mu_N)=\frac12\sum_{n=1}^N\bigl(\sigma_n^2-1-\ln\sigma_n^2\bigr)
   \le C\sum_{n=1}^N(\sigma_n^2-1)^2.
\end{equation*}
But
\begin{equation}\label{sigmn}
   \sigma_n^2-1=2\,\esp\bigl[\ell_n(W)\,\ell_n(X)\bigr]+\esp\bigl[(\ell_n(X))^2\bigr]
   \le C\Bigl(\esp\bigl[(\ell_n(X))^2\bigr]\Bigr)^{1/2}
\end{equation}
(we deduce from $\sigma_n^2\le C_2$ that the variances of $\ell_n(X)$ are uniformly bounded), and
\begin{equation*}
   \inform(\mu_N',\mu_N)\le C\sum_{n=1}^N\esp\Bigl[(\ell_n(X))^2\Bigr]
   =C\sum_{n=1}^N\esp\Bigl[\langle i^\star(\ell_n),X\rangle_\hilbert^2\Bigr]
   \le C\,\esp|X|_\hilbert^2
\end{equation*}
because $i^\star(\ell_n)$ is from \eqref{istarell} an orthonormal sequence in $\hilbert$. Thus the
entropy of the law of $(\ell(W+X);\;\ell\in\wiener_1^\star)$ relative to
$(\ell(W);\;\ell\in\wiener_1^\star)$ is bounded by an expression $C\,\esp|X|_\hilbert^2$ which does
not depend on the choice of the finite-dimensional subspace $\wiener_1^\star$. This implies that
the law in $\wiener$ of $W+X$ is absolutely continuous with respect to the law of $W$, and that the
corresponding relative entropy is also bounded by this expression. \qed
\end{proof}

\begin{remark}
The condition about the equivalence of Cameron-Martin spaces cannot be dropped in Theorem
\ref{feld}, see the counterexample of the Brownian motion $W=(W_t)$ and $X_t=-t\,W_1$.
\end{remark}

\begin{remark}
If $W$ and $X$ are independent, then
\begin{equation*}
   \bigl\|\ell(W+X)\bigr\|_2^2=\bigl\|\ell(W)\bigr\|_2^2+\bigl\|\ell(X)\bigr\|_2^2
   \ge\bigl\|\ell(W)\bigr\|_2^2
\end{equation*}
so $\hilbert\subset\hilbert'$ is automatically satisfied. Moreover the estimation \eqref{sigmn} is
improved and we have $\esp\langle X,h_n\rangle_\hilbert^2$ instead of its square root. This
explains why the laws of $W$ and $W+X$ can be equivalent even when $X$ does not take its values in
$\hilbert$; when $W$ and $X$ consist of sequences of independent variables (and assuming again that
$\hilbert\sim\hilbert'$), this improvement leads to the condition \eqref{kaku}.
\end{remark}

\begin{remark}
More generally, for the comparison of two centred Gaussian measures $\mu$ and $\mu'$ on a separable
Fr\'echet space $\wiener$, a necessary condition for the equivalence of $\mu$ and $\mu'$ is the
equivalence of the Cameron-Martin spaces $\hilbert$ and $\hilbert'$. If this condition holds, there
exists a homeomorphism $Q$ of $\hilbert$ onto itself such that
\begin{equation*}
   \langle h_1,h_2\rangle_{\hilbert'}=\langle h_1,Qh_2\rangle_\hilbert.
\end{equation*}
Then $\mu$ and $\mu'$ are equivalent if and only if $Q-I$ is a Hilbert-Schmidt operator.
\end{remark}

\subsection{Covariance of self-similar processes}\label{covself}

Consider a square integrable $H$-self-similar process for $H>0$; we now explain that if it
satisfies a 0-1 law in small time, then its covariance kernel can be estimated by means of its
behaviour in small time; this is a simple consequence of the Birkhoff ergodic theorem.

\begin{theorem}\label{covgauss}
Let $(\Xi_t;\;t>0)$ be a $H$-self-similar continuous process, and suppose that its filtration
$\tribuf_t(\Xi)$ is such that $\tribuf_{0+}(\Xi)$ is almost surely trivial. Define
\begin{equation*}
   \theta_r\Xi(t)=e^{Hr}\Xi(e^{-r}t),\qquad-\infty<r<+\infty.
\end{equation*}
Then for any measurable functional $f$ on the space of continuous paths such that $f(\Xi)$ is
integrable,
\begin{equation}\label{asconv}
   \lim_{T\to\infty}\frac1T\int_0^Tf(\theta_r\Xi)dr=\esp[f(\Xi)]
\end{equation}
almost surely. In particular, if $\Xi=(\Xi^1,\ldots,\Xi^n)$ is square integrable,
\begin{equation}\label{xiiuxijv}
   \esp[\Xi_u^i\Xi_v^j]=\lim_{t\to0}
   \frac1{|\log t|}\int_t^1\frac{\Xi_{us}^i\Xi_{vs}^j}{s^{2H+1}}ds.
\end{equation}
\end{theorem}

\begin{proof}
One has $\theta_r\theta_{r'}=\theta_{r+r'}$, so $(\theta_r)$ is a family of shifts. Moreover, the
$H$-self-similarity of the process $\Xi$ is equivalent to the shift invariance of its law. Events
which are $(\theta_r)$-invariant are in $\tribuf_{0+}(\Xi)$ which is almost surely trivial, so the
ergodic theorem enables to deduce \eqref{asconv}. Then \eqref{xiiuxijv} is obtained by taking
$f(\Xi)=\Xi_u^i\Xi_v^j$ and by applying the change of variable $r=\log(1/s)$ in the integral. \qed
\end{proof}

\begin{remark}
By using the Lamperti transform defined in \eqref{lamper}, the family $(\theta_r)$ is reduced to
the time translation on stationary processes.
\end{remark}

\begin{remark}\label{blumen}
In the centred Gaussian case, the law is characterised by the covariance kernel, so Theorem
\ref{covgauss} implies that the whole law of $\Xi$ can be deduced from its small time behaviour.
The result can be applied to fractional Brownian motions of index $0<H<1$; by applying the
canonical representation of Section~\ref{represent}, one has indeed
$\tribuf_{0+}(B^H)=\tribuf_{0+}(W)$ and this $\sigma$-algebra is well-known to be almost surely
trivial (Blumenthal 0-1 law). A simple counterexample is the fractional Brownian motion of index
$H=1$; this process (which was always excluded from our study of $B^H$) is given by $B_t^1=t\,B_1$
for a Gaussian variable $B_1$; the assumption about $\tribuf_{0+}(\Xi)$ and the conclusion of the
theorem do not hold.
\end{remark}

\begin{remark}
In the Gaussian case, \eqref{xiiuxijv} is a simple way to prove that the law of $\Xi$ can be
deduced from its small time behaviour. There are however other techniques, such as Corollary 3.1 of
\cite{arcones95} about the law of iterated logarithm.
\end{remark}

\begin{theorem}\label{dichot}
Let $\Xi$ and $\Upsilon$ be two centred continuous $H$-self-similar Gaussian processes on $[0,1]$,
such that $\tribuf_{0+}(\Xi)$ is almost surely trivial. Then the two processes either have the same
law, or have mutually singular laws.
\end{theorem}

\begin{proof}
Gaussian measures are either equivalent, or mutually singular, so suppose that the laws of $\Xi$
and $\Upsilon$ are equivalent. The process $\Xi$ satisfies \eqref{xiiuxijv}, so
\begin{equation*}
   \esp[\Xi_u^i\Xi_v^j]=\lim_{t\to0}
   \frac1{|\log t|}\int_t^1\frac{\Upsilon_{us}^i\Upsilon_{vs}^j}{s^{2H+1}}ds.
\end{equation*}
Moreover, the right hand side is bounded in $L^p(\Omega)$ for any $p$, so we can take the
expectation in the limit, and it follows from the self-similarity of $\Upsilon$ that
\begin{equation*}
   \esp[\Xi_u^i\Xi_v^j]=\lim_{t\to0}
   \frac1{|\log t|}\int_t^1\frac{\esp[\Upsilon_{us}^i\Upsilon_{vs}^j]}{s^{2H+1}}ds
   =\esp[\Upsilon_u^i\Upsilon_v^j].
\end{equation*}
Thus $\Xi$ and $\Upsilon$ have the same law. \qed
\end{proof}

A counterexample of this property is again the fractional Brownian motion with index $H=1$.
Processes corresponding to different variances $\rho=\esp[(B_1)^2]>0$ have equivalent but
different laws.


\end{document}